\documentclass[%
 aip,
 amsmath,amssymb,
 reprint,%
]{revtex4-2}

\usepackage{graphicx}
\usepackage{dcolumn}
\usepackage{bm}

\usepackage[colorlinks,allcolors=blue]{hyperref} 

\usepackage[utf8]{inputenc}
\usepackage[T1]{fontenc}
\usepackage{mathptmx}

\usepackage{subcaption}
\usepackage{caption}

\usepackage{soul} 

\usepackage{calrsfs}

\DeclareMathAlphabet{\pazocal}{OMS}{zplm}{m}{n}
\DeclareMathAlphabet\mathbfcal{OMS}{cmsy}{b}{n}
\SetMathAlphabet\pazocal{bold}{OMS}{zplm}{bx}{n}

\newcommand{\pza}{\pazocal{A}}

\newcommand{\pzh}{\pazocal{H}}
\newcommand{\pzl}{\pazocal{L}}

\newcommand{\pzs}{\pazocal{S}}
\newcommand{\pzz}{\pazocal{Z}}

\newcommand{\pzg}{\pazocal{G}}

\usepackage[dvipsnames]{xcolor}
\definecolor{AUSgreen}{RGB}{0, 132, 61}

\newcommand{\MATHUCC}{School of Mathematical Sciences, University College Cork, Cork T12 XF62, Ireland.}
\newcommand{\INFANT}{INFANT Research Centre, University College Cork, T12 DC4A Cork, Ireland.}

\begin{document}

\preprint{AIP/123-QED}


\title[Confabulation dynamics in a reservoir computer: Filling in the gaps with untrained attractors]{Confabulation dynamics in a reservoir computer: Filling in the gaps with untrained attractors}

\author{Jack O'Hagan}
\affiliation{ 
\MATHUCC
}%

\author{Andrew Keane}
\affiliation{ 
\MATHUCC
}%

\author{Andrew Flynn}
\email{andrew\_flynn@umail.ucc.ie}
\affiliation{ 
\MATHUCC
}
\affiliation{ 
\INFANT
}%


\date{\today}

\begin{abstract}

Artificial Intelligence has advanced significantly in recent years thanks to innovations in the design and training of artificial neural networks (ANNs). Despite these advancements, we still understand relatively little about how elementary forms of ANNs learn, fail to learn, and generate false information without the intent to deceive, a phenomenon known as `confabulation'. To provide some foundational insight, in this paper we analyse how confabulation occurs in reservoir computers (RCs): a dynamical system in the form of an ANN. RCs are particularly useful to study as they are known to confabulate in a well-defined way: when RCs are trained to reconstruct the dynamics of a given attractor, they sometimes construct an attractor that they were not trained to construct, a so-called `untrained attractor' (UA). This paper sheds light on the role played by UAs when reconstruction fails and their influence when modelling transitions between reconstructed attractors. Based on our results, we conclude that UAs are an intrinsic feature of learning systems whose state spaces are bounded, and that this means of confabulation may be present in systems beyond RCs.

\end{abstract}

\maketitle

\begin{quotation}


As artificial intelligence (AI) becomes increasingly integrated into modern society, the reliability and trustworthiness of AI-based systems is a growing concern. One of the main factors driving this concern is that these systems are prone to `confabulation': generating false, but plausible, information. In this paper, we investigate the phenomenon of confabulation using a mathematically well-defined machine learning framework known as reservoir computing. Specifically, we consider training an artificial neural network in the form of a reservoir computer (RC) on tasks involving attractor reconstruction, a process which mimics how memories (i.e., attractors) are encoded in the brain. We investigate how RCs confabulate by generating attractors that they have not been trained to construct, so-called `untrained attractors' (UAs). By analogy, these UAs are similar to false memories. We explore key questions such as, how and why do UAs arise in RCs? What are the most common types of UAs? What design features of RCs contribute to the generation of UAs? Beyond attractor reconstruction, how do RCs fill in gaps in their knowledge when generating transitions between attractors?

\end{quotation}

\section{\label{sec:Intro}Introduction}
The creation of false information by artificial neural networks (ANNs) to fill gaps in their memory has become a well-documented issue, particularly since the rise of generative artificial intelligence (AI) technologies such as OpenAI’s \textit{ChatGPT}. For instance, when asked to perform a literature review on a particular topic, ChatGPT may provide information on a paper whose title appears to be genuine but the paper does not actually exist \cite{emsley2023chatgpt,walters2023fabrication}. This phenomenon is commonly referred to as a `hallucination'. In recognition of the scale of this problem, Cambridge Dictionary chose `hallucinate' as its 2023 word of the year\cite{WordOfTheYear}.

However, the use of hallucination to describe this phenomenon has sparked significant debate since hallucinations typically occur due to false sensory perceptions, such as seeing or hearing things that do not exist in reality. Many argue that `confabulation' provides a more accurate analogy \cite{mcgowan2023chatgpt_confab_misnomer,edwards2023chatgpt_confab,seth_2024}. This term originates from its use to describe the generation of certain types of false memories by humans, where gaps in memory are filled with false, yet usually plausible, information \cite{schnider2003spontaneous}. 
Confabulation often arises when normal neural activity has been disrupted, as seen in cases of Korsakoff’s syndrome \cite{arts2017korsakoff}, Alzheimer’s disease \cite{el2017ConfabAlz}, and other memory-related neurological disorders. These disruptions are believed to prevent the brain from accurately encoding, storing, or retrieving information, thereby leading to false or distorted memories. 
More specifically, there are two types of confabulations in the psychology literature (see \textcite{gilboa2002cognitive} and references therein). ``Momentary'' (or ``provoked'') confabulations are mild distortions of a pre-existing memory and might appear as embellishments or elaborations. ``Fantastic'' (or ``spontaneous'') confabulations are less common and describe unlikely, sometimes bizarre, events.

While the abilities and applications of ANNs in AI-based technologies continue to advance at a remarkable pace, unfortunately the same cannot be said for our ability to explain how even the most elementary forms of ANNs learn, fail to learn, or confabulate. This lack of explainability presents a critical challenge for ensuring the reliability and trustworthiness of AI-based technologies, particularly as their adoption accelerates across scientific, industrial, and societal domains \cite{bengio2024AIrisks,weidinger2022risktaxonomy}. In order to improve our understanding of how confabulations arise in relatively simple ANN designs, in this paper we choose to study the dynamics of confabulation in a `reservoir computer' (RC) \cite{Jaeger01ESN,Maass02_LSM,nakajima2021RCbook}.

A RC is a dynamical system that is typically realised as an ANN with recurrent connections and can be trained to perform a wide variety of tasks, with applications ranging from control \cite{zhai_LWKong2023control} to weather prediction \cite{arcomano2023hybridweather}, to name but a few. These RC-based ANNs typically consist of randomly realised input and internal layers, and an output layer that is trained in order for the RC to perform a particular task. This philosophy of not training the input and internal network weights has inspired researchers to use a variety of physical mediums as the internal layer of a RC; we recommend \textcite{tanaka19PhysRC} and \textcite{Nakajima20phyRC_intro} for more information on the many interesting examples of `physical RCs'.

RCs also provide an ideal framework to systematically study how confabulations occur in ANNs as they are publicly-available, relatively inexpensive to train, and mathematically well-defined. RCs are known to exhibit a distinct form of confabulation: when trained to reconstruct the dynamics of a given attractor, the RC can construct attractors it was not trained to reconstruct, so-called `untrained attractors' (UAs). UAs were first discovered by \textcite{flynn2021multifunctionality} by initialising the trained RC (an autonomous dynamical system) from many different initial conditions and observing that certain initial conditions were within the basin of attraction of the reconstructed attractor and others were not. 
RCs also provide a valuable bridge between the brain and AI. 
A RC's learning process mirrors certain aspects of how the brain learns new tasks \cite{LuBassett20_switching_learning}, the process of training a RC to reconstruct an attractor is also viewed in a similar vein to how the brain stores long-term  episodic memories \cite{kong2024memory}. 
Thus, by extension, when the RC fails to accurately reconstruct an attractor it instead produces a false memory. The case of a poor reconstruction is analogous to a momentary confabulation, where ``moments'' or parts of the trajectory are inconsistent with the original attractor (e.g., see Fig.~\ref{fig:classification_algorithm_illustration_of_mechanics}). An UA is an example of a fantastic confabulation, where the ``complete'' trajectory is inconsistent with the original attractor. Hence, we refer to this case as ``complete confabulation'' for the rest of the paper. As we will see, it mirrors how gaps in a human's memory may be replaced with false information.

Taking into consideration the above associations between the brain, memory, and a RC's learning process, we conduct experiments to investigate the dynamics of confabulation when (i) storing an individual memory/attractor, (ii) interpolating between related memories/attractors, and (iii) interpolating between unrelated memories/attractors (attractors from different systems that do not share the same nonlinearities).
To be more specific, we study the dynamics of confabulation after training RCs to perform the following three tasks: (i) reconstruct the dynamics of the Lorenz attractor, (ii) reconstruct a period-1 and period-4 limit cycle at two different parameter settings, and (iii) reconstruct the Lorenz and Halvorsen attractor at two different parameter settings. For task (i) we train the RC presented in \textcite{LuHuntOtt18RC} and adapt this RC for tasks (ii) and (iii) by adding an additional `parameter-aware' component using a simplified version of the schemes introduced by \textcite{kim2021GlobalLocal} and \textcite{kong21PredCritTrans}.

From these tasks we identify key parameter relationships that influence the likelihood of UAs appearing. Most importantly, we uncover the gap-filling mechanism that gives rise to UAs: the state of the RC has no choice but to approach an UA when it completely fails to reconstruct an attractor, since, like many other learning systems, the RC we work with is an example of a bounded system, i.e. its state cannot approach infinity. 

Another motivation behind these tasks is that they help us to understand how an ANN-based system generalises beyond its training data and what gap filling mechanisms it may use. While these tasks differ from what one would ordinarily associate with the capabilities of highly sophisticated networks of interest, such as the human brain or generative AI models, they are chosen to study some of the circumstances in which confabulations arise when storing and recalling a piece of information (task (i)),  generalising between related pieces of information (task (ii)), and  generalising between unrelated pieces of information (task (iii)). 

When discussing confabulation mechanisms in the wider literature, researchers from both the psychology and machine learning communities argue that confabulations may arise when the brain or some learning system is `filling in the gaps' in its memory in cases where an accurate solution is not found for a particular query. 
For instance, in the APA Dictionary of Psychology, confabulations are defined as ``the falsification of memory in which gaps in recall are filled by fabrications that the individual accepts as fact'' \cite{ConfabPsychDef}. In a machine learning context, \textcite{nejjar2025llms_confab} describe confabulations in terms of ``...where an LLM [large language model] seemingly fills gaps in the information contained in the model with plausibly-sounding words''. 
However, to clearly identify and visualise such `gaps' being filled on a network level is challenging due to the complexity of these network systems. 
In contrast, the value of our RC approach is that it is mathematically tractable and we are able to produce diagrams of intended, learned solutions with gaps that are naturally filled by alternative, sometimes unexpected, solutions. 
More specifically, our results illustrate that what we observe as confabulations is literally the presence of unexpected solutions (UAs) ``filling in the gaps''. 

The rest of the paper is organised as follows, in Sec.~\ref{sec:Experiment} we describe the experiments the conducted and the steps involved in training the RC to do so, in Secs.~\ref{sec:Task1_results}-\ref{sec:Task3_results} we outline the results from training the RC to perform tasks (i)--(iii), 
and in Sec.~\ref{sec:DisConc} we provide some final remarks.

\section[Setup of numerical experiments]{Setup of numerical experiments}\label{sec:Experiment}

\subsection{\label{sec:RCintro} RC setup}

The RC studied throughout this paper was presented by \textcite{LuHuntOtt18RC}. There are two stages involved in training this RC to perform the tasks considered in this paper, a listening stage and a training stage. In both stages the RC, defined as the nonautonomous dynamical system in Eq.\,\eqref{eq:ListenRes},
\begin{align}
    \dot{\boldsymbol{r}}(t) &= \gamma \left[ - \boldsymbol{r}(t) + \tanh{\left( \textbf{M} \, \boldsymbol{r}(t) + \sigma \textbf{W}_{in} \, \boldsymbol{u}(t) \right)} \right], \label{eq:ListenRes}\\
    \boldsymbol{r}\left( 0 \right) &= \boldsymbol{0}^{T},\label{eq:ListenResIC}
\end{align}
is driven by an input $\boldsymbol{u}(t) \in \pza \subset \mathbb{R}^{D}$, where $\pza$ is some attractor that we choose and $D$ is the number of variables of the system that $\pza$ belongs to. We refer to this network driven by the input $\boldsymbol{u}(t)$ as the `open-loop RC'. The response to this driving input is used for training purposes and is given by values of $\boldsymbol{r}(t) \in \mathbb{R}^{N}$, which describes the state of the open-loop RC at a given time $t$ and $N$ is the number of artificial neurons in the network. Values of $\boldsymbol{r}(t)$ are obtained by computing solutions of Eq.\,\eqref{eq:ListenRes} using the 4$^{th}$ order Runge-Kutta method with time step $\tau$. $\gamma$ is a decay-rate parameter arising from the derivation of this RC from the discrete-time system proposed by \textcite{Jaeger01ESN}. The hyperbolic tangent `activation function' is a pointwise operation and is defined as $\tanh\left( \cdot \right) : \mathbb{R}^{N} \to \mathbb{R}^{N}$. The internal connections of the network are described by the adjacency matrix, $\textbf{M} \in \mathbb{R}^{N \times N}$. The input strength parameter, $\sigma$, and the input matrix, $\textbf{W}_{in} \in \mathbb{R}^{N \times D}$, when multiplied together represent the weight given to $\boldsymbol{u}(t)$, as it is fed into the open-loop RC. The steps involved in constructing $\textbf{M}$ and $\textbf{W}_{in}$ are outlined in Appendix\,\ref{app:RCdesign}.

In our experiments, we vary the `spectral radius' of $\textbf{M}$, the maximum of the absolute values of the eigenvalues of a given matrix. We denote this quantity by $\rho \geq 0$. Thus, for $\textbf{M}$ to have a particular $\rho$, we rescale its elements accordingly. We show that $\rho$ plays a crucial role in the confabulation dynamics of this RC design. $\rho$ has also been a key parameter in our previous results on RCs overcoming the influence of UAs \cite{flynn2021multifunctionality,flynn2023seeingdouble,morraflynn23_MF_fly,flynn2024switching} and features prominently in numerous research papers that focus on understanding the dynamics of how a RC learns.

\subsection{\label{sec:AttReconTrain}Training RC to reconstruct an attractor}

Let us first consider the case of training this RC to reconstruct the dynamics of a single attractor, $\pza \subset \mathbb{R}^{D}$, when given access to a trajectory on $\pza$ described by a vector $\boldsymbol{u}(t) \in \pza$. 

In this case, the listening stage is the duration of time where Eq.\,\eqref{eq:ListenRes} is driven by $\boldsymbol{u}(t)$ for $0 \leq t \leq t_{\text{listen}}$ where $t_{\text{listen}}$ is chosen such that $\boldsymbol{r}(t)$ is determined by a history of driving inputs and is no longer dependent on its initial condition. The training stage is the time where Eq.\,\eqref{eq:ListenRes} is driven by $\boldsymbol{u}(t)$ for $t_{\text{listen}} \leq t \leq t_{\text{train}}$. 

Training consists of finding a `readout function/layer' defined as $\hat{\boldsymbol{\psi}}\left( \cdot \right): \mathbb{R}^{2 N} \to \mathbb{R}^{D}$ that replaces $\boldsymbol{u}(t)$ in Eq.\,\eqref{eq:ListenRes} and is written as 
\begin{align}
    \hat{\boldsymbol{\psi}}\left(\boldsymbol{r}(t)\right) = \textbf{W}_{out} \boldsymbol{q}( \boldsymbol{r}(t) ),\label{eq:ReadoutFunction}
\end{align}
where $\textbf{W}_{out} \in \mathbb{R}^{D \times 2N}$ is the `readout matrix' and $\boldsymbol{q}( \boldsymbol{r}(t) ) \in \mathbb{R}^{2 N}$ is given by
\begin{align}
    \boldsymbol{q}(\boldsymbol{r}(t))=\left(\begin{array}{c}
\boldsymbol{r}(t)\\
\boldsymbol{r}^{2}(t)
\end{array}\right),\label{eq:q_square}
\end{align}
where $\boldsymbol{r}^{2}(t) = \left( r_{1}^{2}(t), r_{2}^{2}(t), \ldots, r_{N}^{2}(t) \right)^{T}$. From Eq.\,\eqref{eq:q_square} it follows that Eq.\,\eqref{eq:ReadoutFunction} can be rewritten as
\begin{align}
    \hat{\boldsymbol{\psi}}\left(\boldsymbol{r}(t)\right) = \textbf{W}_{out}^{(1)} \, \boldsymbol{r}(t) + \textbf{W}_{out}^{(2)} \, \boldsymbol{r}^{2}(t), \label{eq:BreakSym_r2}
\end{align}
where $\textbf{W}_{out}^{(1)}$ is the `linear readout matrix', and $\textbf{W}_{out}^{(2)}$ is the `square readout matrix'. The square readout matrix breaks the symmetry in Eq.\,\eqref{eq:ListenRes} when replacing $\boldsymbol{u}(t)$ with $\hat{\boldsymbol{\psi}}\left(\boldsymbol{r}(t)\right)$ after the training stage and prevents the occurrence of `mirror-attractors', which can impede the ability of the RC to reconstruct attractors\cite{herteux2020Symm,flynn2021symmetry}.

$\textbf{W}_{out}$ is determined by the ridge regression technique, which consists of minimising the expression
\begin{align}
    \frac{1}{t^{*}-l^{*}}\left[\sum_{i=l^{*}}^{t^{*}} || \textbf{W}_{out}\,\boldsymbol{q}(\boldsymbol{r}[i]) - \boldsymbol{u}[i]||_{2}^{2} + \beta \, || \textbf{W}_{out} ||_{2}^{2}\right]
    \label{eq:minimiseWout}
\end{align}
with respect to $\textbf{W}_{out}$, where $|| \cdot ||_{2}$ is the Euclidean norm, $l^{*} = t_{\text{listen}}/\tau$ and $t^{*} = t_{\text{train}}/\tau$, $\boldsymbol{r}[i]$ and $\boldsymbol{u}[i]$ are discrete-time samples of the continuous-time variables $\boldsymbol{r}(t)$ and $\boldsymbol{u}(t)$ at discrete time $i = t/\tau$. 
The corresponding time series, $\left[ \boldsymbol{r}[i]\right]_{i=l^{*}}^{t^{*}}$ and $\left[ \boldsymbol{u}[i]\right]_{i=l^{*}}^{t^{*}}$, are constructed by sampling $\boldsymbol{r}(t)$ and $\boldsymbol{u}(t)$ at intervals of length $\tau$ (discrete-time intervals of length $1$). $\beta$ is the regularisation parameter and the purpose of the $\beta \, || \textbf{W}_{out} ||_{2}^{2}$ term in Eq.\,\eqref{eq:minimiseWout} is to modify the linear least-squares regression to reduce the magnitudes of elements in $\textbf{W}_{out}$ in order to discourage overfitting.

Minimising Eq.\,\eqref{eq:minimiseWout} involves computing its partial derivative with respect to $\textbf{W}_{out}$ and setting the resulting expression equal to $0$. The $\textbf{W}_{out}$ that minimises Eq.\,\eqref{eq:minimiseWout} is given by
\begin{align}
    \textbf{W}_{out} = \textbf{Y} \textbf{X}^{T} \left( \textbf{X} \textbf{X}^{T} + \beta \, \textbf{I} \right)^{-1},\label{eq:WoutRegression}
\end{align}
where
\begin{align}
    \textbf{X} = \left[ \begin{array}{cccc}
    \boldsymbol{q}(\boldsymbol{r}[l^{*}]) & \boldsymbol{q}(\boldsymbol{r}[l^{*}+1]) 
    &
    \cdots 
    &
    \boldsymbol{q}(\boldsymbol{r}[t^{*}])
    \end{array} \right]
    \label{eq:Xmat}
\end{align}
is the `response data matrix',
\begin{align}
    \textbf{Y} = \left[ \begin{array}{cccc}
        \boldsymbol{u}[l^{*}] & \boldsymbol{u}[l^{*}+1] & \cdots & \boldsymbol{u}[t^{*}]
    \end{array} \right]
    \label{eq:Ymat}
\end{align}
is the `input data matrix' and $\textbf{I}$ is the identity matrix.

After the training, $\boldsymbol{u}(t)$ in Eq.\,\eqref{eq:ListenRes} is replaced by $\hat{\boldsymbol{\psi}} \left( \boldsymbol{r}(t) \right)$. In Eq.\,\eqref{eq:PredRes} we now define the `closed-loop' RC as the following autonomous dynamical system,
\begin{align}
    \hspace{-0.1cm}\dot{\hat{\boldsymbol{r}}}(t) &= \gamma \left[ - \hat{\boldsymbol{r}}(t) + \tanh{\left( \textbf{M} \, \hat{\boldsymbol{r}}(t) + \sigma \textbf{W}_{in}  \textbf{W}_{out} \, \boldsymbol{
    q}(\hat{\boldsymbol{r}}(t)) \right)} \right], \label{eq:PredRes}\\
    \hspace{-0.1cm}\hat{\boldsymbol{r}}\left( 0 \right) &= \boldsymbol{r}\left( t_{train} \right),\label{eq:PredResIC}
\end{align}
where $\hat{\boldsymbol{r}}(t)$ denotes the state of the closed-loop RC at a given time $t$. To distinguish between the $N$-dimensional vectors, $\hat{\boldsymbol{r}}(t)$ and $\boldsymbol{r}(t)$, we define $\hat{\boldsymbol{r}}(t) \in \mathbb{S}$, where $\mathbb{S}$ is referred to as the `RC's state space' and is used henceforth when discussing the dynamics of the closed-loop RC. By computing solutions of Eq.\,\eqref{eq:PredRes}, predictions of $\boldsymbol{u}(t)$ for $t>t_{train}$, denoted as $\hat{\boldsymbol{u}}(t)$, are given by
\begin{align}
    \hat{\boldsymbol{u}}(t) = \hat{\boldsymbol{\psi}}\left( \hat{\boldsymbol{r}}(t)\right).\label{eq:uhat}
\end{align}
Again, while both $\boldsymbol{u}(t)$ and $\hat{\boldsymbol{u}}(t)$ are $D$-dimensional vectors, to distinguish the dynamics of $\hat{\boldsymbol{u}}(t)$ from $\boldsymbol{u}(t)$, we define $\hat{\boldsymbol{u}}(t) \in \mathbb{P}$, where $\mathbb{P}$ is referred to as the `projected state space' and is used henceforth when discussing the dynamics of the closed-loop RC projected via $\hat{\boldsymbol{\psi}}$.

We say that the closed-loop RC has achieved attractor reconstruction when the long-term dynamical characteristics of $\hat{\boldsymbol{u}}(t)$ are `indistinguishable' from $\boldsymbol{u}(t)$, i.e. there exists an attractor $\pzs \subset \mathbb{S}$ such that when the state of the closed-loop RC approaches $\pzs$ and is projected from $\mathbb{S}$ to $\mathbb{P}$ via $\hat{\boldsymbol{\psi}}\left( \cdot \right)$, the dynamics of the `reconstructed attractor', $\hat{\pza} \subset \mathbb{P}$, resembles the dynamics of $\pza$ in the long-term. By resembling the long-term dynamics it is meant that, for instance, $\pza$ and $\hat{\pza}$ will have nearly identical Poincar{\'e} sections when computed for the same region of $\mathbb{R}^{D}$ and $\mathbb{P}$ as $t \to \infty$ for $t > t_{\text{train}}$. The method we use to determine this is specified in Appendix~\ref{apx:OutputClassification}.

\subsection{\label{sec:BifReconTrain}Training a RC to reconstruct attractors at different parameter settings}

We use a simplified version of the training procedure outlined in \textcite{kim2021GlobalLocal} and \textcite{kong21PredCritTrans} to train a RC to reconstruct two different attractors, $\pza_{1}$ and $\pza_{2}$, each at a different parameter settings of the RC. To do this, we adapt the $\tanh{\left( \cdot \right)}$ activation function in Eq.~\eqref{eq:ListenRes} to include a constant bias vector, $\boldsymbol{b} \in \mathbb{R}^{N}$, which is set to different values when reconstructing $\pza_{1}$ or $\pza_{2}$. The resulting `parameter-aware' open-loop RC is described as follows,
\begin{align}
    \dot{\boldsymbol{r}}(t) &= \gamma \left[ - \boldsymbol{r}(t) + \tanh{\left( \textbf{M} \, \boldsymbol{r}(t) + \sigma \textbf{W}_{in} \, \boldsymbol{u}(t) + \boldsymbol{b} \right)} \right]. \label{eq:ListenResPA}
\end{align}
The same initial condition specified in Eq.~\eqref{eq:ListenResIC} and integration scheme are used to generate solutions of Eq.~\eqref{eq:ListenResPA}.

There are also two stages involved in training this RC, a listening and training stage. In both stages, we first drive the parameter-aware open-loop RC in Eq.\,\eqref{eq:ListenResPA} with input $\boldsymbol{u}_{\left(\pza_{1}\right)}(t) \in \pza_{1}$ for $0 < t \leq t_{\text{train}}$ with $\boldsymbol{b} = \boldsymbol{b}_{1} = (b) \boldsymbol{1}^{T}$ and then repeat for $\boldsymbol{u}_{\left(\pza_{2}\right)}(t) \in \pza_{2}$ with $\boldsymbol{b} = \boldsymbol{b}_{2} = -(b) \boldsymbol{1}^{T}$ for a given $b > 0$. This choice of setting $\boldsymbol{b}_{1} = - \boldsymbol{b}_{2}$ is to further simplify the training procedure.

The responses of the open-loop RC to these driving inputs are denoted by $\boldsymbol{r}_{\left(\pza_{1}\right)}(t)$ and $\boldsymbol{r}_{\left(\pza_{2}\right)}(t)$.
It is important to highlight that $\textbf{M}$, $\textbf{W}_{in}$, and all training parameters remain identical when generating $\boldsymbol{r}_{\left(\pza_{1}\right)}(t)$ and $\boldsymbol{r}_{\left(\pza_{2}\right)}(t)$ in this setup.

The same readout function design in Eq.\,\eqref{eq:ReadoutFunction} is used. To distinguish between the $\textbf{W}_{out}$ obtained from training the RC to reconstruct a single attractor and a collection of attractors, we denote the readout matrix used in the parameter-aware RC setup as $\textbf{W}_{out}^{C}$ and the ridge regression approach for computing $\textbf{W}_{out} \in \mathbb{R}^{D \times 2N}$ in Eq.\,\eqref{eq:minimiseWout} is modified to account for additional attractors. This consists of minimising the following with respect to $\textbf{W}_{out}^{C}$,
\begin{align}
    &\frac{1}{t^{*}-l^{*}}\Bigg[\sum_{i=l^{*}}^{t^{*}} || \textbf{W}_{out}^{C}\,\boldsymbol{q}(\boldsymbol{r}_{\left(\pza_{1}\right)}[i]) - \boldsymbol{u}_{\left(\pza_{1}\right)}[i]||_{2}^{2} \nonumber \\
    & + \sum_{i=l^{*}}^{t^{*}} || \textbf{W}_{out}^{C}\,\boldsymbol{q}(\boldsymbol{r}_{\left(\pza_{2}\right)}[i]) - \boldsymbol{u}_{\left(\pza_{2}\right)}[i]||_{2}^{2} + \beta \, || \textbf{W}_{out}^{C} ||_{2}^{2}\Bigg],\label{eq:minimiseWoutMF}
\end{align}
where $\boldsymbol{r}_{\left(\pza_{1}\right)}[i]$, $\boldsymbol{u}_{\left(\pza_{1}\right)}[i]$, $\boldsymbol{r}_{\left(\pza_{2}\right)}[i]$, and $\boldsymbol{u}_{\left(\pza_{2}\right)}[i]$ are discrete-time samples of $\boldsymbol{r}_{\left(\pza_{1}\right)}(t)$, $\boldsymbol{u}_{\left(\pza_{1}\right)}(t)$, $\boldsymbol{r}_{\left(\pza_{2}\right)}(t)$, and $\boldsymbol{u}_{\left(\pza_{2}\right)}(t)$, constructed using the same convention outlined in Sec.\,\ref{sec:RCintro}. $\beta$ plays the same role as discussed in Sec.\,\ref{sec:RCintro}. The $\textbf{W}_{out}^{C}$ that minimises Eq.\,\eqref{eq:minimiseWoutMF} is given by
\begin{align}
    \textbf{W}_{out}^{C} = \textbf{Y}_{C} \textbf{X}_{C}^{T} \left( \textbf{X}_{C} \textbf{X}_{C}^{T} + \beta \, \textbf{I} \right)^{-1}.\label{eq:WoutRegressionMF}
\end{align}
$\textbf{X}_{C} = \left[ \textbf{X}_{\left(\pza_{1}\right)}, \, \textbf{X}_{\left(\pza_{2}\right)} \right]$ and $\textbf{Y}_{C} = \left[  \textbf{Y}_{\left(\pza_{1}\right)}, \, \textbf{Y}_{\left(\pza_{2}\right)} \right]$, where $\textbf{X}_{\left(\pza_{1}\right)}$ and $\textbf{X}_{\left(\pza_{2}\right)}$ are constructed as in Eq.\,\eqref{eq:Xmat} for the corresponding $\boldsymbol{r}_{\left(\pza_{1}\right)}[i]$ and $\boldsymbol{r}_{\left(\pza_{2}\right)}[i]$ and similarly for $\textbf{Y}_{\left(\pza_{1}\right)}$ and $\textbf{Y}_{\left(\pza_{2}\right)}$ as in Eq.\,\eqref{eq:Ymat} for $\boldsymbol{u}_{\left(\pza_{1}\right)}[i]$ and $\boldsymbol{u}_{\left(\pza_{2}\right)}[i]$. In Eq.\,\eqref{eq:WoutRegressionMF}, $\textbf{I}$ is the identity matrix of the appropriate size.

The parameter-aware closed-loop RC is then described by the same autonomous continuous-time dynamical system as in Eq.\,\eqref{eq:PredRes} with the additional $\boldsymbol{b}$ term included in the activation function and using the $\textbf{W}_{out}^{C}$ obtained from Eq.\,\eqref{eq:WoutRegressionMF}. For completeness, the parameter-aware closed-loop RC setup we study is described as follows,
\begin{align}
    \hspace{-0.2cm}\dot{\hat{\boldsymbol{r}}}(t\hspace{-0.02cm}) \hspace{-0.055cm}=\hspace{-0.055cm} \gamma \hspace{-0.055cm}\left[ - \hat{\boldsymbol{r}}(t\hspace{-0.02cm}) \hspace{-0.03cm}+\hspace{-0.03cm} \tanh{\left( \textbf{M} \hspace{-0.01cm} \hat{\boldsymbol{r}}(t\hspace{-0.02cm}) \hspace{-0.03cm}+\hspace{-0.03cm} \sigma \textbf{W}_{in}  \textbf{W}_{out}^{C} \, \boldsymbol{
    q}(\hat{\boldsymbol{r}}(t\hspace{-0.02cm})) \hspace{-0.03cm}+\hspace{-0.03cm} \boldsymbol{b} \right)} \right]. \label{eq:ParamAwarePredRes}
\end{align}

To reconstruct the dynamics of either $\pza_{1}$ or $\pza_{2}$, the $\boldsymbol{b}$ term in the activation function needs to be chosen to correspond with the value used during the training and the parameter-aware closed-loop RC needs to be initialised with the corresponding $\hat{\boldsymbol{r}}\left( 0 \right) = \boldsymbol{r}_{\left(\pza_{1}\right)}\left( t_{train} \right)$ or $\boldsymbol{r}_{\left(\pza_{2}\right)}\left( t_{train} \right)$ (or some other initial condition that is in the basin of attraction of the corresponding $\pzs_{1}$ or $\pzs_{2}$).

The same steps outlined above can be extended in order to produce a parameter-aware closed-loop RC that reconstructs $n$ attractors, \{$\pza_{1}, \, \pza_{2}, \, \ldots, \, \pza_{n}$\}, at parameter settings $\{\boldsymbol{b}_{1}, \boldsymbol{b}_{2}, \ldots, \boldsymbol{b}_{n}\}$, respectively. All that is required is to produce the corresponding $\textbf{X}_{C} = \left[ \textbf{X}_{\left(\pza_{1}\right)}, \, \textbf{X}_{\left(\pza_{2}\right)}, \, \ldots, \, \textbf{X}_{\left(\pza_{n}\right)} \right]$ and $\textbf{Y}_{C} = \left[  \textbf{Y}_{\left(\pza_{1}\right)}, \, \textbf{Y}_{\left(\pza_{2}\right)}, \, \ldots, \, \textbf{Y}_{\left(\pza_{n}\right)} \right]$ and solve for $\textbf{W}_{out}^{C}$ as in Eq.\,\eqref{eq:WoutRegressionMF}.

It is also worth mentioning that if the same process is repeated with $\boldsymbol{b}=\textbf{0}$, and the RC reconstructs each attractor in a co-existing fashion, then the RC it said to exhibit multifunctionality.

\subsection{Defining untrained and generated attractors}\label{ssec:UA_define}

Before proceeding any further, it is important to define what we consider to be an UA in terms of training a RC to reconstruct an attractor or to reconstruct attractors at different parameter settings via the steps outlined in Secs.~\ref{sec:AttReconTrain} and \ref{sec:BifReconTrain}.

\textbf{Definition\,1:}\,\textit{An attractor, $\pzz$, is considered to be an `untrained attractor' if (i) $\pzz$ was not present during training, and (ii) $\pzz$ does not originate from a continuous deformation/displacement or bifurcation of a given reconstructed attractor by varying certain parameters of the given reservoir computer.}

This definition is consistent with our previous studies of UAs in the context of multifunctionality\cite{flynn2021multifunctionality,flynn2023seeingdouble,morraflynn23_MF_fly,flynn2024switching}, as we referred to the attractors that the RC was explicitly trained to reconstruct as the `reconstructed attractors', and the attractors that the RC was not explicitly trained to reconstruct as `untrained attractors'. However, when it comes to discussing the dynamics of a parameter-aware closed-loop RC, one must take further caution on how to refer to attractors that the RC was not explicitly trained to reconstruct that are also not UAs. 

Let us consider the following, a parameter-aware RC is explicitly trained to reconstruct attractors at different parameter settings, in our case the values of $\boldsymbol{b}$. Thus, if the training is successful, then for each of the chosen values of $\boldsymbol{b}$, there exists a corresponding reconstructed attractor. However, if the parameter-aware RC interpolates/extrapolates the dynamics of the reconstructed attractors between/beyond the chosen values of $\boldsymbol{b}$, which involves the reconstructed attractors undergoing either a continuous deformation/displacement or a series of bifurcations, then it is fair to say that the RC was not explicitly trained to reconstruct such attractors. However, to call any of these attractors an UA would contradict Definition 1. To resolve this issue we introduce the following definition.

\textbf{Definition\,2:}\,\textit{An attractor, $\pzg$, is considered to be a `generated attractor' if (i) $\pzg$ was not present during training, and (ii) $\pzg$ originates from a continuous deformation/displacement or bifurcation of a given reconstructed attractor by varying certain parameters of the given reservoir computer.}

We choose the term `generated' to avoid having to refer to such an attractor as an interpolated/extrapolated attractor based on the value of $\boldsymbol{b}$. To clarify, we only discuss generated attractors (GAs) in the context of parameter-aware RCs.

We would also like to remark that Definition~2 is designed to compliment Definition~1 in terms of their subtle differences and also to apply to attractors of all types. As a result, Definition~2 may seem relatively trivial in the context of fixed-points and limit cycles. However, we believe that the strength of Definition~2 is its generality. 

We denote GAs by $\hat{\pza}^{*}$ to emphasise its connection to the reconstructed attractor $\hat{\pza}$. For completeness, when presenting our results, the attractor $\pza_{i}$ that we explicitly train the RC to reconstruct at parameter setting $\boldsymbol{b}_{i}$ is labelled by $\hat{\pza}_{i}$. The corresponding family of GAs that exist at values of $\boldsymbol{b}$ not used during training are labelled as $\hat{\pza}_{i}^{*}$. If, for instance, there is a family of GAs that connect a pair of reconstructed attractors, $\hat{\pza}_{i}$ and $\hat{\pza}_{j}$, then we refer to this family as $\hat{\pza}_{i,j}^{*}$.

\subsection{Description of tasks and experiments}\label{ssec:OutlineOfExperiments}

\begin{figure*}
    \centering
    \includegraphics[width=\linewidth]{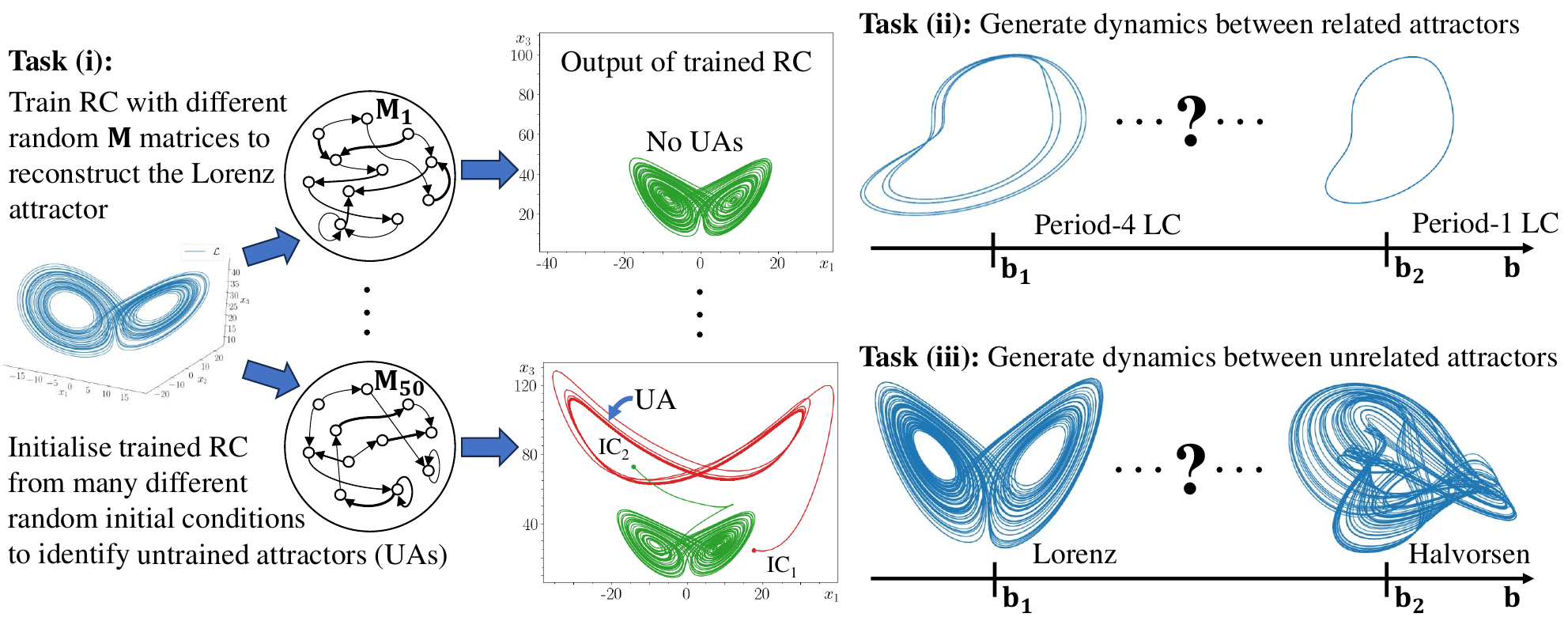}
    \caption{Overview of tasks considered to explore confabulation dynamics in a RC.}
    \label{fig:ConfDynm_Experiment_Overview}
\end{figure*}

In Secs.~\ref{sssec:Task1_description}-\ref{sssec:Task3_description}, we describe how the training data is generated for each of the three tasks mentioned in Sec.~\ref{sec:Intro} and outline the specifics of the experiments that are carried out in relation to each task. A schematic of the main aspects involved in each task is provided in Fig.~\ref{fig:ConfDynm_Experiment_Overview}.

\subsubsection{Task (i) description}\label{sssec:Task1_description}

For task (i) we focus on exploring the dynamics of confabulation when training the open-loop RC in Eq.~\eqref{eq:ListenRes} to reconstruct the dynamics of the Lorenz attractor, $\pzl$, by finding a suitable $\textbf{W}_{out}$ matrix according to the steps outlined in Sec.~\ref{sec:AttReconTrain}. 
For training data, we obtain a trajectory on $\pzl$ by generating solutions of the Lorenz system\cite{lorenz63model},
\begin{align}
    \begin{array}{ccl}
        \dot{x}_{1} &=& 10 \left( x_{2} + x_{1} \right),\\
        \vspace{-0.3cm}
        & & \\
        \dot{x}_{2} &=& x_{1} \left( 28 - x_{3} \right) - x_{2},\\
        \vspace{-0.3cm}
        & & \\
        \dot{x}_{3} &=& x_{1} x_{2} - \frac{8}{3} x_{3},
    \end{array}\label{eq:Lorenz}
\end{align}
until $t=t_{predict}$ using the 4$^{th}$ order Runge-Kutta method with time step $\tau$. For this task, we set $t_{listen} = 100$, $t_{train} = t_{listen}~+~100$, $t_{predict} =  t_{train}~+~100$, $\sigma = 0.2$, $\beta = 0.001$, $\gamma = 10$, and consider several different values for $\rho$.

\begin{figure*}
    \centering
    \includegraphics[width=0.95\linewidth]{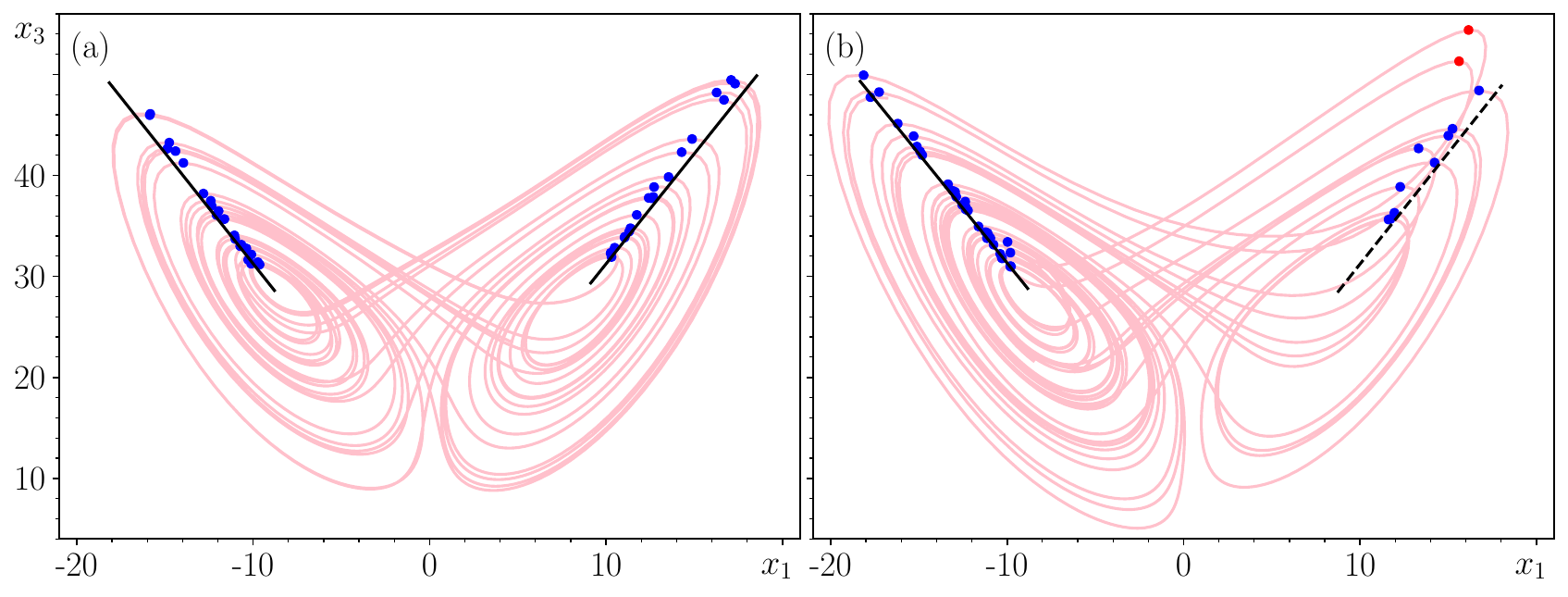}
    \caption{Illustration of decision process in distinguishing between a `good' reconstruction of $\pzl$ ($\hat{\pzl}$) in (a) and a `bad' reconstruction of $\pzl$ ($\tilde{\pzl}$) in (b). Blue points indicate local maxima that satisfy classification criteria, red point indicate local maxima that do not. Solid black lines indicate a good fit between $\pzl$ and its reconstruction, dashed black lines indicate a poor fit.}
    \label{fig:classification_algorithm_illustration_of_mechanics}
\end{figure*}

To assess how well the closed-loop RC reconstructs $\pzl$, we use the method outlined in Appendix~\ref{apx:OutputClassification}. From this we classify the output from a given initial condition as a `good' reconstruction of $\pzl$, denoted by $\hat{\pzl}$, a `poor' reconstruction of $\pzl$, denoted by $\tilde{\pzl}$. 
An example of good and poor reconstructions of $\pzl$, as determined through the steps outlined in Appendix~\ref{apx:OutputClassification}, is provided in Fig.~\ref{fig:classification_algorithm_illustration_of_mechanics}.
If there is no reconstruction of $\pzl$, i.e. there is no stable attractor located nearby $\pzl$ with properties resembling $\pzl$, we then say that the state of the closed-loop RC approaches an UA, corresponding to complete confabulation. 
We find a wide variation of dynamics corresponding to $\tilde{\pzl}$ and later show examples of these attractors in Sec.~\ref{sec:Task1_results} and Appendix~\ref{apx:AdditionalReconRoutes}.

$\textbf{M}$ and $\textbf{W}_{in}$ are realised randomly in this RC setup, therefore, RCs with different randomly realised $\textbf{M}$ and $\textbf{W}_{in}$ may result in different training outcomes. Thus, as indicated in Fig.~\ref{fig:ConfDynm_Experiment_Overview}, for a given parameter setting of the RCs we study, UAs may be present for one realisation of $\textbf{M}$ and $\textbf{W}_{in}$ and not present in other realisations. From this, many questions arise, such as how often do UAs appear and what is the most common type of UA that appears when reconstructing $\pzl$? 
However, exploring all realisations of $\textbf{M}$ and $\textbf{W}_{in}$ may not be so practical, depending on the size of the network and the degree of weight variation. 

To gain some insight into the questions posed above, we consider training an ensemble of RCs with 50 different randomly realised versions of $\textbf{M}$ and one random realisation of $\textbf{W}_{in}$ for $N=100$. Once these 50 realisations of $\textbf{M}$ have been generated, we rescale the weights of each $\textbf{M}$ to have a specific spectral radius, $\rho$. To stress this point further, we do not generate 50 new realisations of $\textbf{M}$ each time we change $\rho$. 
As indicated in Fig.~\ref{fig:ConfDynm_Experiment_Overview}, we initialise each of the closed-loop RCs from 100 different random initial conditions to determine whether there are any UAs present in each respective $\mathbb{P}$. We then classify the output of a given closed-loop RC (dynamics in $\mathbb{P}$) depending on the different attractors it approaches using the approach specified in Appendix~\ref{apx:OutputClassification}. 
Based on the resulting classifications, we identify five different scenarios on the type of dynamics seen in $\mathbb{P}$ that we specify in Table~\ref{tab:5scenarios}. In Sec.~\ref{sec:Task1_results} we show that the frequency of these scenarios vary depending on the choice of $\rho$.

\begin{table}[t]
    \centering
    \renewcommand{\arraystretch}{1.5}
    \begin{tabular}{|c|l|}
        \hline
        1. & The only attractor present in $\mathbb{P}$ is $\hat{\pzl}$.\\
        \hline
        2. & The only attractor present in $\mathbb{P}$ is $\tilde{\pzl}$.\\
        \hline
        3. & The only attractor(s) present in $\mathbb{P}$ are UA(s).\\
        \hline
        4. & $\hat{\pzl}$ coexists with one or several UAs in $\mathbb{P}$.\\
        \hline
        5. & $\tilde{\pzl}$ coexists with one or several UAs in $\mathbb{P}$.\\
        \hline
    \end{tabular}
    \caption{Attractors found in $\mathbb{P}$ after training Eq.~\eqref{eq:ListenRes} to reconstruct the Lorenz attractor, $\pzl$.}
    \label{tab:5scenarios}
\end{table}

Building on this, we investigate some of the most common ways in which $\hat{\pzl}$ appears. We find that these mainly appear to involve saddle type bifurcations. Interestingly, another common case in which $\hat{\pzl}$ appears is through a period-doubling route to chaos. What is also interesting is that we do not find any cases where $\hat{\pzl}$ comes into existence through the well-known sequence of events that unfold by varying the constant term (28) in the $\dot{x}_{2}$ equation in Eq.~\eqref{eq:Lorenz}. We discuss the ramifications of these results in greater detail in Sec.~\ref{sec:DisConc}.

We would like to clarify that the reason why we chose to conduct our experiments with 50 different realisations of $\textbf{M}$ as opposed to 50 different realisations of both $\textbf{M}$ and $\textbf{W}_{in}$ was to simply only vary one component of the RC at a given time. 
From experiments not presented in the present paper, we found that different UAs can also appear for several different random realisations of $\textbf{W}_{in}$ and a given random $\textbf{M}$. In future work we intend to investigate the relationship between confabulation dynamics and $\textbf{W}_{in}$ in much greater detail. 

\subsubsection{Task (ii) description}\label{sssec:Task2_description}

The purpose of task (ii) is to investigate the dynamics of confabulation that can arise after the parameter-aware RC has been trained to reconstruct two different, but \textit{related attractors} at two different parameter settings. By `related' we mean that there ordinarily exists a trivial sequence of transitions between these attractors. For task (ii) we consider training the parameter-aware open-loop RC in Eq.~\eqref{eq:ListenResPA} to reconstruct a period-1 and period-4 limit cycle at two different parameter settings according to the steps outlined in Sec.~\ref{sec:BifReconTrain}.

For training data in task (ii), we choose to generate solutions of the system studied in \textcite{lisprott2016crisis} to obtain a trajectory on a period-1 and period-4 limit cycle,
\begin{align}
    \begin{array}{ccl}
        \dot{x}_{1} &=& x_{2} \left( 1 + x_{3} \right),\\
        \vspace{-0.3cm}
        & & \\
        \dot{x}_{2} &=& x_{3} \left( x_{2} - a x_{1} \right),\\
        \vspace{-0.3cm}
        & & \\
        \dot{x}_{3} &=& 0.55 ~ x_{3}^{2} - x_{2}^{2}.
    \end{array}\label{eq:sprott}
\end{align}
For both limit cycles, we generate solutions of Eq.~\eqref{eq:sprott} up to $t=t_{predict}$ using the 4$^{th}$ order Runge-Kutta method with time step $\tau$. For this task, we set $~t_{listen} = 100$, $t_{train} = t_{listen}~+~200$, and $t_{predict} = t_{train}~+~200$. When $a=17$ we obtain a period-4 limit cycle and when $a=27$ we obtain a period-1 limit cycle, trajectories on these limit cycles are illustrated in the portion of Fig.~\ref{fig:ConfDynm_Experiment_Overview} associated with task (ii). 

\begin{figure}
    \centering
    \includegraphics[width=\linewidth]{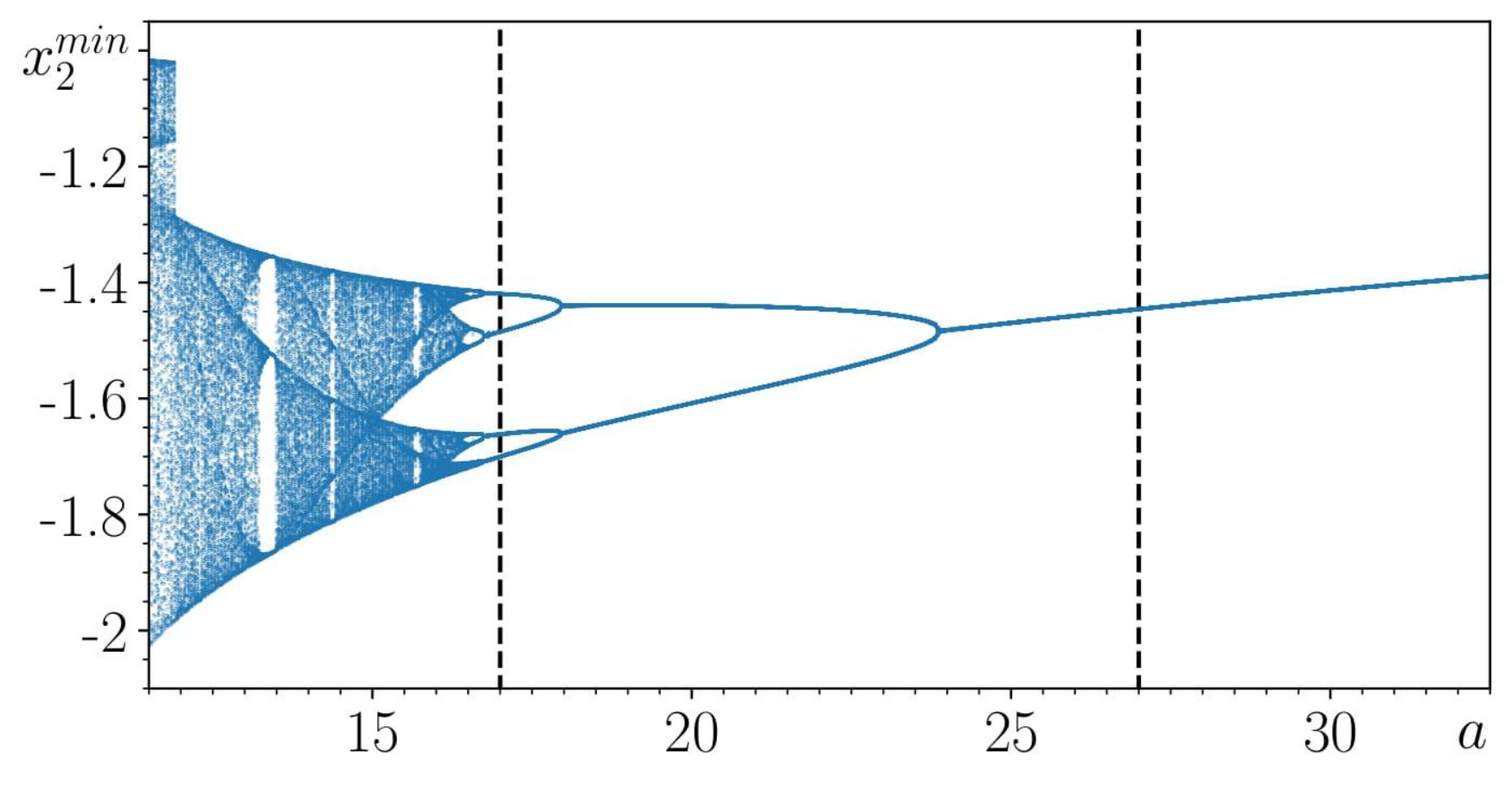}
    \caption{Bifurcation diagram to be reconstructed in task (ii). Changes in local $x_{2}$ minima of attractor from Eq.~\eqref{eq:sprott} with respect to changes in $a$. Black lines indicate values of $a$ used as training data.}
    \label{fig:SprottSys_BifDiagram}
\end{figure}

When driving the parameter-aware open-loop RC setup in Eq.~\eqref{eq:ListenResPA} with input from these limit cycles, we set $\boldsymbol{b} = \boldsymbol{b}_{1} = (b) \textbf{1}^{T}$ for the period-4 limit cycle, and $\boldsymbol{b} = \boldsymbol{b}_{2} = -(b) \textbf{1}^{T}$ for the period-1 limit cycle. Thus, through the convention used in Sec.~\ref{sec:BifReconTrain}, we denote the reconstructed period-4 limit cycle by $\hat{\pza_{1}}$ and its related GAs by $\hat{\pza_{1}}^{*}$. Similarly, we denote the reconstructed period-1 limit cycle by $\hat{\pza_{2}}$ and its related GAs by $\hat{\pza_{2}}^{*}$.

In Fig.~\ref{fig:SprottSys_BifDiagram} we show that for decreasing $a$, the period-1 limit cycle undergoes a sequence of period-doubling bifurcations. Therefore, if the parameter-aware closed-loop RC successfully interpolates the dynamics of Eq.~\eqref{eq:sprott} between the period-1 and period-4 limit cycles, then the RC should generate a period-2 limit cycle at some intermediate parameter setting between $b$ and $-b$. In addition, the RC should generate a continuous transition between this period-2 limit cycle and the reconstructed period-1 and period-4 limit cycles, where each transition should occur due to a period-doubling bifurcation. Furthermore, the RC should continue to undergo successive period-doubling bifurcations beyond the period-4 limit cycle to produce chaotic dynamics. However, as our results show, it is not so straightforward for the RC to fill the gaps in its knowledge of Eq.~\eqref{eq:sprott} by producing GAs that mimic the dynamics shown in Fig.~\ref{fig:SprottSys_BifDiagram}.

In our experiments, we investigate how different choices of the parameter $b$ affects the ability of the parameter-aware closed-loop RC in Eq.~\eqref{eq:ParamAwarePredRes} to reconstruct the period-1 and period-4 limit cycles at the respective $\boldsymbol{b}_{1}$ and $\boldsymbol{b}_{2}$ settings of $\boldsymbol{b}$. We find empirically that by setting $\rho = 1.2$, $\sigma = 1.6$, $\beta = 0.01$, and $\gamma = 10$, we obtain reasonable reconstruction of the period-1 and period-4 limit-cycles for $b=\{\pm0.2, \pm0.3, \pm0.4\}$. We find that reconstruction worsens for larger and smaller $b$, particularly for the period-4 limit cycle. After training, we track the changes in the reconstructed attractors in $b$-space for $b \in \left[ 0.42, -0.42\right]$. To do this we increase/decrease the value of $b$ used during training by a small amount and initialise Eq.~\eqref{eq:ParamAwarePredRes} at this new $b$ value with a point on the attractor from the previous $b$ value and continue this process until we can no longer track the subsequent GAs. We map the points on the bifurcation diagram in Fig.~\ref{fig:SprottSys_BifDiagram} to the set of $b$ values specified above for the case of $b = \pm0.2$. Thus, for $b = \pm0.3$ and $b=\pm0.4$, the resulting bifurcation diagrams from the RC correspond to rescalings of the original bifurcation diagram shown in Fig.~\ref{fig:SprottSys_BifDiagram}.

Another reason why we chose to train the RC using limit cycles from Eq.~\eqref{eq:sprott} is because there is another limit cycle that undergoes period-doubling route to chaos for the same values of $a$, with both attractors merging to form one larger attractor when $a \approx 12$ (evidenced by the sudden growth in the range of the $x_{2}$ local minima in Fig.~\ref{fig:SprottSys_BifDiagram}). This other limit cycle is a symmetric copy of the limit cycle we study, it exists due to the following symmetry relation of Eq.~\eqref{eq:sprott}, $\left( x_{1}, x_{2}, x_{3} \right) \rightarrow \left( -x_{1}, -x_{2}, x_{3} \right)$. In some of the parameter settings we explore, we find that the parameter-aware closed-loop RC generates an attractor with similar properties to this other limit cycle.

\subsubsection{Task (iii) description}\label{sssec:Task3_description}

In contrast to task (ii), the purpose of task (iii) is to investigate the dynamics of confabulation that can arise after the parameter-aware RC has been trained to reconstruct two different \textit{unrelated attractors} at two different parameter settings. By `unrelated' we mean that there exists no pre-defined sequence of transitions or continuous deformation from one attractor to the other. For task (iii) we focus on exploring the dynamics of confabulation when training the parameter-aware RC in Eq.~\eqref{eq:ListenResPA} to reconstruct the Lorenz attractor, $\pzl$, and the Halvorsen attractor, $\pzh$, at two different parameter settings according to the steps outlined in Sec.~\ref{sec:BifReconTrain}.

For training data, we obtain a trajectory on $\pzl$ according to the steps specified in Sec.~\ref{sssec:Task1_description}, we obtain a trajectory on $\pzh$ with the same procedure by generating solutions of the following system,
\begin{align}
    \begin{array}{ccl}
        \dot{x}_{1} &=& -1.3 x_{1} - 4 \left( x_{2} + x_{3} \right) - x_{2}^{2},\\
        \vspace{-0.3cm}
        & & \\
        \dot{x}_{2} &=& -1.3 x_{2} - 4 \left( x_{3} + x_{1} \right) - x_{3}^{2},\\
        \vspace{-0.3cm}
        & & \\
        \dot{x}_{3} &=& -1.3 x_{3} - 4 \left( x_{1} + x_{2} \right) - x_{1}^{2},
    \end{array}\label{eq:Halvorsen}
\end{align}
until $t=t_{predict}$ using the 4$^{th}$ order Runge-Kutta method with time step $\tau$. For this task, we set $t_{listen} = 100$, $t_{train} =  t_{listen}~+~200$, and $t_{predict} = t_{train}~+~200$.

When training the parameter-aware RC setup in Eq.~\eqref{eq:ListenResPA} to reconstruct $\pzl$, we set $\boldsymbol{b} = \boldsymbol{b}_{1} = (b) \textbf{1}^{T}$. For $\pzh$ we set $\boldsymbol{b} = \boldsymbol{b}_{2} = -(b) \textbf{1}^{T}$. We denote the reconstructed $\pzl$ at parameter setting $\boldsymbol{b}_{1}$ by $\hat{\pzl}$, and its corresponding GAs by $\hat{\pzl}^{*}$. Similarly, we denote the reconstruction of $\pzh$ at parameter setting $\boldsymbol{b}_{2}$ by $\hat{\pzh}$, and its corresponding GAs by $\hat{\pzh}^{*}$. Similar to task (ii), we investigate how different choices of the parameter $b$ impacts the ability of the parameter-aware closed-loop RC to reconstruct these attractors for $\boldsymbol{b} = \boldsymbol{b}_{1}$ and $\boldsymbol{b}_{2}$. 

\begin{figure}
    \centering
    \includegraphics[width=0.9\linewidth]{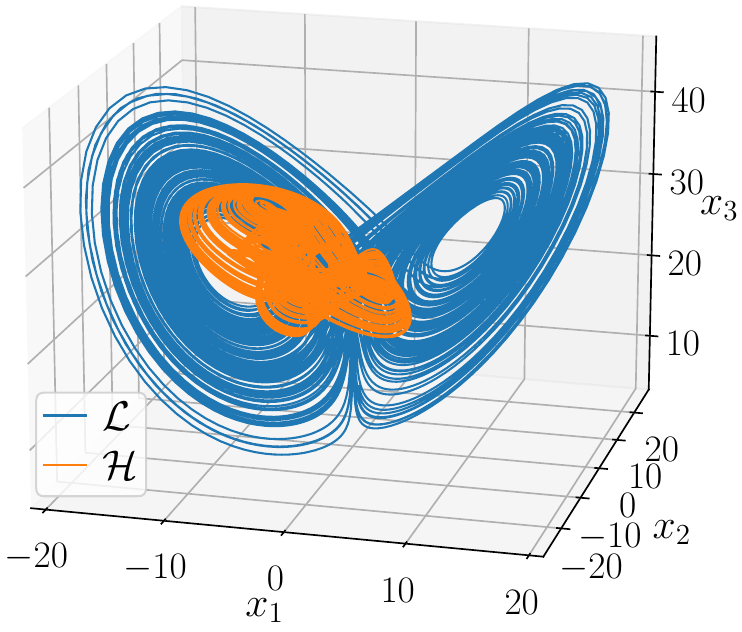}
    \caption{Training data for task (iii).}
    \label{fig:LorHal_trainingdata}
\end{figure}

In our investigations we concentrate on a set of RC training parameters for which the parameter-aware closed-loop RC reconstructs $\pzl$ and $\pzh$ reasonably well for $\boldsymbol{b} = \boldsymbol{b}_{1}$ and $\boldsymbol{b}_{2}$ for different choices of the constant $b$. We find that by setting $\rho = 1.2$, $\sigma = 0.2$, $\beta = 0.1$, and $\gamma = 10$, we obtain reasonable reconstruction of $\pzl$ and $\pzh$ for $b=\{\pm0.2, \pm0.3, \pm0.4\}$. 
In order to find such a set of training parameters, we had to manually move the training data corresponding to $\pzh$ to overlap with $\pzl$, the location of the attractors in the $(x_{1}, x_{2}, x_{3})$-plane is shown in Fig.~\ref{fig:LorHal_trainingdata}. Similar to task (ii), we also find that reconstruction of both attractors worsens for larger and smaller $b$. 

After training, we track the changes in the reconstructed attractors in $b$-space for $b \in \left[ 0.42, -0.42\right]$ using the same procedure mentioned in Sec.~\ref{sssec:Task2_description}. From this we uncover the different ways in which the RC interpolates and extrapolates beyond its training dataset. We display these results in diagrams that show the changes in the local minima of the $x_{1}$ variable in $\mathbb{P}$ for changes in $b$. We choose to illustrate these results in terms of the $x_{1}$ variable for clarity of presentation as the dynamics of both reconstructed attractors and their families of GAs are confined to the most similar range of values. From this we identify the following three different ways in which, depending on the choice of $b$, the parameter-aware closed-loop RC fills in gaps in its memory by constructing:
\begin{enumerate}
    \item An UA that coexists with $\hat{\pzl}^{*}$ and $\hat{\pzh}^{*}$ at different values of $\boldsymbol{b}$.
    \item A region of bistability between $\hat{\pzl}^{*}$ and $\hat{\pzh}^{*}$ for a small range of $\boldsymbol{b}$.
    \item A continuous transition between $\hat{\pzl}^{*}$ and $\hat{\pzh}^{*}$.
\end{enumerate}


\subsection{Summary of related work}\label{ssec:PreviousWork}

As mentioned during Sec.~\ref{sec:RCintro}, in our previous investigations of multifunctionality in a RC, we illustrated that in many cases changing the spectral radius, $\rho$, of the internal connections strongly influences the prevalence of UAs \cite{flynn2021multifunctionality,flynn2023seeingdouble,morraflynn23_MF_fly,flynn2024switching}. We show that this can also depend on the random realisation of $\textbf{M}$ from the results of task (i) in this paper. Interestingly, in studies of multistability in input driven RCs \cite{ceni2020interpreting,fadera2024arbitrary}, UAs also appeared, with the authors referring to these as `spurious attractors'. Upon further investigation, we found that the same term has been used by researchers studying behaviours of Hopfield networks that were not learned by the network, like in \textcite{robins2004spuriousatt}, however these studied have mainly focused on trained and untrained attractors in the form of fixed points.

While the presence of UAs in RC systems are being increasingly studied and acknowledged \cite{pecora2024statistics,kong2024memory,terasaki2024thermodynamic,Kabayama25_UAbif,tomiokanakajima2025SoftBodyUAs}, 
and have recently been detected in a physical RC system\cite{ryo25_multifunctionalphysicalRC}, 
it is quite likely that UAs have lurked in the background of many researchers dealings with RCs but have simply gone unnoticed. For instance, in one of the most influential papers on reservoir computing in recent times, \textcite{LuHuntOtt18RC}, the authors present an interesting result in their Fig.~5 where they show that for a relatively small change in the input scaling parameter, $\sigma$, there is an abrupt transition in the RC's dynamics from a period-2 limit cycle to the chaotic Lorenz attractor. Based on the results presented in the present paper, it is not too far-fetched to speculate that the limit cycle is an UA that coexists with the reconstructed Lorenz attractor. Other examples of UAs potentially hiding in plain sight can be seen in \textcite{deJong2025EvidenceOfUAs}. 


The connection between tasks (i)--(iii) and the process of storing and recalling a memory has been established by \textcite{LuBassett20_switching_learning} and \textcite{kong2024memory} who relate different processes of how the brain stores multiple memories simultaneously to results in reservoir computing. More specifically, Kong \textit{et al.} relate `context-addressable memory' to how a multifunctional RC stores different attractors in different regions of state space, and relate `location-addressable memory' to how a parameter-aware RC stores different attractors at different parameter settings. 
Another motivation behind tasks (ii) and (iii) is to build on the above by relating how a parameter-aware RC generalises beyond its training data by filling in gaps between these attractors to how the brain generalises between memories and also to studies of generalisation in the wider machine learning research community. 
These tasks are analogous to studies of how generative AI models interpolate in embedding or latent spaces (semantic spaces learned during training where similar concepts cluster together). The results we present may provide some further insight towards the dynamics exhibited by these models before the `grokking transition' occurs, when a model's ability to generalise beyond its training data changes suddenly from a very poor to a perfect generalisation \cite{liuTegmark2022UnderstandGrok}. 
In the context of our RC set-up, poor generalisations correspond to UAs and GAs that do not resemble the dynamics they are supposed to. 
Task (ii) enables us to identify factors that influence the generative properties of the parameter-aware RC by testing whether the transition from one attractor to another is representative of how transitions between these attractors naturally occur, which shares many similarities with the work presented in \textcite{struski2023Interpolate}. 
On the other hand, task (iii) provides us with a means to study these generative properties in a situation when there is no reference for how these transitions should occur. 
Moreover, tasks (ii) and (iii) represent studies of how different memories are linked together in a neural network via the dynamics of GAs and UAs, which in some sense relates to the ability of the hippocampus to form relationships between different memories in decision making \cite{zeithamova2012hippocampus}. 
It has also been reported that related and unrelated events can become linked in memory when they happen close together \cite{uitvlugt2019temporalproximity}. 
In the case of tasks (ii) and (iii), two events occurring within a brief window of time is similar to storing the attractors close together in $\boldsymbol{b}$-space, and events becoming linked in memory is similar to there being a continuous transition from one attractor to the other.

\section{Task (i) results}\label{sec:Task1_results}

\subsection{Frequency of scenarios in Table~\ref{tab:5scenarios} across 50 different \texorpdfstring{$\textbf{M}$}{TEXT}}\label{ssec:FreqOfScenarios}

\begin{figure}
    \centering
    \includegraphics[width=\linewidth]{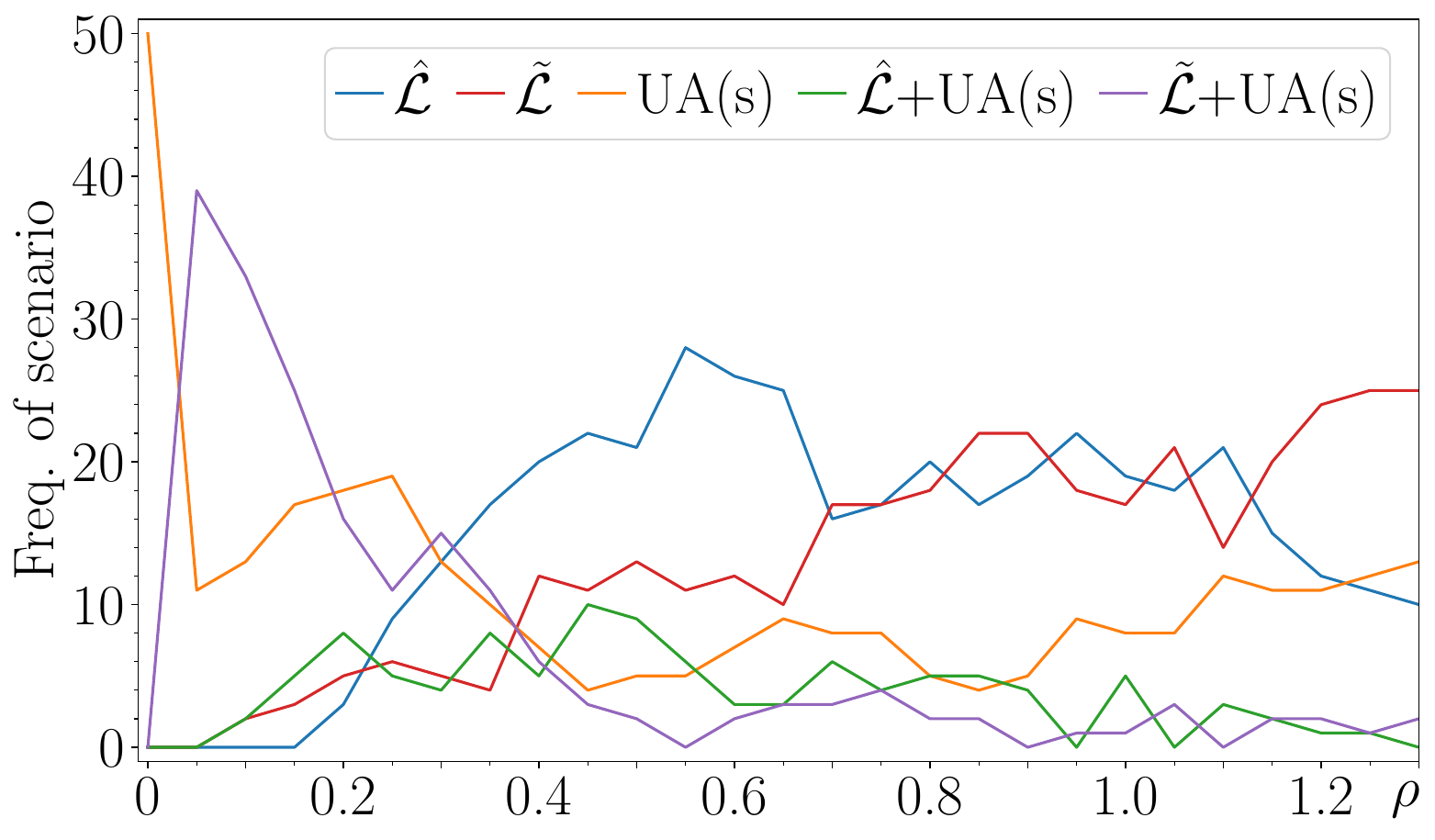}
    \caption{Frequency of each scenario specified in Table~\ref{tab:5scenarios} for a given $\rho$ when training Eq.~\eqref{eq:ListenRes} to reconstruct $\pzl$ using 50 different random realisations of $\textbf{M}$.}
    \label{fig:rho_vs_ScenarioFreq_updated}
\end{figure}

Fig.~\ref{fig:rho_vs_ScenarioFreq_updated} illustrates the frequency of each of the scenarios specified in Table~\ref{tab:5scenarios} for a given spectral radius, $\rho$, across 50 different random realisations of the internal connectivity matrix, $\textbf{M}$, whose elements have been scaled such that each $\textbf{M}$ has the same $\rho$. 

Fig.~\ref{fig:rho_vs_ScenarioFreq_updated} shows that for $\rho = 0$, the only attractor present in each $\mathbb{P}$ is an UA. By increasing $\rho$ to $0.05$, there is a dramatic change in $\mathbb{P}$ as the case of a poor reconstruction of $\pzl$ coexisting with UAs becomes the dominant scenario. For $\rho = 0.1$ we see the first occurrences of RCs that are able to reconstruct $\pzl$, albeit poorly, without any UAs in $\mathbb{P}$. This trend continues for increasing $\rho$ and when $\rho = 0.2$ we see the first instances of RCs that are able to provide a good reconstruction of $\pzl$ without any UAs present in $\mathbb{P}$. By increasing $\rho$ further we see this becomes the dominant scenario for $0.3 \lesssim \rho \lesssim 0.7$. Within this range of $\rho$ values, UAs are found for a decreasing number of $\textbf{M}$. 
For $0.7 \lesssim \rho \lesssim 1.1$, Fig.~\ref{fig:rho_vs_ScenarioFreq_updated} shows that there is a relatively equal likelihood that the RC will provide either a good or bad reconstruction of $\pzl$, and that it is much less likely for these reconstructions to coexist with UAs. However, there is a steady increase in the number of RCs that fail to provide any reconstruction of $\pzl$ within this range of $\rho$ values, thus, resulting in an increase in the UA only scenario. These trends continue for $\rho > 1.1$ with poor reconstructions of $\pzl$ becoming the dominant scenario.

\begin{figure*}
    \centering
    \includegraphics[width=0.95\linewidth]{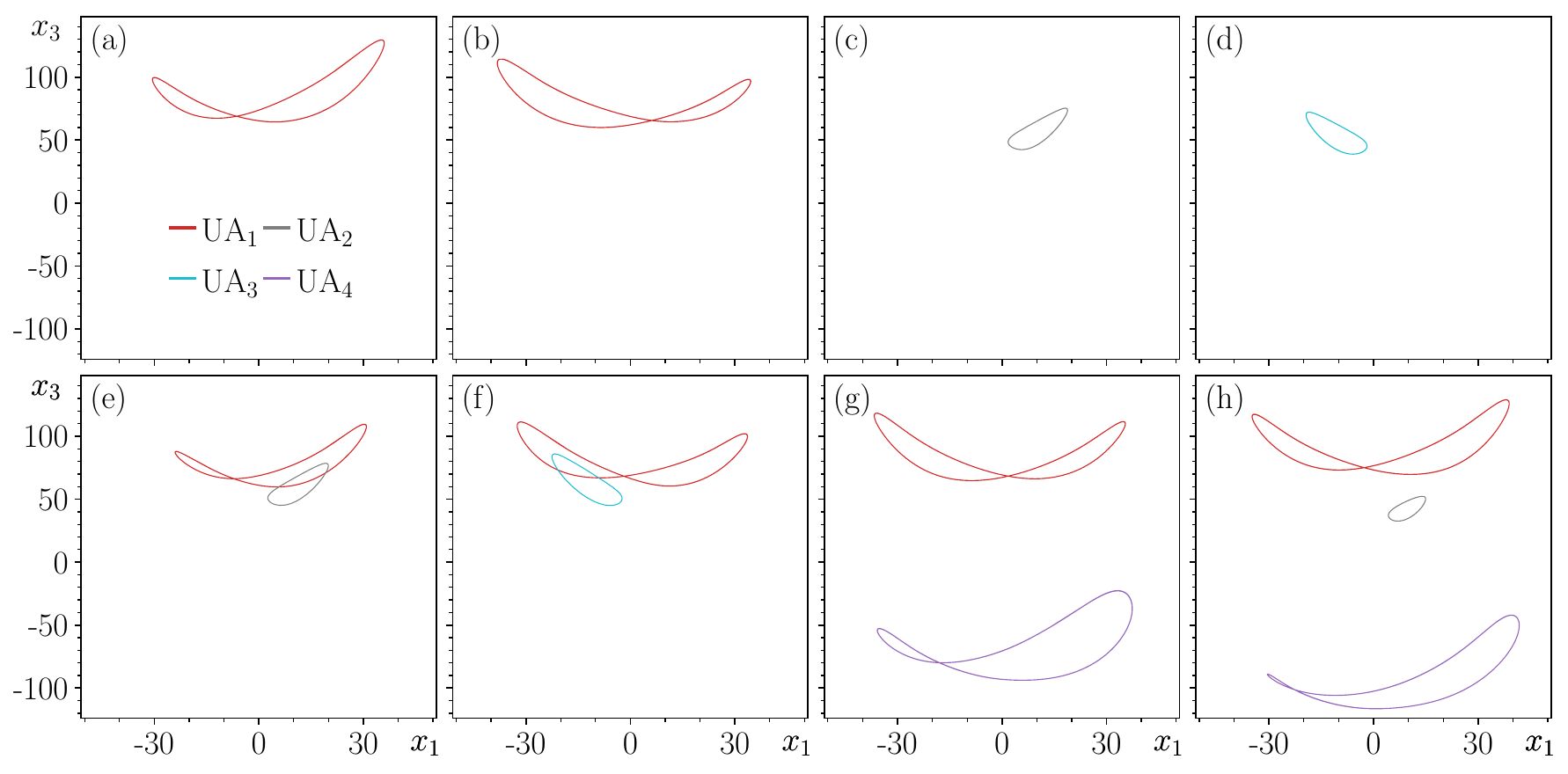}
    \caption{Examples of UAs that frequently appear for different random realisations of $\textbf{M}$ when $\rho=0.15$.}
    \label{fig:UA_examples_rho015_}
\end{figure*}

When generating the results shown in Fig.~\ref{fig:rho_vs_ScenarioFreq_updated} we encounter a variety of UAs. As an example, in Fig.~\ref{fig:UA_examples_rho015_} we illustrate the UAs that most frequently appear when setting $\rho = 0.15$. We label each of these UAs according to the region of $\mathbb{P}$ that they occupy. We omit any reconstructions of $\pzl$ from Fig.~\ref{fig:UA_examples_rho015_} in order to focus solely on the UAs that appear. Panels (a) and (b) illustrate the different behaviours of the most commonly found UA, denoted by UA$_{1}$. The difference between these illustrations is whether the maximum $x_{3}$ value occurs for either positive or negative $x_{1}$. The next most frequently found UAs are shown in panels (c) and (d), denoted by UA$_{2}$ and UA$_{3}$. We also find cases where UA$_{1}$ coexists with either UA$_{2}$ or UA$_{3}$, like in panels (e) and (f). Panel (g) shows that another UA, denoted by UA$_{4}$, can also coexist with UA$_{1}$. Less commonly, we find instances where three UAs coexist, like in panel (h).

\begin{figure*}
    \centering
    \includegraphics[width=\textwidth]{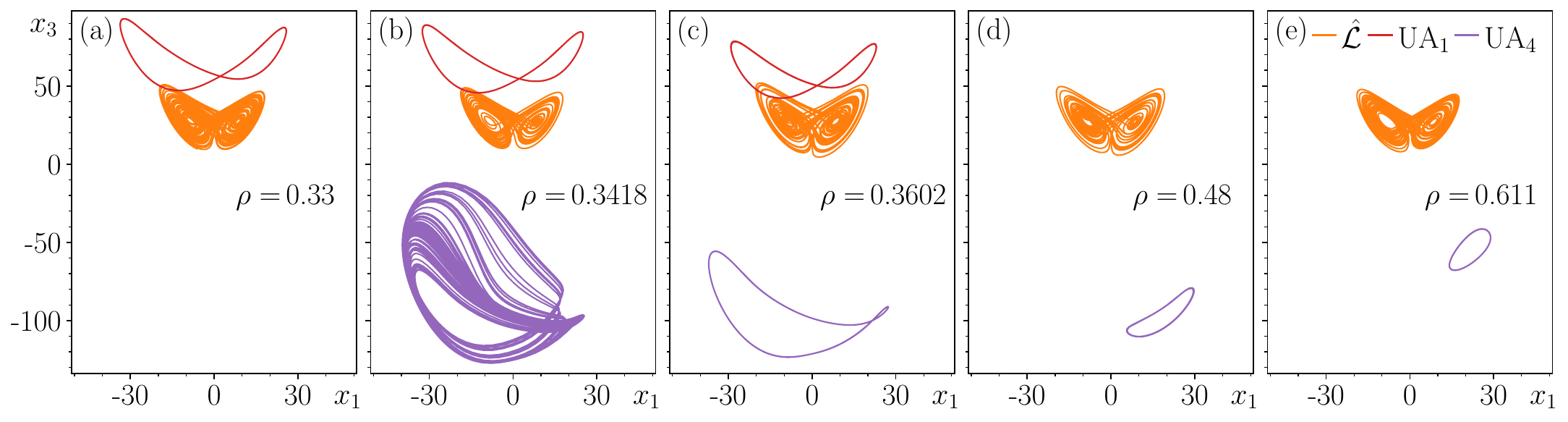}
    \caption{Example where UA$_{1}$ becomes unstable for increasing $\rho$ and UA$_{4}$ becomes chaotic and then unstable for decreasing $\rho$.}
    \label{fig:UA_examples_CA_}
\end{figure*}

Interestingly, for one of the random $\textbf{M}$ realisations in our ensemble, we show in Fig.~\ref{fig:UA_examples_CA_} that while UA$_{4}$ is a limit cycle for large $\rho$, panels (e), (d), (c), and (b) show that by decreasing $\rho$ this limit cycle becomes chaotic before becoming unstable in panel (a). On the other hand, for increasing $\rho$, panels (a)--(c) show that UA$_{1}$ moves closer and closer to $\hat{\pzl}$ until it becomes unstable and no longer directly observable. From this perspective, this behaviour indicates that for increasing $\rho$, UA$_{1}$ either collides with a nearby saddle or with the basin boundary of $\hat{\pzl}$.

\begin{figure}
    \centering
    \includegraphics[width=\linewidth]{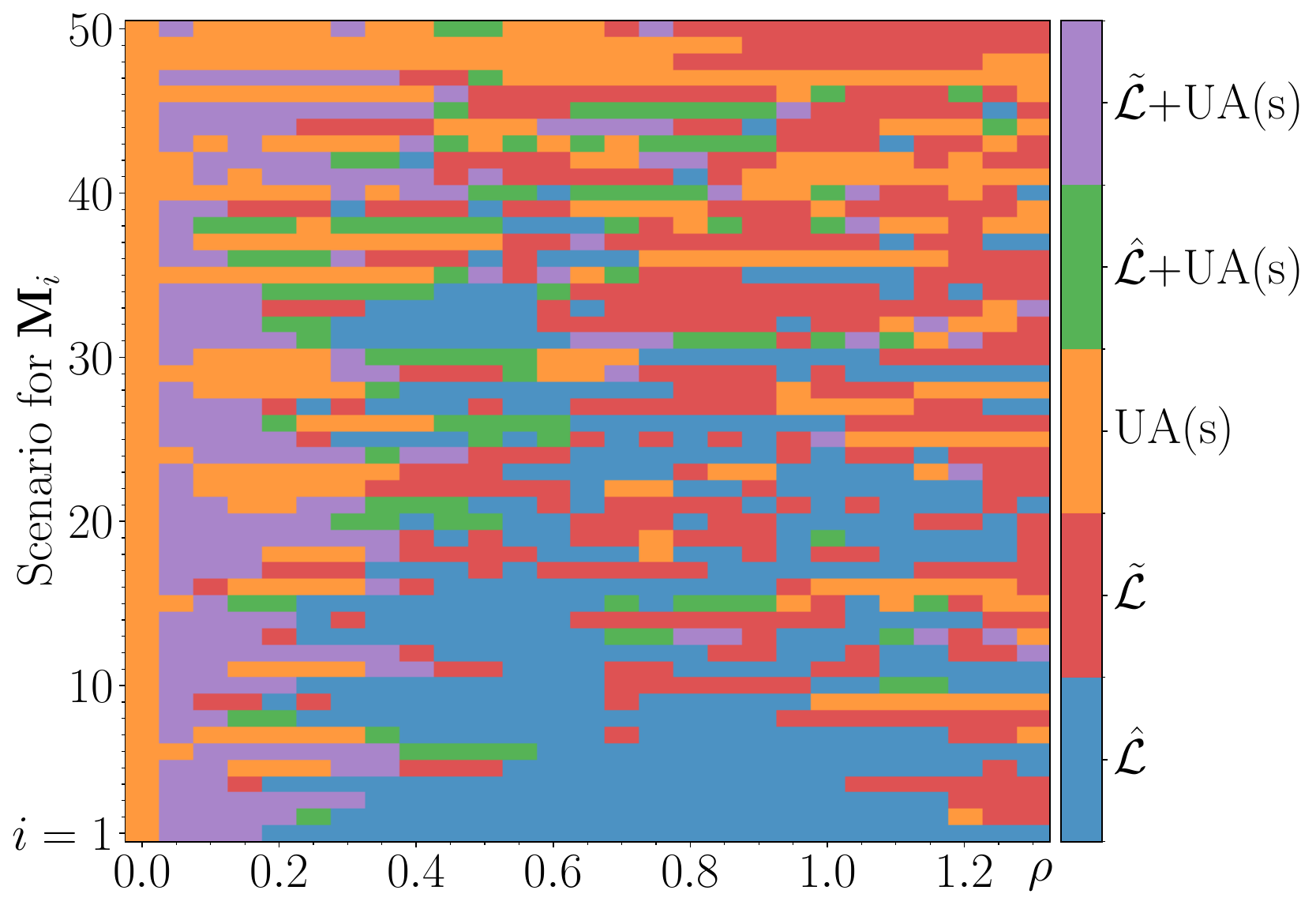}
    \caption{Resulting scenario for each of the 50 different $\textbf{M}$ used when training Eq.~\eqref{eq:ListenRes} to reconstruct $\pzl$ for a given $\rho$.}
    \label{fig:rho_vs_M_colorplot_updated}
\end{figure}

In Fig.~\ref{fig:rho_vs_M_colorplot_updated} we re-present the results of Fig.~\ref{fig:rho_vs_ScenarioFreq_updated} from the point of view of each $\textbf{M}$ in order to illustrate the changes that occur in each of the corresponding $\mathbb{P}$ when changing $\rho$. We label each $\textbf{M}$ by $\textbf{M}_{i}$, from $i=1$ to $50$, in order of the number of $\rho$ values that the corresponding closed-loop RC in Eq.~\eqref{eq:PredRes} provides a good reconstruction of $\pzl$, with $i=1$ being the most amount of $\rho$ values. Fig.~\ref{fig:rho_vs_M_colorplot_updated} provides us with additional information, such as: some random realisations of $\textbf{M}$ quickly overcome the influence of UAs, like $\textbf{M}_{1}$ and $\textbf{M}_{4}$, whereas others require relatively large $\rho$, like $\textbf{M}_{43}$ needing $\rho > 0.9$.

\subsection{UAs filling in the gaps when reconstruction of \texorpdfstring{$\pzl$}{TEXT} fails}

Fig.~\ref{fig:rho_vs_M_colorplot_updated} shows that when $\rho = 0.05$ the majority of the different $\textbf{M}$ matrices are able to provide a poor reconstruction of $\pzl$, denoted by $\tilde{\pzl}$, which coexists with at least one UA. 
In this subsection, we illustrate an interesting route in which the RC reconstructs $\pzl$ that involves period-doubling bifurcations. 
In Appendix~\ref{apx:LorenzTransientStable} we illustrate the most common way that the RC reconstructs $\pzl$ when increasing $\rho$ from $0$.
In Appendix~\ref{apx:AdditionalReconRoutes} we provide additional examples of interesting ways in which the RC reconstructs $\pzl$.

\subsubsection{Reconstruction of \texorpdfstring{$\pzl$}{TEXT} via period-doubling route to chaos}

\begin{figure*}
    \centering
    \includegraphics[width=0.8\linewidth]{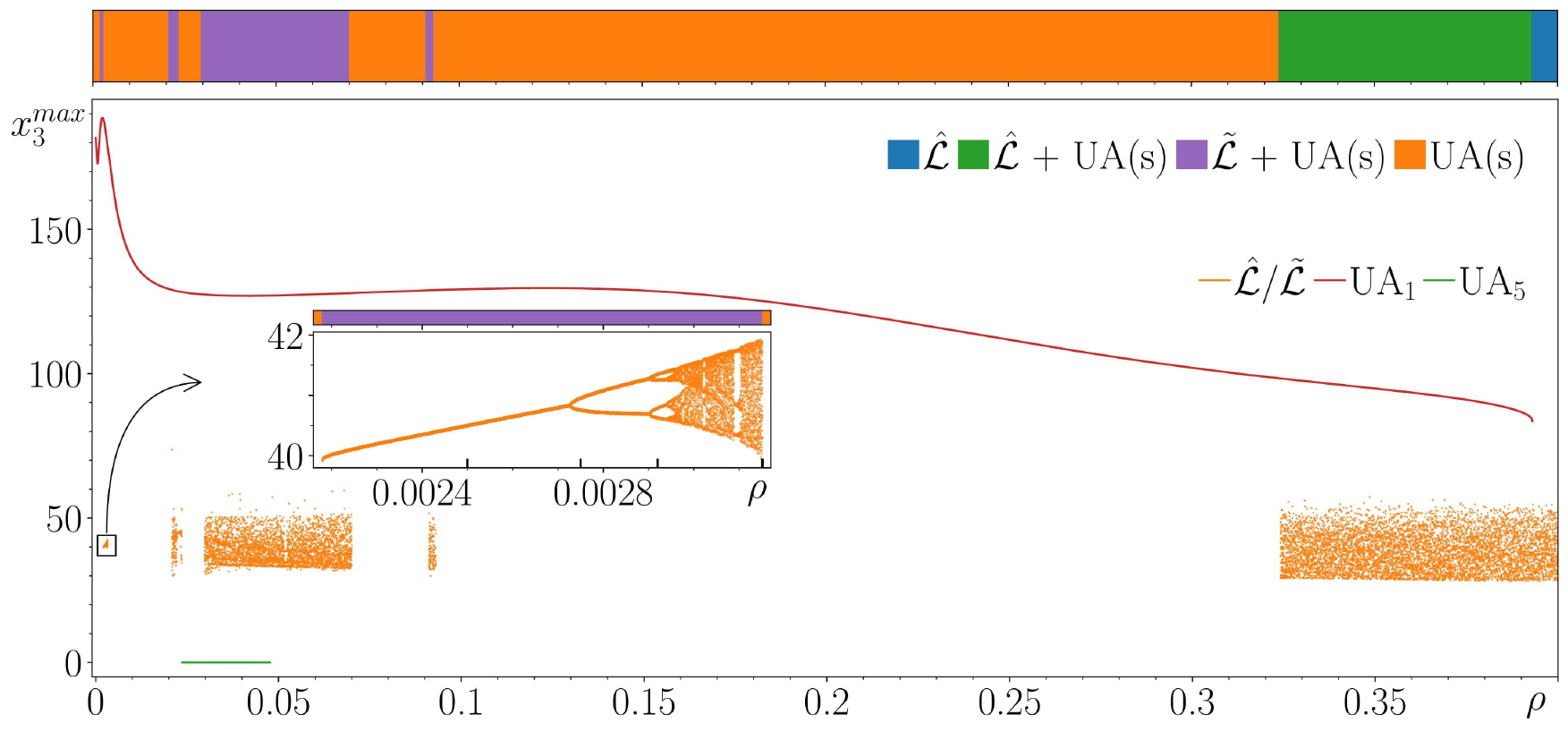}
    \caption{Period-doubling route for reconstruction of $\pzl$ for $\textbf{M}_{28}$. Plotting changes in local maxima of $x_{3}$ in $\mathbb{P}$ for $x_{1}>0$ for changes in $\rho$ for all attractors found for $\rho \in \left[0,0.4\right]$. Inset plot's inward ticks indicate $\rho$ values used to generate Fig.~\ref{fig:Lor_PD_snapshots_}. Coloured bars above bifurcation diagram indicate the algorithm's classification of the RC's output for a given $\rho$.}
    \label{fig:Lor_PDBifexample_rho_x3_}
\end{figure*}
\begin{figure*}
    \centering
    \includegraphics[width=\textwidth]{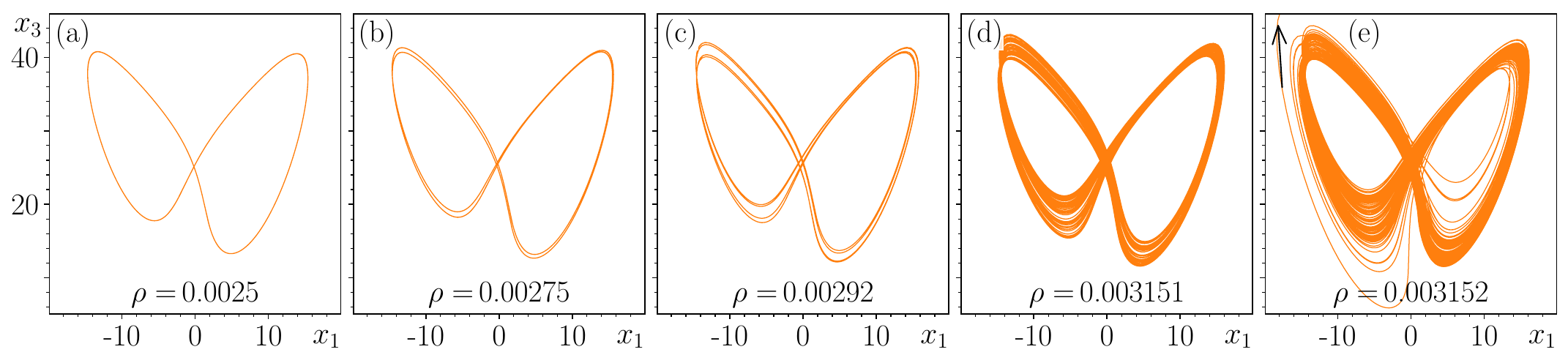}
    \caption{Snapshots of period-doubling route for reconstruction of $\pzl$ shown in Fig.~\ref{fig:Lor_PDBifexample_rho_x3_} at specified values of $\rho$ in panels (a)--(e).}
    \label{fig:Lor_PD_snapshots_}
\end{figure*}

In Fig.~\ref{fig:Lor_PDBifexample_rho_x3_}, we show that by increasing the parameter $\rho$ of $\textbf{M}_{28}$ from $0$, the first time the RC reconstructs a Lorenz-like attractor, it reconstructs $\pzl$ as a limit cycle. In the inset panel we see that this limit cycle undergoes a period-doubling route to chaos before becoming unstable. To provide further insight, in Fig.~\ref{fig:Lor_PD_snapshots_} we illustrate the dynamics of this limit cycle as it undergoes these bifurcations at different $\rho$ values. As indicated by the arrow in Fig.~\ref{fig:Lor_PD_snapshots_}~(e), after this reconstruction becomes unstable the state of the RC then approaches UA$_{1}$ as this is the only stable state in $\mathbb{P}$ for these $\rho$ values. Prior to becoming unstable, this reconstruction grows in its resemblance to $\pzl$, which is why our algorithm classifies this as a poor reconstruction of $\pzl$ rather than an UA. At larger $\rho$ values, we see this poor reconstruction of $\pzl$ reappear and disappear within different ranges of $\rho$ values for $\rho < 0.1$ in a similar way to the events shown in Fig.~\ref{fig:UAs_pred_rho0000_0009} in Appendix~\ref{apx:LorenzTransientStable}. Within this range of $\rho$ values we also find another UA, a fixed point that we denote by UA$_{5}$, which the state of the RC approaches when reconstruction fails for $0.025 \lesssim \rho \lesssim 0.03$. For $0.1 < \rho <0.325$ we find that the reconstruction exists as a transient and, depending on the value of $\rho$, this transient activity can last for significantly long durations of time. In order to determine whether the RC provides a stable reconstruction of $\pzl$ within this large range of $\rho$ values, we had to increase $t_{predict}$ to $30,000$. The RC achieves a good reconstruction of $\pzl$ for $\rho = 0.325$ and maintains this for much larger $\rho$ values.

Above all, Fig.~\ref{fig:Lor_PDBifexample_rho_x3_} demonstrates the benefit of performing this type of analysis of tracking the changes in these attractors. It provides a strong visual aid to show that when the RC fails to provide any reconstruction of $\pzl$, there is already an UA waiting for the RC's state to approach. Furthermore, this is not the only time we see period-doubling bifurcations playing a role in how the RC reconstructs $\pzl$. In Appendix~\ref{apx:AdditionalReconRoutes} we illustrate that period-doubling bifurcations of another Lorenz-like period-1 limit cycle and also a saddle-focus homoclinic orbit result in the reconstruction of $\pzl$. 
We would also like to mention that we find no sequence of events where $\hat{\pzl}$ comes into existence through the well-known sequence of events that unfold when varying the constant term in the $\dot{x}_{2}$ equation in Eq.~\eqref{eq:Lorenz}\cite{barrioShilnikov2012Lorenz}.

\section{Task (ii) results}\label{sec:Task2_results}

\subsection{Reconstructed attractors for \texorpdfstring{$b=\pm0.4, \pm0.3, \pm0.2$}{TEXT}}\label{ssec:Task2_reconstruction_b04_03_02}

\begin{figure}
    \centering
    \includegraphics[width=\linewidth]{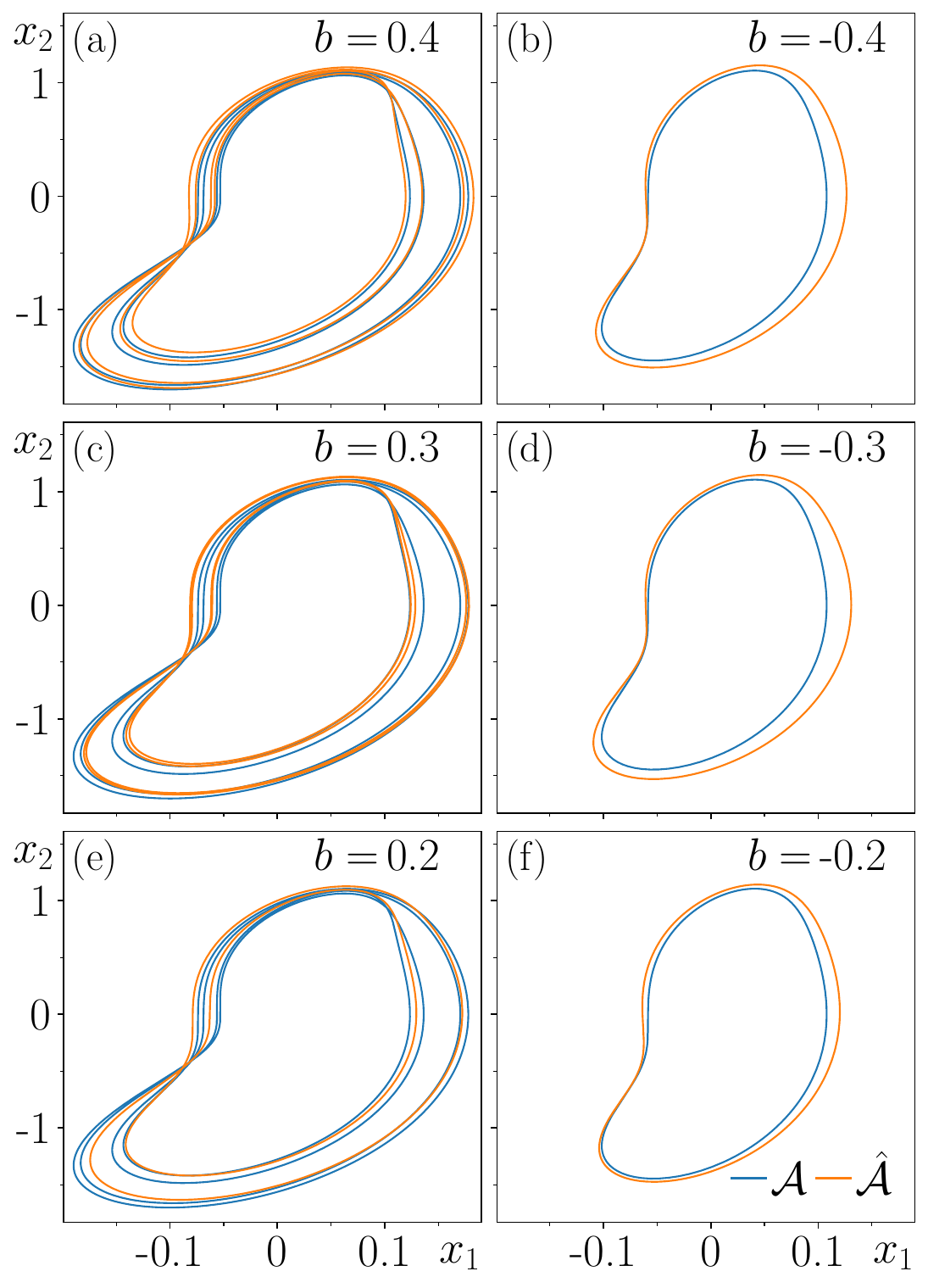}
    \caption{Reconstructed attractors (in orange) when parameter-aware RC is trained on task (ii) for $b=\pm0.4$ in (a)--(b), $b=\pm0.3$ in (c)--(d), and $b=\pm0.2$ in (e)--(f).}
    \label{fig:Task2_attrecon_b04_03_02_}
\end{figure}

We first show in Fig.~\ref{fig:Task2_attrecon_b04_03_02_} that the closed-loop parameter-aware RC (Eq.~\eqref{eq:ParamAwarePredRes}) achieves a reasonably good reconstruction of the period-1 and period-4 limit cycles for each choice of $b$ used during the training. The parameter $b$ is used to change the distance in parameter-space ($\boldsymbol{b}$-space) between these attractors. The other training parameters are specified in Sec.~\ref{sssec:Task2_description}.

\begin{figure}
    \centering
    \includegraphics[width=\linewidth]{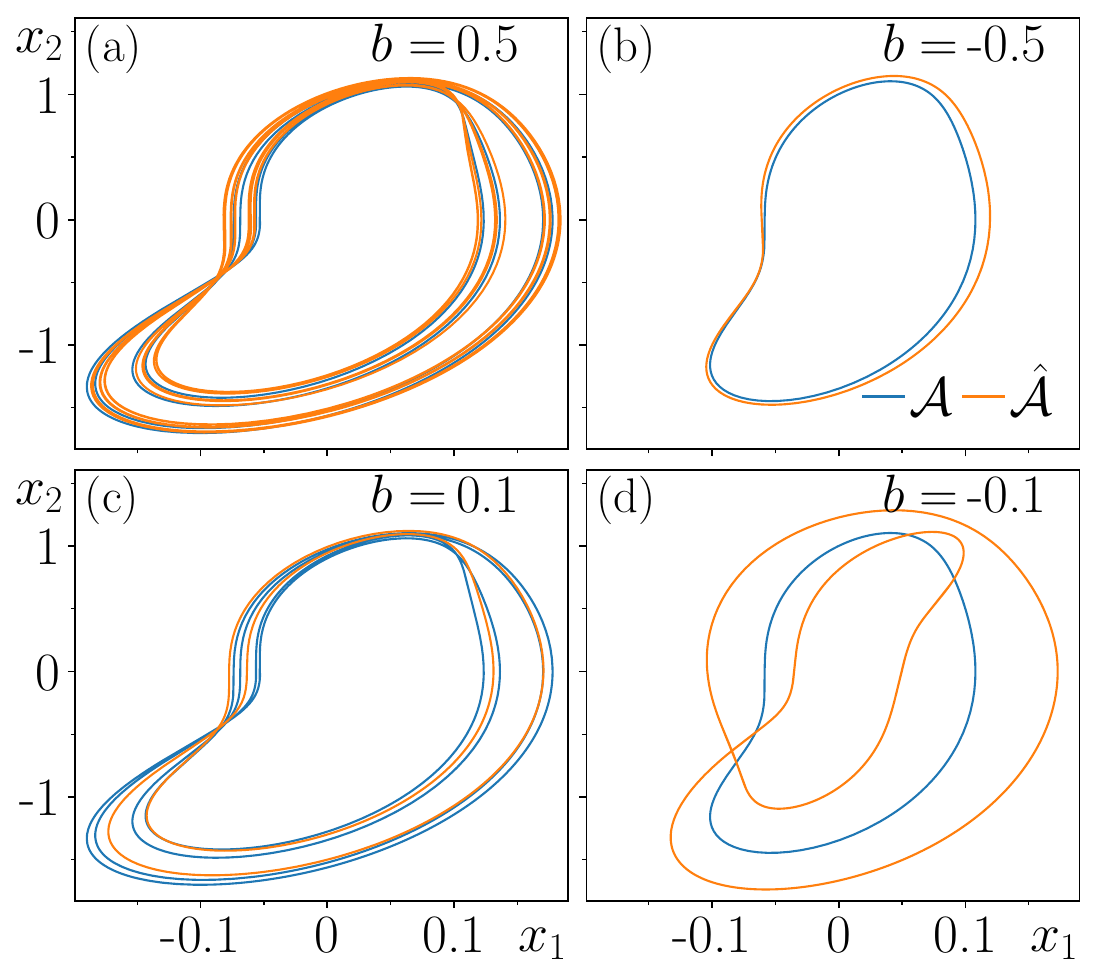}
    \caption{Reconstructed attractors (in orange) when parameter-aware RC is trained to reconstruct both limit cycles from task (ii) for $b=\pm0.5$ in (a)--(b) and $b=\pm0.1$ in (c)--(d).}
    \label{fig:Task2_attrecon_b05_01_}
\end{figure}

It is clear from Fig.~\ref{fig:Task2_attrecon_b04_03_02_} that there is an improvement in the RC's ability to reconstruct both limit cycles at their respective $b$ values when the trained limit cycles are further apart in $\boldsymbol{b}$-space. Panels~(a) and (c) show period-4 limit cycles, while panel~(e) reveals a similar appearing limit cycle that is actually period-2. Nonetheless, Figs.~\ref{fig:Task2_attrecon_b05_01_}~(a)--(b) shows there is a limit to this: if the limit cycles are too far apart, like for $b= \pm0.5$, then there is a decline in the RC's ability to reconstruct both limit cycles at their respective $b$ values. The RC fails to reconstruct the period-4 limit cycle and instead produces a chaotic attractor. Further, panels~(c)--(d) show what can happen when the limit cycles are too close together in parameter space. For $b= \pm0.1$ neither limit cycles are reconstructed properly and the RC instead produces period-2 limit cycles at both values of $b$.

\subsection{Analysis of GAs and UAs for \texorpdfstring{$b=\pm0.4, \pm0.3, \pm0.2$}{TEXT}}\label{ssec:Task2_gen_b04_03_02}

\begin{figure*}
    \centering
    \includegraphics[width=0.73\linewidth]{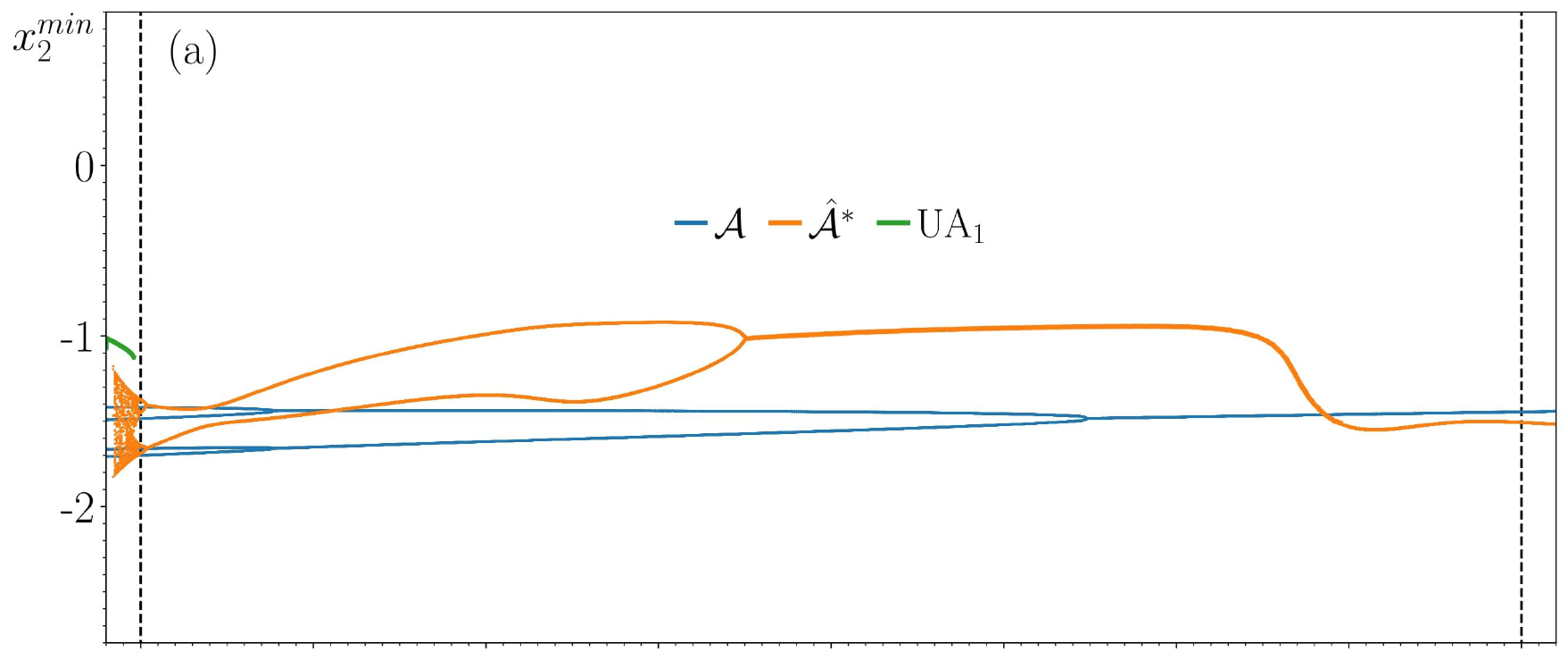}
    \includegraphics[width=0.73\linewidth]{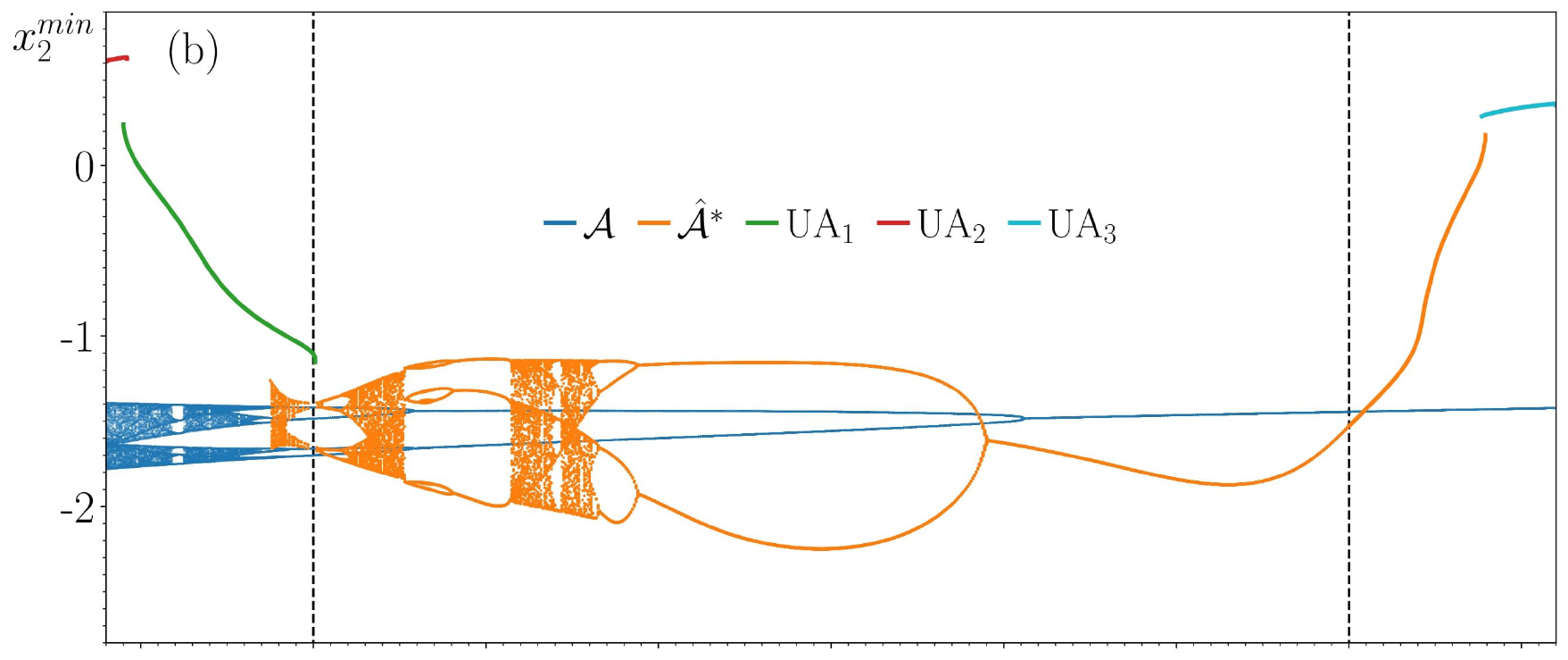}
    \includegraphics[width=0.73\linewidth]{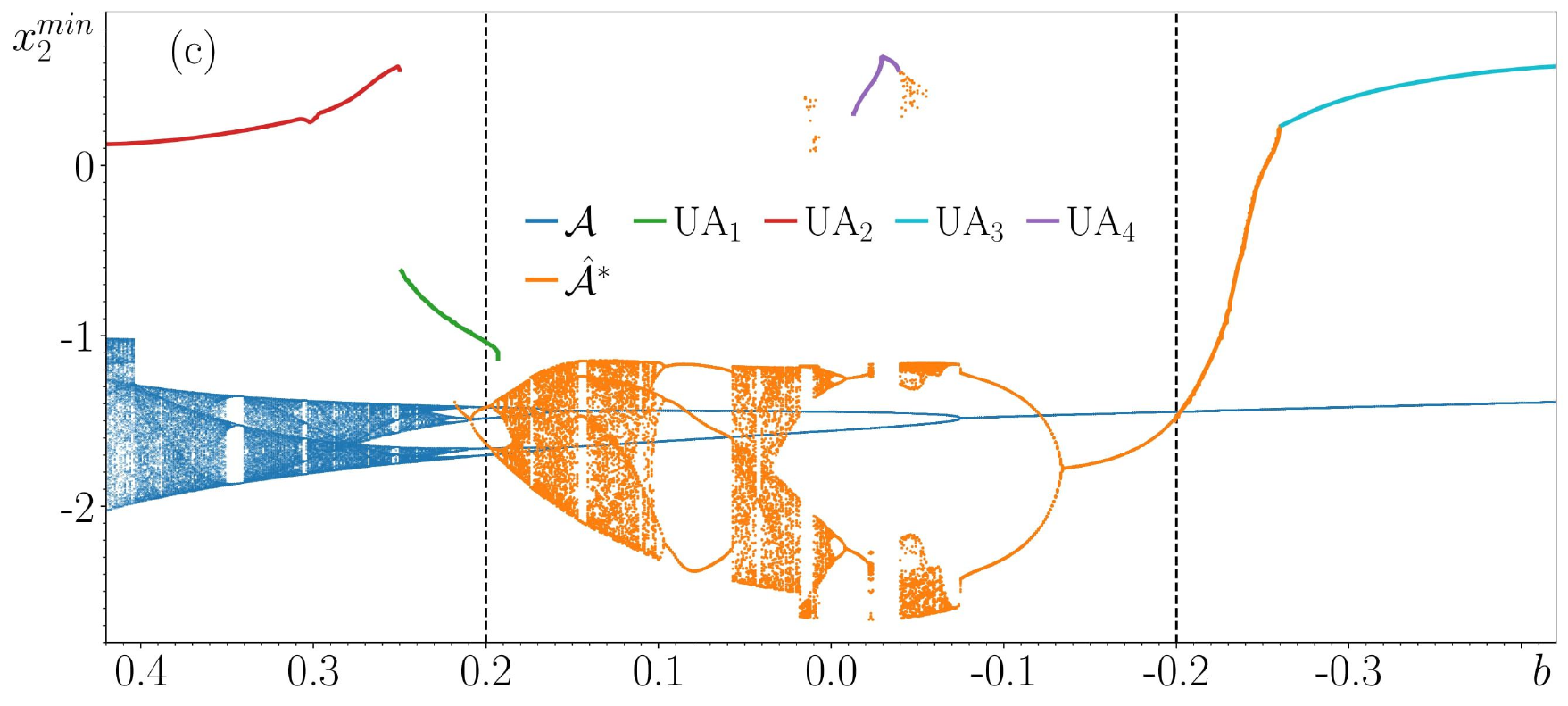}
    \caption{Resultant bifurcation diagram from training parameter-aware RC to perform task (ii) with $b=\pm0.4$ in (a), $b=\pm0.3$  in (b), or $b=\pm0.2$  in (c). Legend indicates the attractors that are tracked.}
    \label{fig:Task2_bbifplot040302_compare_}
\end{figure*}

\begin{figure*}
    \centering
    \includegraphics[width=0.85\linewidth]{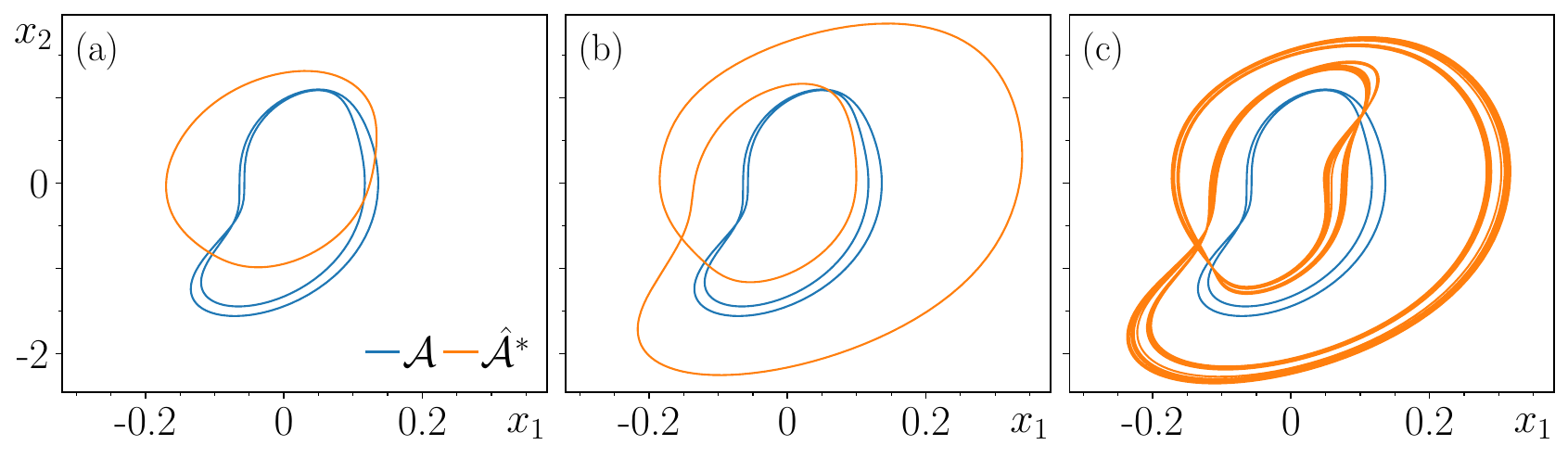}
    \caption{Examples of GAs (in orange) for $\boldsymbol{b}=0$ when trained for $b=\pm0.4$ in (a), $b=\pm0.3$ in (b), and $b=\pm0.2$ in (c).}
    \label{fig:attrecon_x1x2_b04_03_02_b0_}
\end{figure*}

In each panel of Fig.~\ref{fig:Task2_bbifplot040302_compare_} we show how the closed-loop parameter-aware RC in Eq.~\eqref{eq:ParamAwarePredRes} generates the dynamics between and beyond the reconstructed limit cycles shown in Fig.~\ref{fig:Task2_attrecon_b04_03_02_} by tracking the changes in the dynamics of all attractors present in the RC's prediction state space, $\mathbb{P}$, for $b\in \left[ 0.42, -0.42 \right]$. In Fig.~\ref{fig:attrecon_x1x2_b04_03_02_b0_} we show the attractors generated by the RC in each panel of Fig.~\ref{fig:Task2_bbifplot040302_compare_} for $b=0$. We discuss each of the corresponding panels of Figs.~\ref{fig:Task2_bbifplot040302_compare_} and \ref{fig:attrecon_x1x2_b04_03_02_b0_} for values of $b$ under the following headings. 

\subsubsection{\texorpdfstring{Case $b=\pm0.4$}{TEXT}}\label{ssec:Task2_gen_b04}

We first discuss the results for $b = \pm 0.4$ shown in Fig.~\ref{fig:Task2_bbifplot040302_compare_}~(a). The most prominent feature here is the continuous transition between the trained period-1 and period-4 limit cycles. This continuous transition is consistent with the period-doubling bifurcation that occurs in Eq.~\eqref{eq:sprott}. However, the RC's period-doubling bifurcation does not occur nearby the corresponding bifurcation in Eq.~\eqref{eq:sprott}. Also, the generated period-2 limit cycle bears little resemblance to the limit cycle from Eq.~\eqref{eq:sprott}, as evidenced by the difference in the local minima of the $x_{2}$ variable. Fig.~\ref{fig:attrecon_x1x2_b04_03_02_b0_}~(a) further illustrates this difference between the RC's generated dynamics for $b=0$, generating a period-1 limit cycle instead of a period-2 limit cycle. 

For $b>0.4$, Fig.~\ref{fig:Task2_bbifplot040302_compare_}~(a) shows that the period-4 limit cycle undergoes a period-doubling route to chaos, and that the resultant chaotic attractor becomes unstable for $b \approx 0.415$. To fill in the gaps, Eq.~\eqref{eq:ParamAwarePredRes} generates an UA, that we denote by UA$_{1}$ (which has nothing to do with the UA$_{1}$ studied in task (i), nor do any of the UAs mentioned in this section). From further inspection, UA$_{1}$ is a period-1 limit cycle and its dynamics are consistent with the dynamics of the limit cycle that coexists with the limit cycle from Eq.~\eqref{eq:sprott} that we use for training, which to some extent indicates that this RC learns properties of Eq.~\eqref{eq:sprott} without being explicitly trained to. This phenomenon is described in \textcite{rohm2021unseen}.

\subsubsection{\texorpdfstring{Case $b=\pm0.3$}{TEXT}}\label{ssec:Task2_gen_b03}

The most prominent feature in Fig.~\ref{fig:Task2_bbifplot040302_compare_}~(b) is that the continuous transition between the trained period-1 and period-4 limit cycles is not consistent with the period-doubling bifurcation that occurs in Eq.~\eqref{eq:sprott}. Instead of undergoing a period-doubling bifurcation from the reconstructed period-1 and period-4 limit cycles to a GA in the form of a period-2 limit cycle, the RC generates infinite sequences of period-doubling and halving bifurcations to and from generated chaotic attractors. Fig.~\ref{fig:attrecon_x1x2_b04_03_02_b0_}~(b) illustrates the period-2 limit cycle generated by the RC for $b=0$, which occupies a much larger region of the state space than the corresponding dynamics from Eq.~\eqref{eq:sprott}. 

Similar to Fig.~\ref{fig:Task2_bbifplot040302_compare_}~(a), Fig.~\ref{fig:Task2_bbifplot040302_compare_}~(b) shows that for $b>0.3$, the reconstructed period-4 limit cycle undergoes a period-doubling route to chaos, and the resultant chaotic attractor becomes unstable for $b \approx 0.325$. UA$_{1}$ fills in the resulting gap. However, the changes in the dynamics of UA$_{1}$ can no longer be tracked for $b \gtrsim 0.41$. The RC fills in the gap with another UA, which for this range of $b$ values is a fixed point that we denote by UA$_{2}$. Note, there is a small range of $b$ values where UA$_{1}$ and UA$_{2}$ coexist, specifically, $b \in [0.408,0.410]$. 

Fig.~\ref{fig:Task2_bbifplot040302_compare_}~(b) shows that GAs of the period-1 limit cycle can no longer be tracked for $b<-0.379$. The RC fills in the resulting gap with UA$_{3}$, which is a fixed point for the $b$ values shown here and coexists with the GAs of the reconstructed period-1 limit cycle for $b \in [-0.377,-0.379]$.

\subsubsection{\texorpdfstring{Case $b=\pm0.2$}{TEXT}}\label{ssec:Task2_gen_b02}

Fig.~\ref{fig:Task2_bbifplot040302_compare_}~(c) shows there is no continuous transition between the reconstructed period-1 and period-4 limit cycles when training for $b= \pm 0.2$. The attractors that are generated from the reconstructed period-1 and period-4 limit cycles become chaotic and subsequently unstable for decreasing/increasing $b$. No GAs exist for $b \in \left( -0.024, -0.04 \right)$ and UA$_{4}$, which is a fixed point, fills in this gap. UA$_{4}$ coexists with the generated dynamics from the reconstructed period-4 limit cycle for $b \in \left[ -0.013, -0.024 \right]$, and coexists with the generated dynamics from the reconstructed period-1 limit cycle for a much smaller range of $b$ values. Fig.~\ref{fig:attrecon_x1x2_b04_03_02_b0_}~(c) illustrates the dynamics of the chaotic attractor shown in Fig.~\ref{fig:Task2_bbifplot040302_compare_}~(c) for $b=0$. 

Similar to Fig.~\ref{fig:Task2_bbifplot040302_compare_}~(b),  Fig.~\ref{fig:Task2_bbifplot040302_compare_}~(c) shows that UA$_{1}$ and UA$_{2}$ fill in the gaps from $b > 0.2$, and likewise UA$_{3}$ fills in the gaps for $b < -0.2$. The small dip in the branch of points describing the changes in the dynamics of UA$_{2}$ corresponds to a bifurcation from a fixed point to a limit cycle for decreasing $b$.

\subsection{Does adding more data help?}\label{ssec:Task2.5_gen_b02}

\begin{figure*}
    \centering
    \includegraphics[width=0.85\linewidth]{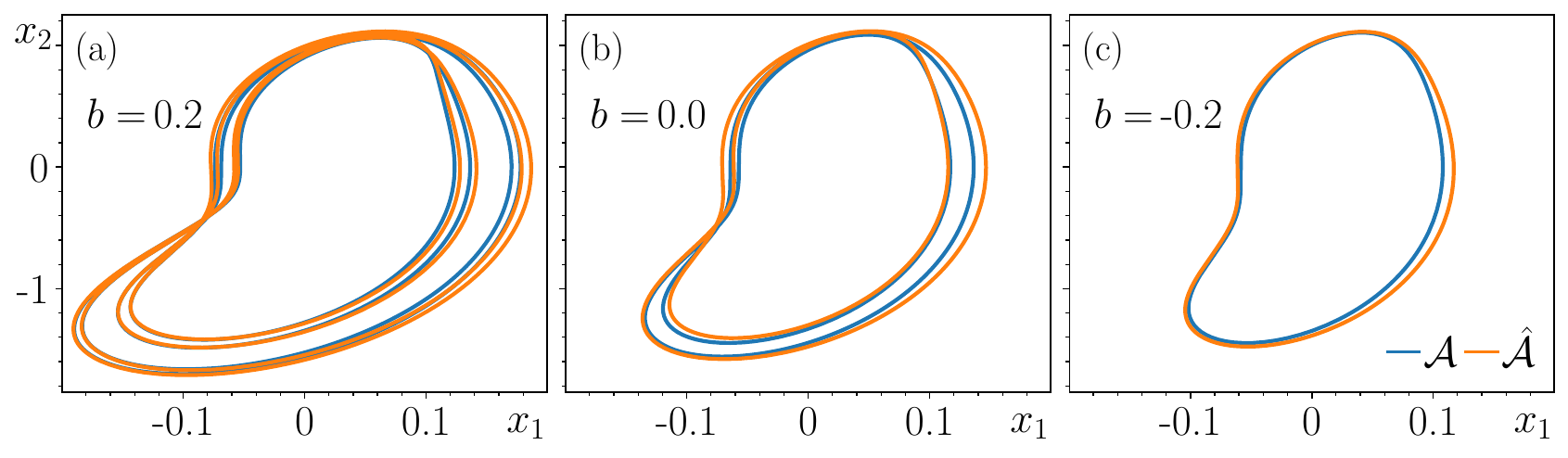}
    \caption{Examples of reconstructed attractors (in orange) when trained for $b=0.2,0.0,-0.2$.}
    \label{fig:attrecon_x1x2_b02_train3_}
\end{figure*}
\begin{figure*}
    \centering
    \includegraphics[width=0.73\linewidth]{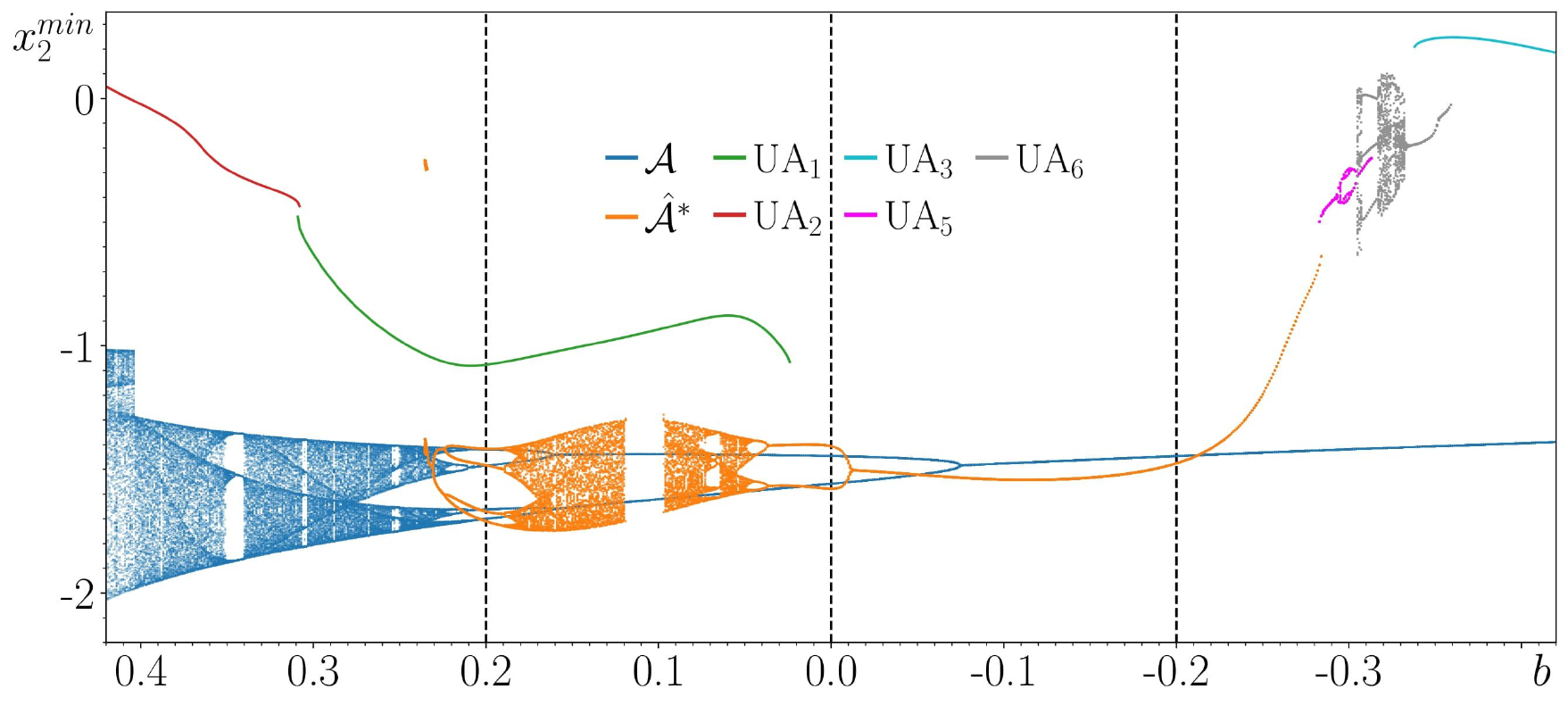}
    \caption{Resultant bifurcation diagram from training parameter-aware RC to reconstruct period-1 limit cycle at $b=-0.2$, period-2 limit cycle at $b=0.0$, and period-4 limit cycle at $b=0.2$. Legend indicates the attractors that are tracked.}
    \label{fig:bbifplot02_train3}
\end{figure*}

A natural question arises from the results presented so far on task (ii), what would happen if the parameter-aware RC was also trained to reconstruct the period-2 limit cycle from Eq.~\eqref{eq:sprott}? Does this improve the RCs ability to generate transitions, i.e. period-doubling bifurcations, between the period-1, 2, and 4 limit cycles? Would the resultant dynamics be more consistent with the dynamics of Eq.~\eqref{eq:sprott}? To shed light on these questions, in this subsection we illustrate our results from training the parameter-aware RC to reconstruct the period-1 and period-4 limit cycles for $\boldsymbol{b}=(\pm 0.2) \boldsymbol{1}^{T}$ as before and (i) to also reconstruct the intermediary period-2 limit cycle for $\boldsymbol{b}=(0) \boldsymbol{1}^{T}$, and (ii) reconstruct the period-2 limit cycle for $\boldsymbol{b}=(0) \boldsymbol{1}^{T}$ as well as reconstruct the corresponding intermediary period-1 and period-2 limit cycles for $\boldsymbol{b}=(\pm 0.1) \boldsymbol{1}^{T}$. We do so using the extended training method mentioned in the second last paragraph of Sec.~\ref{sec:BifReconTrain}.

\begin{figure*}
    \centering
    \includegraphics[width=\linewidth]{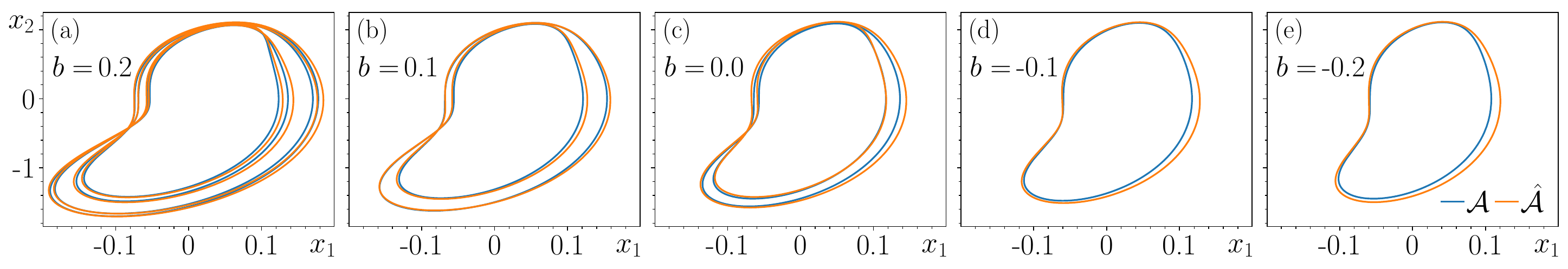}
    \caption{Examples of reconstructed attractors (in orange) when trained for $b=0.2,0.1,0.0,-0.1,-0.2$.}
    \label{fig:attrecon_x1x2_b02_train5_}
\end{figure*}

\begin{figure*}
    \centering
    \includegraphics[width=0.73\linewidth]{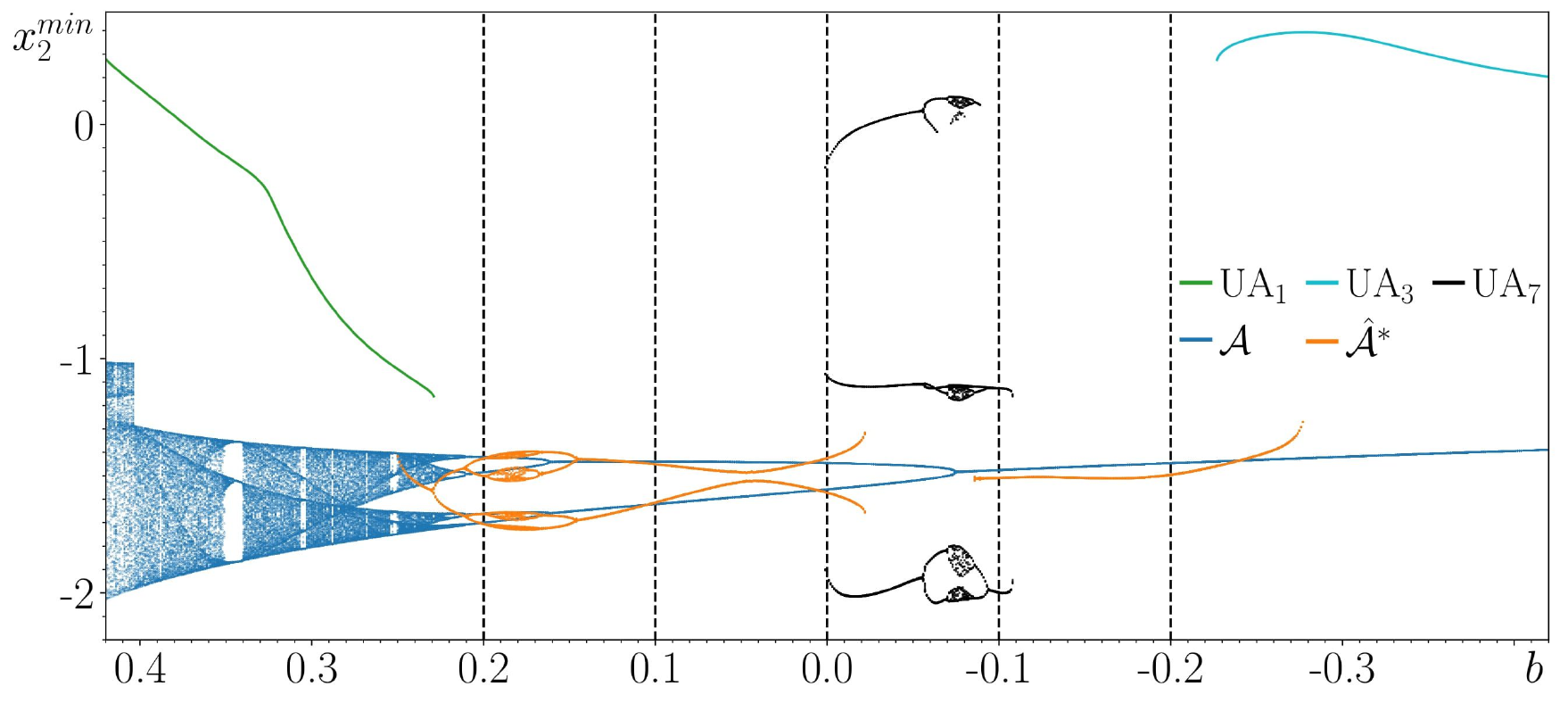}
    \caption{Resultant bifurcation diagram from training parameter-aware RC to reconstruct period-1 limit cycle at $b=-0.2$ and $-0.1$, period-2 limit cycle at $b=0.0$ and $0.1$, and period-4 limit cycle at $b=0.2$. Legend indicates the attractors that are tracked.}
    \label{fig:bbifplot02_train5}
\end{figure*}

In both cases we found that the closed-loop parameter-aware RC was unable to successfully reconstruct the limit cycles at their respective $b$ values with the set of RC training parameters specified in Sec.~\ref{sssec:Task2_description} that were used so far in task (ii). Thus, to follow the same rationale as before, we find a new set of RC training parameters that provides reasonably good reconstruction of all limit cycles for the values of $b$ that were used during the training. This allows us to concentrate our discussion on how well the RC fills in gaps in its memory for memories (i.e. attractors) that are stored correctly. For the first case of three limit cycles we achieved this by setting $\rho = 1.4$, $\sigma = 1.3$, without changing the other parameters. For the second case of five limit cycles we achieved this by setting $\rho = 1.3$, $\sigma = 1.0$, $\beta = 0.1$, and the other parameters unchanged. We illustrate in Figs.~\ref{fig:attrecon_x1x2_b02_train3_} and \ref{fig:attrecon_x1x2_b02_train5_} that for these new choices of training parameters, the closed-loop parameter-aware RC achieves a good reconstruction of each limit cycle. This is further exemplified in Figs.~\ref{fig:bbifplot02_train3} and \ref{fig:bbifplot02_train5} by the overlap between the blue and orange data points at the values of $b$ used during the training, as indicated by the vertical dashed black lines.

In Fig.~\ref{fig:bbifplot02_train3} we examine how the closed-loop parameter-aware RC in Eq.~\eqref{eq:ParamAwarePredRes} generates the dynamics between and beyond the three reconstructed limit cycles shown in Fig.~\ref{fig:attrecon_x1x2_b02_train3_}. More specifically, Fig.~\ref{fig:bbifplot02_train3} shows the result of tracking the changes in the dynamics of all attractors generated by Eq.~\eqref{eq:ParamAwarePredRes} for $b\in \left[ 0.42, -0.42 \right]$. While Fig.~\ref{fig:attrecon_x1x2_b02_train3_} shows that the RC achieves a good reconstruction of each limit cycle at their respective $b$ values, the same cannot be said for how the RC generates the dynamics between these limit cycles. While the transition from the period-1 to period-2 limit cycle is relatively consistent with the dynamics of Eq.~\eqref{eq:sprott}, there is no continuous transition between the period-2 and period-4 limit cycles. The period-2 and period-4 limit cycles both become chaotic when tracking the changes in the dynamics of these attractors for increasing $b$ from $0$ and decreasing $b$ from $0.2$, respectively. The chaotic attractor that is born from the reconstructed attractors reaches a point where it can no longer be tracked and thus becomes unstable for $b \in \left(0.12, 0.06\right)$. Like before, UA$_{1}$ fills this gap, while UA$_{1}$ and UA$_{2}$ fill in the gap for $b > 0.235$. A different scenario from Fig.~\ref{fig:Task2_bbifplot040302_compare_}~(c) emerges in Fig.~\ref{fig:bbifplot02_train3} for $b < -0.2$. In addition to UA$_{3}$, there are two further UAs required to fill in the gaps: the limit cycle UA$_{5}$ and the mostly chaotic dynamics of UA$_{6}$.


In Fig.~\ref{fig:bbifplot02_train5} we apply the same procedure to examine how Eq.~\eqref{eq:ParamAwarePredRes} generates the dynamics between and beyond the five reconstructed limit cycles shown in Fig.~\ref{fig:attrecon_x1x2_b02_train5_}. Similarly, while Fig.~\ref{fig:attrecon_x1x2_b02_train5_} shows that the RC achieves a good reconstruction of each limit cycle at their respective $b$ values, there is a slightly better agreement between the attractors that the RC generates and the dynamics of Eq.~\eqref{eq:sprott} than in Fig.~\ref{fig:attrecon_x1x2_b02_train3_}, aside for the branch of GAs between $b=0.2$ and $0.1$ which shows the RC's dynamics become chaotic nearby $b=0.18$. Furthermore, there is still no continuous transition between all reconstructed attractors. The periodic and chaotic dynamics of UA$_{7}$ fills the gap in this case, coexisting with the corresponding GAs for a small range of $b$ values.
Just as in many previous cases, UA$_{1}$ and UA$_{3}$ fill in the gaps for $b > 0.25$ and $b < -0.27$.

In summary, the main message from Figs.~\ref{fig:bbifplot02_train3} and \ref{fig:bbifplot02_train5} is that even if you try to account for gaps in the system's memory then complete confabulations may still appear. However, by taking greater measures to account for these gaps, complete confabulations may become less prevalent and there may be a greater resemblance between the generated dynamics and the dynamics of the system where the training data comes from. 

In the next section, we investigate how the RC fills in gaps in its memory when there exists no natural or simple way to account for gaps in its memory.


\section{Task (iii) results}\label{sec:Task3_results}

\subsection{Reconstructed attractors for \texorpdfstring{$b=\pm0.4, \pm0.3, \pm0.2, \pm0.1$}{TEXT}}\label{ssec:Task3_reconstruction_b04_03_02_01}

\begin{figure}
    \centering
    \includegraphics[width=\linewidth]{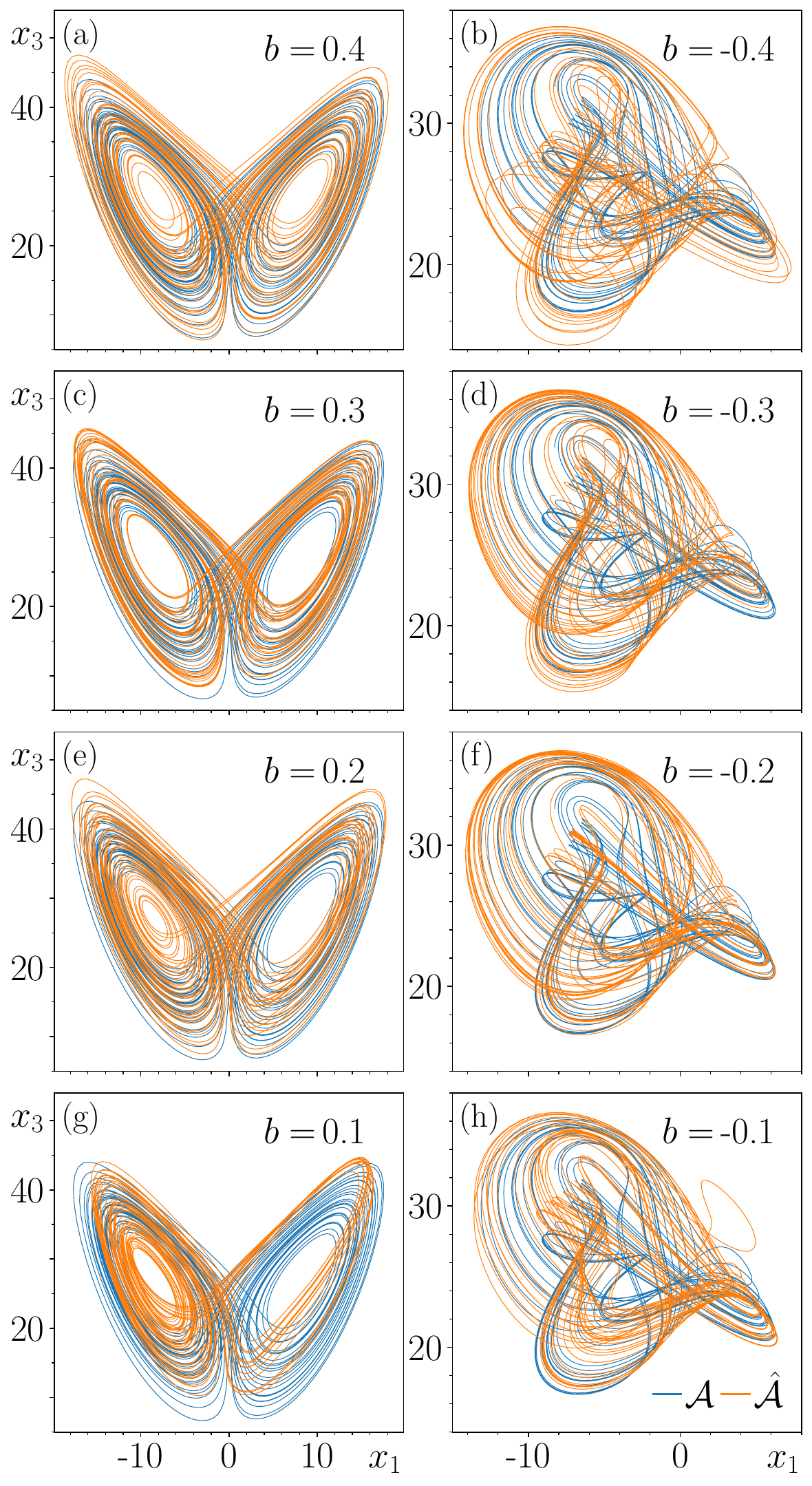}
    \caption{Reconstructed attractors (in orange) when parameter-aware RC is trained on task (iii) for $b=\pm0.4$ in (a)--(b), $b=\pm0.3$ in (c)--(d), $b=\pm0.2$ in (e)--(f), and $b=\pm0.1$ in (g)--(h).}
    \label{fig:Task3_attrecon_b04_03_02_01_}
\end{figure}

We first show in Fig.~\ref{fig:Task3_attrecon_b04_03_02_01_} that the closed-loop parameter-aware RC in Eq.~\eqref{eq:ParamAwarePredRes} achieves a reasonably good reconstruction of the Lorenz ($\pzl$) and Halvorsen ($\pzh$) chaotic attractors for $b=\pm0.2$, $\pm0.3$, and $\pm0.4$, and relatively poor reconstruction for $b=\pm0.1$, in particular for $\pzl$. Similar to Fig.~\ref{fig:Task2_attrecon_b04_03_02_}, Fig.~\ref{fig:Task3_attrecon_b04_03_02_01_} shows that as the attractors are moved further apart in $b$-space then reconstruction of both attractors improves. However, there is also a limit to this as well, comparing Fig.~\ref{fig:Task3_attrecon_b04_03_02_01_}~(b) to (c) we see that reconstruction of $\pzh$ begins to worsen as $b$ is increased from $\pm0.3$ to $\pm0.4$.

\begin{figure*}
    \centering
    \includegraphics[width=0.73\linewidth]{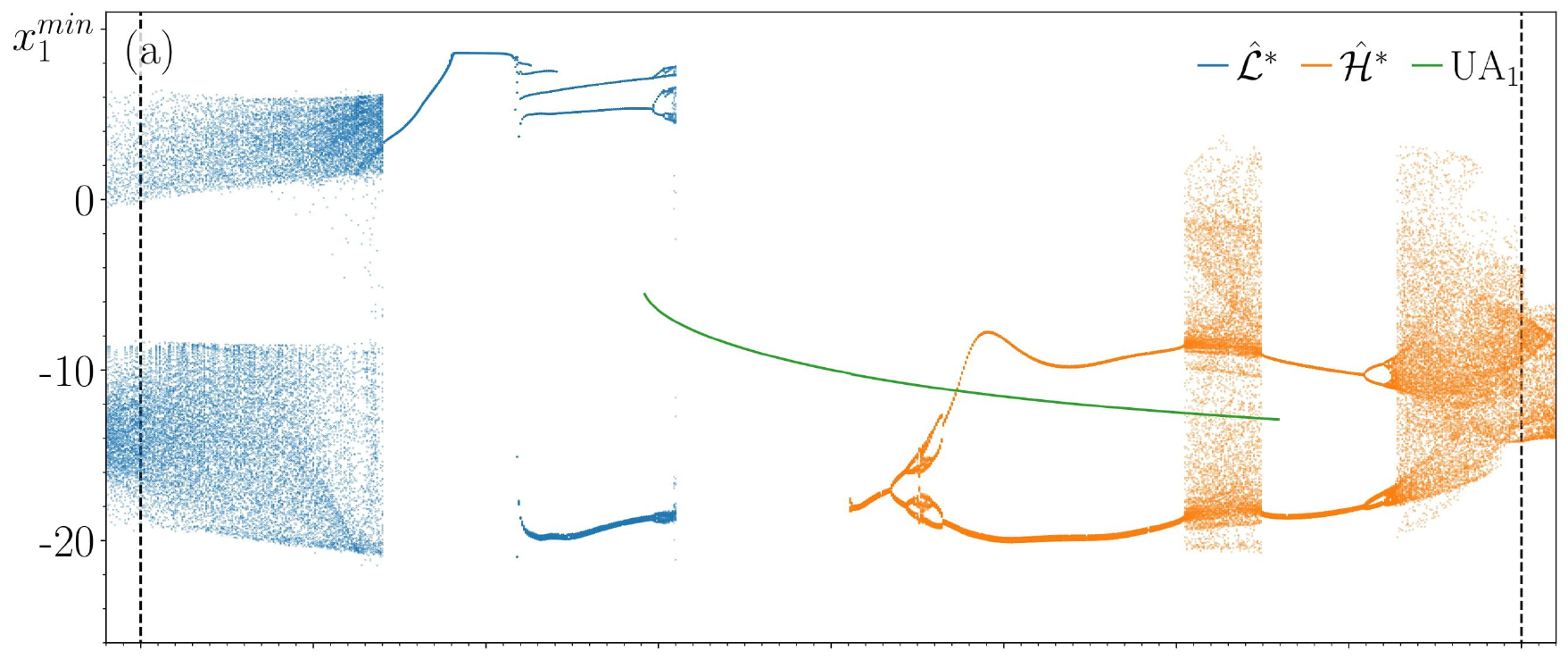}
    \includegraphics[width=0.73\linewidth]{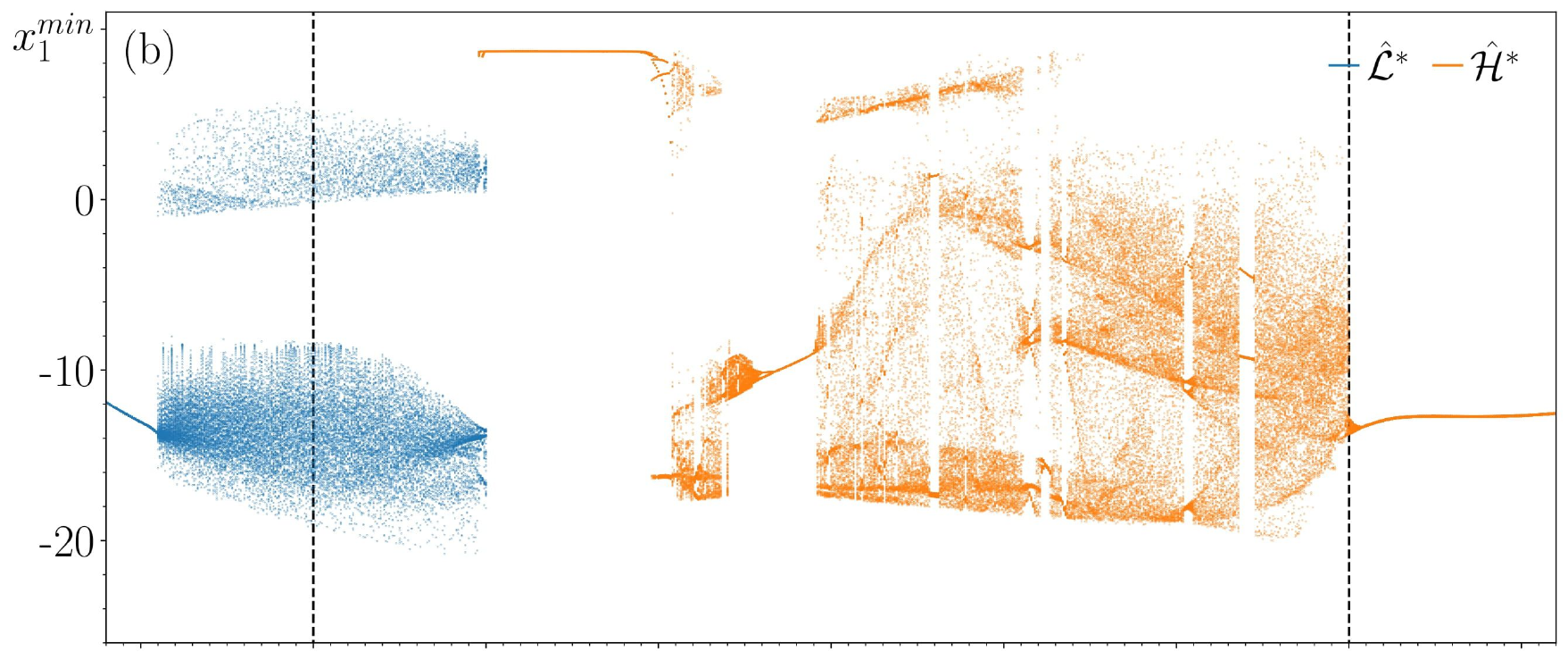}
    \includegraphics[width=0.73\linewidth]{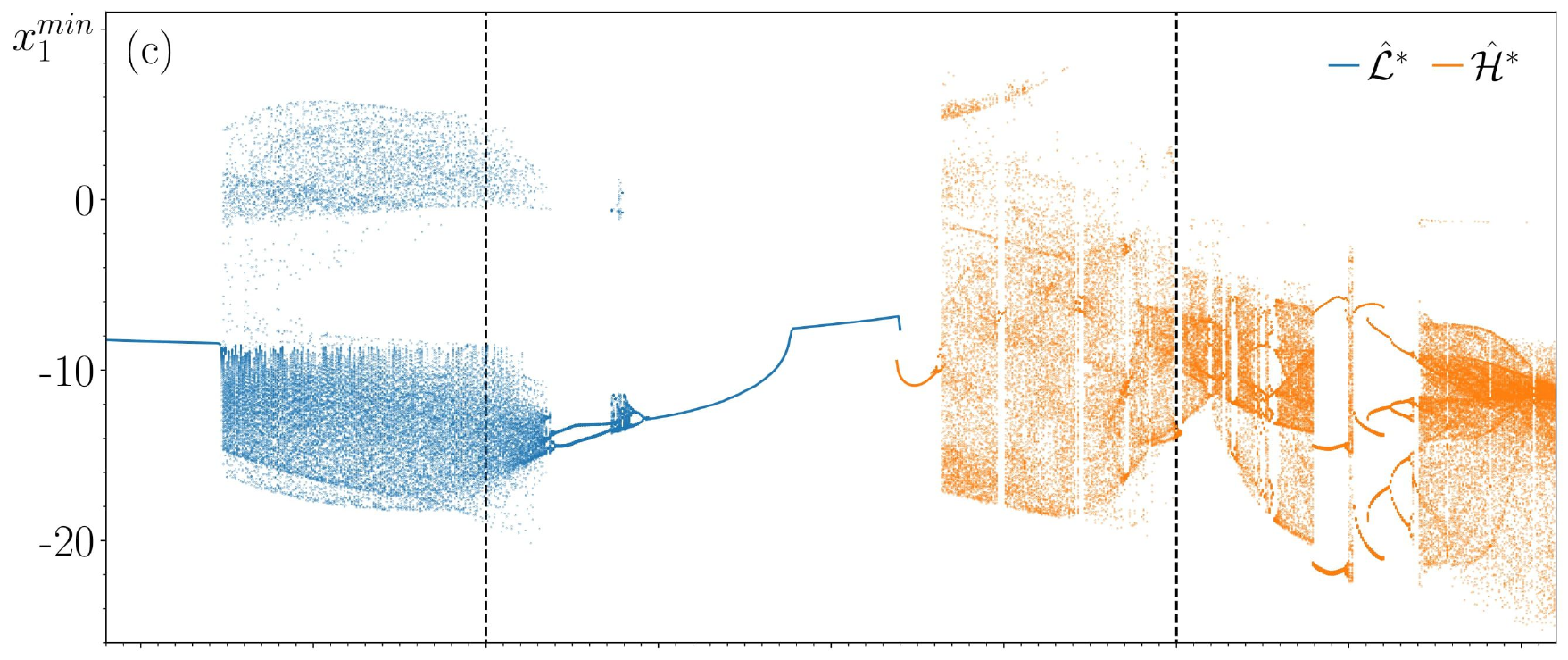}
    \includegraphics[width=0.73\linewidth]{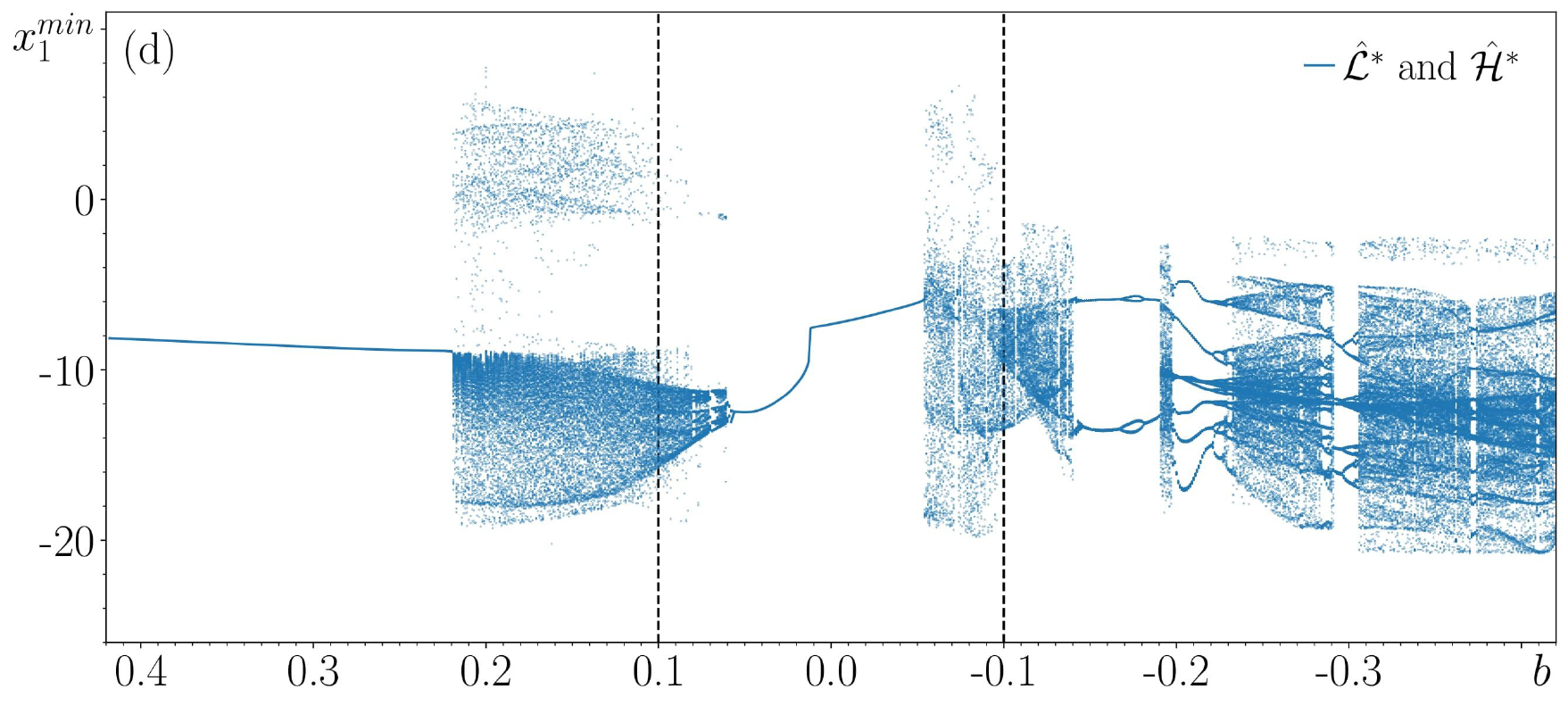}
    \caption{Resultant bifurcation diagram from training parameter-aware RC to perform task (iii) with $b=\pm0.4$, $\pm0.3$, $\pm0.2$, and $\pm0.1$. Legend indicates the attractors that are tracked.}
    \label{fig:Task3_bbifplot04030201_compare_}
\end{figure*}

\subsection{Analysis of GAs and UAs for \texorpdfstring{$b=\pm0.4, \pm0.3, \pm0.2, \pm0.1$}{TEXT}}\label{ssec:Task3_gen_b04_03_02_01}

In each panel of Fig.~\ref{fig:Task3_bbifplot04030201_compare_} we show the result of tracking the changes in the dynamics of all attractors generated by Eq.~\eqref{eq:ParamAwarePredRes} for $b\in \left[ 0.42, -0.42 \right]$. We do this in order to examine how the closed-loop parameter-aware RC in Eq.~\eqref{eq:ParamAwarePredRes} generates the dynamics between and beyond the reconstructed chaotic attractors, $\pzl$ and $\pzh$, shown in Fig.~\ref{fig:Task3_attrecon_b04_03_02_01_}. We discuss each of the corresponding panels of Figs.~\ref{fig:Task3_attrecon_b04_03_02_01_} and \ref{fig:Task3_bbifplot04030201_compare_} for values of $b$ under the following headings.

\subsubsection{\texorpdfstring{Case $b=\pm0.4$}{TEXT}}\label{ssec:Task3_gen_b04}

We first discuss the results for $b = \pm 0.4$ shown in Fig.~\ref{fig:Task3_bbifplot04030201_compare_}~(a). The most prominent feature here is that there is no continuous transition between the two attractors. To fill the gap between the attractors, the RC creates UA$_{1}$ (of no connection to the previously mentioned UA$_{1}$s), a fixed point which coexists for a wide range of $b$ values with the attractors that are generated from the reconstruction of $\pzl$ and $\pzh$, namely $\hat{\pzl}^{*}$ and $\hat{\pzh}^{*}$. Despite Figs.~\ref{fig:Task3_attrecon_b04_03_02_01_}~(a) and (b) illustrating that the RC achieves a reasonably good reconstruction of $\pzl$ and $\pzh$, Fig.~\ref{fig:Task3_bbifplot04030201_compare_}~(a) shows that both attractors alternate between chaotic and periodic dynamics at several $b$ values when tracking the changes in the dynamics of $\hat{\pzl}^{*}$ and $\hat{\pzh}^{*}$ for decreasing and increasing $b$. Both attractors can no longer be tracked for $b \in \left[ 0.09, -0.01 \right]$ and UA$_{1}$ fill this gap. 

\subsubsection{\texorpdfstring{Case $b=\pm0.3$}{TEXT} and \texorpdfstring{$b=\pm0.2$}{TEXT}}\label{ssec:Task3_gen_b03_02}

In Figs.~\ref{fig:Task3_bbifplot04030201_compare_}~(b) and (c) we see that while both $\pzl^{*}$ and $\pzh^{*}$ become unstable at certain $b$ values, there is no need for the RC to create an UA as there is a small region of bistability between the GAs. Furthermore, despite that Figs.~\ref{fig:Task3_attrecon_b04_03_02_01_}~(c)--(d) and (e)--(f) show that the RC achieves a reasonably good reconstruction of $\pzl$ and $\pzh$ at the respective $b$ values, Figs.~\ref{fig:Task3_bbifplot04030201_compare_}~(b) and (c) show that when tracking the changes in the dynamics of $\hat{\pzl}^{*}$ and $\hat{\pzh}^{*}$ for decreasing and increasing $b$, the dynamics of both GAs also alternate between chaotic and periodic dynamics at several $b$ values like in Fig.~\ref{fig:Task3_bbifplot04030201_compare_}~(a). 

In Fig.~\ref{fig:Task3_bbifplot04030201_compare_}~(b), we see that when tracking the changes in the dynamics of both attractors from $b = \pm 0.3$, $\hat{\pzl}^{*}$ can no longer be tracked for $b < 0.2$ and $\hat{\pzh}^{*}$ can no longer be tracked for $b > 0.21$. Thus, the region of bistability between $\hat{\pzl}^{*}$ and $\hat{\pzh}^{*}$ exists for $0.2 \le b \le 0.21$. When tracking $\hat{\pzl}^{*}$ for decreasing $b$ from $b=0.3$, we find that these GAs are deformations of the original reconstructed attractor $\hat{\pzl}$. On the other hand, when tracking $\hat{\pzh}^{*}$ for increasing $b$ from $b=-0.3$ we find that the dynamics of the GAs share much less similarities to $\pzh$. When tracking $\hat{\pzl}^{*}$ for increasing $b$ from $b=0.3$ and $\hat{\pzh}^{*}$ for decreasing $b$ from $b=-0.3$, both GAs become period-1 limit cycles.

In Fig.~\ref{fig:Task3_bbifplot04030201_compare_}~(c), we see that when tracking the changes in the dynamics of both attractors from $b=\pm0.2$, $\hat{\pzl}^{*}$ can no longer be tracked from $b < -0.04$ and $\hat{\pzh}^{*}$ can no longer be tracked for $b > -0.038$. Thus, the region of bistability between $\hat{\pzl}^{*}$ and $\hat{\pzh}^{*}$ exists for $-0.04 \le b \le -0.038$. When tracking $\hat{\pzl}^{*}$ for decreasing $b$ from $b=0.2$ we find that the deformations of the original reconstructed attractor $\hat{\pzl}$ occur for a much shorter range of $b$ values. Instead, $\hat{\pzl}^{*}$ is a limit cycle for $b \in \left( 0.11, 0.02 \right)$ and fixed point for $b \in \left( 0.02, -0.038 \right]$. 
When tracking $\hat{\pzh}^{*}$ for increasing $b$ from $b=-0.2$, $\hat{\pzh}^{*}$ becomes a limit cycle for $b \in \left[ -0.06, -0.038 \right]$. When tracking $\hat{\pzl}^{*}$ for increasing $b$ from $b=0.2$, at $b=0.35$ this GAs undergoes a continuous transition from a chaotic attractor to a fixed point. 
When tracking the changes in $\hat{\pzh}^{*}$ for decreasing $b$ from $b=-0.2$, the GAs alternate between regions of chaotic and periodic behaviour. Note that the different branches of $x_{1}^{min}$ nearby $b=-0.315$ do not correspond to a region of bistability between two limit cycles, instead there is one limit cycle that has two local minima at different $b$ values.


\subsubsection{\texorpdfstring{Case $b=\pm0.1$}{TEXT}}\label{ssec:Task3_gen_b01}

The most interesting result from task (iii) is illustrated in Fig.~\ref{fig:Task3_bbifplot04030201_compare_}~(d) which shows that the closed-loop parameter-aware RC in Eq.~\eqref{eq:ParamAwarePredRes} creates a continuous transition between attractors from two different systems. Thus, using the terminology developed in Sec.~\ref{ssec:UA_define}, since there is a continuous transition then the generated dynamics between both attractors aligns with our definition of $\hat{\pza}_{i,j}^{*}$. However, we find it more convenient in the plot legend to refer to the information presented here as `$\hat{\pzl}^{*}$ and $\hat{\pzh}^{*}$', since there is a continuous transition between all GAs in Fig.~\ref{fig:Task3_bbifplot04030201_compare_}~(d). 

When tracking the changes in $\hat{\pzl}^{*}$ for decreasing $b$ from $b=0.1$, we find that $\hat{\pzl}^{*}$ undergoes a series of bifurcations to become a period-1 limit cycle near $b=0.06$ and undergoes a further bifurcation to a fixed point near $b=0.01$. 
When tracking the changes in $\hat{\pzh}^{*}$ for increasing $b$ from $b=-0.1$, $\hat{\pzh}^{*}$ undergoes what appears to be a more abrupt transition from chaotic behaviour to the same fixed point that is found when tracking $\hat{\pzl}^{*}$. 
We study these transitions in greater detail in Appendix~\ref{ssec:CloserExam_L_H}. 

Similar to Fig.~\ref{fig:Task3_bbifplot04030201_compare_}~(c), Fig.~\ref{fig:Task3_bbifplot04030201_compare_}~(d) shows that when tracking $\hat{\pzl}^{*}$ for increasing $b$ from $b=0.1$, at $b=0.22$ this GA undergoes a continuous transition from a chaotic attractor to a fixed point. 
Likewise, when tracking $\hat{\pzh}^{*}$ for decreasing $b$ from $b=-0.1$, the GAs alternate between regions of chaotic and periodic dynamics. 
While this continuous transition emerges for decreasing $b$, Figs.~\ref{fig:Task3_attrecon_b04_03_02_01_} shows that the RC's reconstruction of $\pzl$ and $\pzh$ worsens for decreasing $b$. From further experiments not shown in this paper, we find this trend continues for decreasing $b$.

\section{\label{sec:DisConc}Conclusion}

The primary motivation of this paper was to develop a foundational understanding of the reasons why and how ANNs confabulate and generate false information to fill in gaps in their memory. Our broader objective was to use these insights to inform the study of analogous processes in other learning systems, including the brain. As a first step, we study how confabulations arise in a relatively simple ANN design, a dynamical system in the form of a RC. 

Previous work mentioned in Sec.~\ref{ssec:PreviousWork} has shown that when RCs are trained to reconstruct the dynamics of a given attractor, the RC can construct an attractor it was not trained to reconstruct, a so-called `untrained attractor' (UA). 
Building on the work of \textcite{LuBassett20_switching_learning} and \textcite{kong2024memory}, we interpret the phenomenon of UAs through a memory-theoretic lens, associating trained attractors with encoded memories and UAs with \textit{complete confabulations}.

We designed three tasks to explore the reasons why and how complete confabulations in the form of UAs emerge. In task (i) (Sec.~\ref{sec:Task1_results}) we studied the UAs that appear when a RC is trained to store and recall an individual attractor/memory. In tasks (ii) and (iii), (Secs.~\ref{sec:Task2_results} and \ref{sec:Task3_results}), we studied the UAs that appear when a RC is trained to generalise between related and unrelated attractors/memories. 
The common result across all tasks is that the appearance of UAs consistently coincided with failure modes, either in recall or generalisation, appearing whenever the RC failed to encode or sensibly generalise beyond its training data. This aligns with the known association between the emergence of confabulations and a decline in the brain's ability to store and/or recall memories\cite{arts2017korsakoff,el2017ConfabAlz}.

Several additional key findings emerge from our experiments. First, confabulations are not necessarily a system error but are an intrinsic feature of how learning systems behave. For instance, the results presented in Sec.~\ref{sec:Task1_results} show that confabulations can appear as soon as the system is trained and disappear as the system improves its performance, in these cases the UAs become unstable and thus their presence is no longer directly visible. 
Second, our results challenge the assumption that confabulations can be mitigated solely through training the system with additional data. As shown in Sec.~\ref{sec:Task2_results}, increasing the amount of training data does not guarantee that the system provides a more sensible generalisation between learned states and UAs may still appear. Furthermore, while it has been shown that a reasonably good agreement can be achieved between the bifurcation diagrams of parameter-aware RCs and certain systems, like in  \textcite{kim2021GlobalLocal}, \textcite{kong21PredCritTrans}, \textcite{kong2023digitaltwins}, \textcite{koglmayr2024tipping} and even when using experimental data from systems like in \textcite{panahiYCLai2024AMOC}, here we show that the amount of attractors used and the spacing between the parameter values where the attractors are stored (the choices of $\boldsymbol{b}$) are two factors that influence the level of agreement. 
Third, our findings showcase the explanatory power of applying concepts from dynamical systems theory to better understand how machines learn by analysing the stability and bifurcation structure of trained and untrained states in ANNs. From this analysis we show that UAs may become stable or unstable through collisions with transients or reconstructions of attractors as certain training parameters are varied. To augment this analysis, our next step is to extend current continuation software, such as AUTO \cite{AUTO_Doedel}, in order to study the unstable dynamics of trained, untrained, and generated states and how these interact with their stable counterparts.

While the results presented in this paper are based on `toy examples' (synthetic data, not real-world data), in the future we intend to conduct follow-up research to the work presented here using real-world data. 
However, we would like to take the opportunity to highlight investigations of confabulation in learning systems that include a real-world component. 
For instance, in terms of confabulation in physical RCs, in \textcite{ryo25_multifunctionalphysicalRC} the authors conduct an in-depth study of UAs that appear in a multifunctional physical RC in the form of a soft tensegrity robot. More specifically, the authors trained the robot to reconstruct attractors that correspond to two different motions (rapid crawl and shake). However, after further investigation of the robots state space, the authors found that the robot has five additional UAs present in its state space which correspond to completely different motions that it was not trained to exhibit. 
In terms of confabulation in large language models (LLMs), in \textcite{castilho2025_LLM_Irish} the authors present an interesting study of how different versions of ChatGPT confabulate when used to translate sentences from the English to Irish language (Gaeilge). Interestingly, the authors identify that these LLMs confabulate in a variety of ways by replacing unknown words with made-up words that, for instance, do or do not follow plausible translation rules, or disregard conventions of the Irish alphabet (which does not include certain letters from the English alphabet like `v').
In terms of confabulation in the human brain, in \textcite{el2017_AD_confab} the authors present a study of complete confabulations in patients with Alzheimer's disease, concluding that complete confabulations are correlated with patients that have spatial and temporal disorientation and that momentary confabulations occur more often than complete confabulations which, interestingly and perhaps coincidentally, corresponds with our findings in Fig.~\ref{fig:rho_vs_ScenarioFreq_updated}. 






\subsection{\label{sec:ConfabConj} Complete confabulation conjecture}

Our results suggest that complete confabulations in the form of UAs are not a random occurrence but a by-product of the system being \textit{bounded} (i.e. a system whose state cannot approach infinity). In the RC we study, this boundedness property is induced by the $\tanh{\left( \cdot \right)}$ activation function as it confines the elements of the RC's state to be valued between $-1$ and $1$. Thus, with both the decay term and bounded activation function the state of the system $\boldsymbol{r}(t)$ cannot grow to infinity and the state space of the RC is bounded. 
As a result, when reconstruction fails, the system does not become unstable with its state tending towards infinity. Instead, since the system cannot approach infinity, it must approach an alterative attractor, which is the UA. 
In contrast, when reconstruction fails using a so-called `next-generation RC', the system may tend to infinity because there is no guarantee of boundedness in the design of a next-generation RC. Examples of this occurrence can be found in Fig.~1 of \textcite{Flynn22_LimitsMF}.

Here we suggest that the boundedness property of a system may serve as a unifying principle for the intrinsic existence of complete confabulations in artificial and biological learning systems:
\begin{itemize}
    \item \textbf{\textit{The complete confabulation conjecture}:} any bounded learning system has the capacity to create complete confabulations.
\end{itemize} 



Many learning systems share this property of boundedness. 
The state of most deep neural networks will be bounded due to the use of common activation functions that are either bounded or asymptotically linear \cite{szandala2020_DLactivfn}, such that any bounded input will generate a bounded output.
Similarly, biological neurons in the brain are subject to physiological constraints and can only respond to stimuli within a finite dynamical range \cite{miller2019_NeuralFiringRate}. 
In our future research we will test the validity of the complete confabulation conjecture across different types of RCs, other machine learning paradigms, different tasks, and levels of biological realism. 




\begin{acknowledgments}
This work is a continuation of work that started during A.F.'s time as a PhD student which was funded by the Irish Research Council Enterprise Partnership Scheme (Grant No. EPSPG/2017/301) and features elements of work from J.O'H.'s M.Sc. thesis. 
We would like to thank Andreas Amann for his influential conversations and input when discussing the contents of this paper as it has evolved over the past few years. A.F. would also like to thank Hiroshi Kokubu, Masato Hara, Ryo Terajima, and Yuichiro Terasaki for their influential conversations when discussing the concept of untrained and generated attractors at the recent DEDS conference in Kyoto.
\end{acknowledgments}

\appendix

\section{RC design and training parameters}\label{app:RCdesign}

In this paper, $\textbf{M}$ is constructed with an Erd\"{o}s-Renyi topology where each of the non-zero elements are then replaced with a random number between $-1$ and $1$, the matrix is subsequently scaled to a specific spectral radius, $\rho$. To elaborate, $\textbf{M}$ is designed such that each element is chosen independently to be nonzero with probability $P$ (i.e. sparsity $= P$ or degree $= N/P$) and these nonzero elements are chosen uniformly from $\left( -1, 1 \right)$. This random sparse matrix is rescaled such that the spectral radius, which is the maximum of the absolute values of its eigenvalues, is $\rho$. 
Similarly, the input matrix, $\textbf{W}_{in}$, is designed such that each row has only one nonzero randomly assigned element, chosen uniformly from $\left( -1, 1 \right)$.

The numerical experiments in task (ii) and (iii) are performed using the same $\textbf{M}$ and $\textbf{W}_{in}$. From further analysis
and results which are not shown here, we remark that there are relatively small quantitative changes to our results in task (ii) and (iii) when using different initialisations of $\textbf{M}$ and $\textbf{W}_{in}$ however the main characteristics of our results remain similar.

\section{Classifying output of closed-loop RC in Task (i)}\label{apx:OutputClassification}

The algorithm used to classify the output of the closed-loop RC (Eq.~\eqref{eq:PredRes}) after it was trained to reconstruct the chaotic Lorenz attractor, $\pzl$, according to the steps outlined in Sec.~\ref{sssec:Task1_description} involves the following three main components: The algorithm checks, if after a time $t_{trans}$, 
\begin{enumerate}
    \item[C1:] The output remains stationary (fixed point) or periodic (limit cycle). 
    \item[C2:] The output remains in a subspace of $\mathbb{P}$ that we associate with $\pzl$ and denote by $\mathbb{P}_{\pzl}$, i.e., whether the output remains within a range of $x_{1}$, $x_{2}$, and $x_{3}$ values.
    \item[C3:] The location of the local maxima of the $x_{3}$ variable on both wings of the reconstructed $\pzl$ is less than a distance $\alpha$ away from the lines fitted through the corresponding local maxima on both wings of $\pzl$.
\end{enumerate}
Thus our algorithm has eight parameters, from empirical testing we set these as, $t_{trans} = t_{train} + 70$, $\alpha = 3.75$, $\mathbb{P}_{\pzl} = \left\{\left( x_{1}, x_{2}, x_{3}\right) | -22 \le x_{1} \le 22,\, -32 \le x_{2} \le 32,\, 0 \le x_{3} \le 55 \right\}$.

When C1 is satisfied, we preliminarily classify the output as an `untrained attractor' (UA). By preliminarily classify we mean that we later manually inspect whether this classification is consistent with our definition of an UA as specified in Sec.~\ref{ssec:UA_define}. If C1 is not satisfied, the output is classified as chaotic/quasi-periodic and we then check whether C2 or C3 are satisfied.

We classify the output as a `good reconstruction of $\pzl$', that we denote by $\hat{\pzl}$ if C2 and C3 are satisfied. 
A reconstructed attractor that satisfies this criteria is shown in Fig.~\ref{fig:classification_algorithm_illustration_of_mechanics}~(a). This shows that all local maxima of the reconstructed attractor (blue points) are sufficiently close to the line fitted through the corresponding local maxima values of $\pzl$, indicated by the solid black line.

We classify the output of the closed-loop RC as a `poor reconstruction of $\pzl$', that we denote by $\tilde{\pzl}$ if C2 is satisfied but C3 is not.
A reconstructed attractor which falls into this category is shown in Fig.~\ref{fig:classification_algorithm_illustration_of_mechanics}~(b). This shows that some of the local maxima of the reconstructed attractor (red points) are too far away from the line fitted through the corresponding local maxima values of $\pzl$, indicated by the dashed black line.

If C2 is not satisfied and C3 is, or if C1, C2, or C3 are not satisfied then we also manually inspect the output to decide whether the output is more consistent with $\tilde{\pzl}$ or a UA.



\section{Most common reconstruction route for \texorpdfstring{$\pzl$}{TEXT}\label{apx:LorenzTransientStable}}

\begin{figure*}
    \centering
    \includegraphics[width=0.95\textwidth]{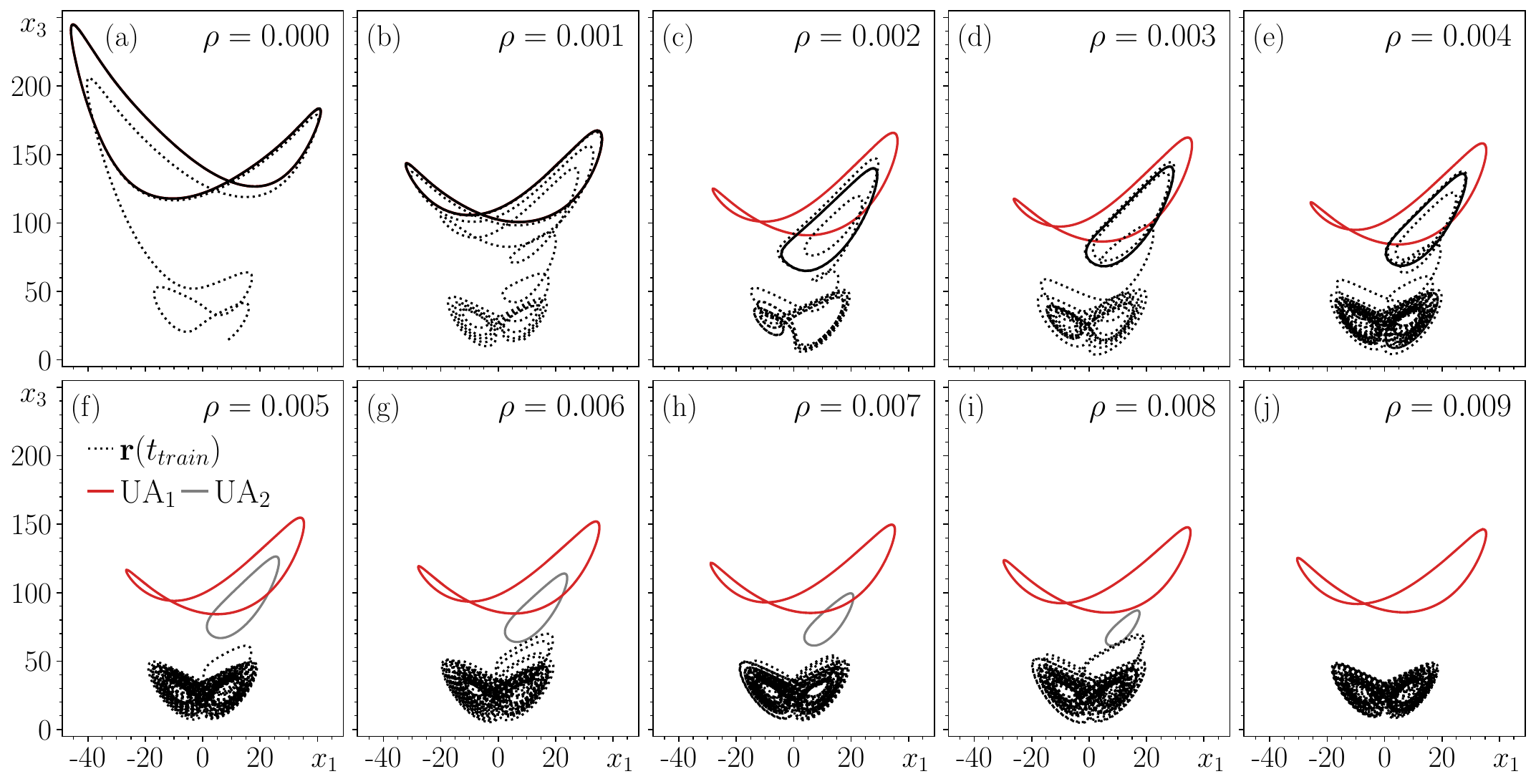}
    \caption{Lorenz transient becoming stable for increasing $\rho$ shown in the $\left( x_{1}, x_{3} \right)$ projection of $\mathbb{P}$. Each black dashed curve corresponds to the trajectory taken by Eq.~\eqref{eq:PredRes} when initialised with Eq.~\eqref{eq:PredResIC}. Curves of different colours correspond to UAs that are present in $\mathbb{P}$ for the specified $\rho$ values.}
    \label{fig:UAs_pred_rho0000_0009}
\end{figure*}

In Fig.~\ref{fig:UAs_pred_rho0000_0009} we show the most common way that the RC provides a reconstruction of $\pzl$ when increasing $\rho$ from $0$ for one of the matrices $\textbf{M}_{i}$ from our ensemble. In each panel of Fig.~\ref{fig:UAs_pred_rho0000_0009}, the black dotted curves correspond to the output of the closed-loop RC (Eq.~\eqref{eq:PredRes}) when initialised with the last point of the training, $\boldsymbol{r}(t_{train})$. The same colour scheme from Fig.~\ref{fig:UA_examples_rho015_} is used to illustrate different UAs that we find. 
In panel (a) we show that while the closed-loop RC is unable to provide any reconstruction of $\pzl$, the brief transient behaviour of the RC's state approaching UA$_{1}$ resembles Lorenz-like behaviour. 
In panels (b)--(e) we see that before approaching UA$_{1}$ or UA$_{2}$, the duration of this transient activity increases and the transient also grows in its resemblance to $\pzl$. 
In panel (f) we no longer see this transient behaviour as the RC is able to provide a poor reconstruction of $\pzl$. However, this reconstruction has some noticeable features on the right-hand wing ($x_{1} > 0$ portion of $\mathbb{P}$) that are not consistent with $\pzl$ which is why our algorithm classifies this attractor as a poor reconstruction of $\pzl$. In panels (g)--(i) we see that as $\rho$ increases, UA$_{2}$ shrinks in size. 
Interestingly, in panel (j), after UA$_{2}$ becomes unstable, this activity on the right-hand wing of the reconstructed $\pzl$ is no longer visible, and this attractor is classified as a good reconstruction of $\pzl$. The RC maintains its good reconstruction of $\pzl$ for larger values of $\rho$. 
Without too much speculation, one reason for this could be that since UA$_{2}$ is located near $\pzl$ in this projection, it may also be in $\mathbb{S}$, and therefore the saddle which separates these states may play a role in preventing the RC from achieving a good reconstruction of $\pzl$, and after UA$_{2}$ becomes unstable this enables the RC to improve its reconstruction.

In summary, Fig.~\ref{fig:UAs_pred_rho0000_0009} shows a sequence of events where we see a transition from an unstable reconstruction of $\pzl$ to a stable reconstruction of $\pzl$ with a stable UA filling in the gaps. From further investigation we find transitions from a stable UA and unstable reconstruction of $\pzl$ to an unstable UA and stable reconstruction of $\pzl$. We also find switching dynamics via transitions from a stable UA and unstable reconstruction of $\pzl$ to an unstable UA and unstable reconstruction of $\pzl$ before another transition to an unstable UA and stable reconstruction of $\pzl$. As shown in \textcite{flynn2024switching}, when two attractors become unstable in quick succession of one another and there are no other attractors present in $\mathbb{P}$, then the state of the RC will switch between transients of both attractors in finite time and remain in either transient for different durations of time.


\section{Additional reconstruction routes for \texorpdfstring{$\pzl$}{TEXT}}\label{apx:AdditionalReconRoutes}

In this Appendix we illustrate some less common routes which we find the RC takes to reconstruct $\pzl$ in task (i). The routes of reconstruction we illustrate here correspond to the classification of the RC's output at different $\rho$ values for $\textbf{M}_{19}$ and $\textbf{M}_{27}$ as shown in Fig.~\ref{fig:rho_vs_M_colorplot_updated}. We use the same method that was used in Fig.~\ref{fig:Lor_PDBifexample_rho_x3_} to track the changes in the dynamics of the reconstructed attractor. For both $\textbf{M}_{19}$ and $\textbf{M}_{27}$ we encounter no UAs for the values of $\rho$ where the tracking is performed. More specifically, we plot $x_{3}^{max}$ (local maxima of the $x_{3}$ variable) and the algorithm's classification at specified values of $\rho$ in Figs.~\ref{fig:Lor_M49_Bifexample_rho_x3_} and \ref{fig:Lor_M15_PDBifexample_rho_x3_}. In Figs.~\ref{fig:Lor_snapshots_M49_} and \ref{fig:Lor_PD_snapshots_M15_} we illustrate the RC's reconstruction of $\pzl$ at key points indicated by the inward ticks in Figs.~\ref{fig:Lor_M49_Bifexample_rho_x3_} and \ref{fig:Lor_M15_PDBifexample_rho_x3_}.

\subsubsection{Reconstruction routes found from \texorpdfstring{$\textbf{M}_{19}$}{TEXT}}

\begin{figure*}
    \centering
    \includegraphics[width=0.8\linewidth]{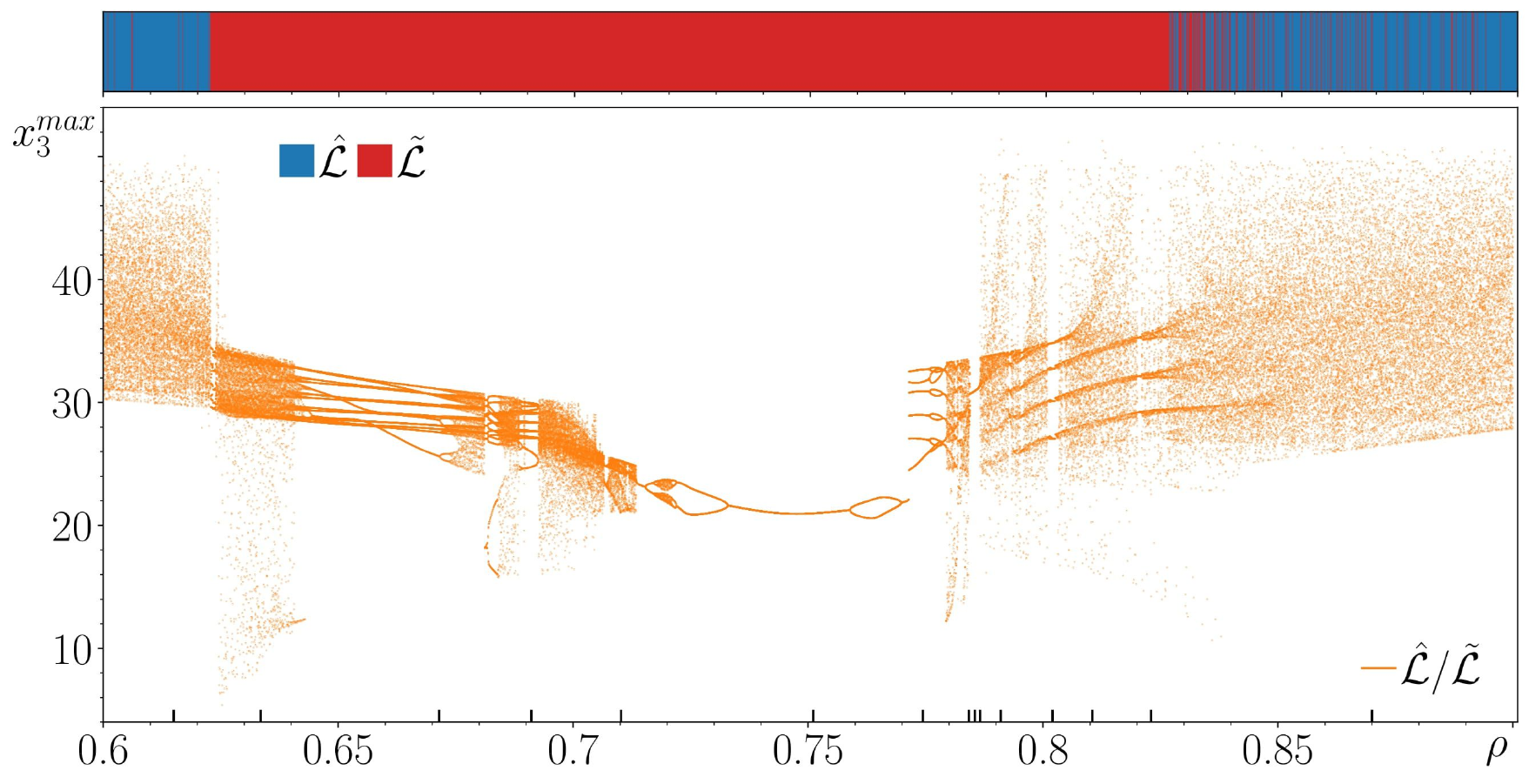}
    \caption{Route for reconstruction of $\pzl$ for $\textbf{M}_{19}$. Plotting changes in local maxima of $x_{3}$ in $\mathbb{P}$ for $x_{1}>0$ for changes in $\rho$ for all attractors found for $\rho \in \left[0.6,0.9\right]$. Inward ticks indicate $\rho$ values used to generate Fig.~\ref{fig:Lor_snapshots_M49_}. Coloured bars above bifurcation diagram indicate the algorithm's classification of the RC's output for a given $\rho$.}
    \label{fig:Lor_M49_Bifexample_rho_x3_}
\end{figure*}

\begin{figure*}
    \centering
    \includegraphics[width=\textwidth]{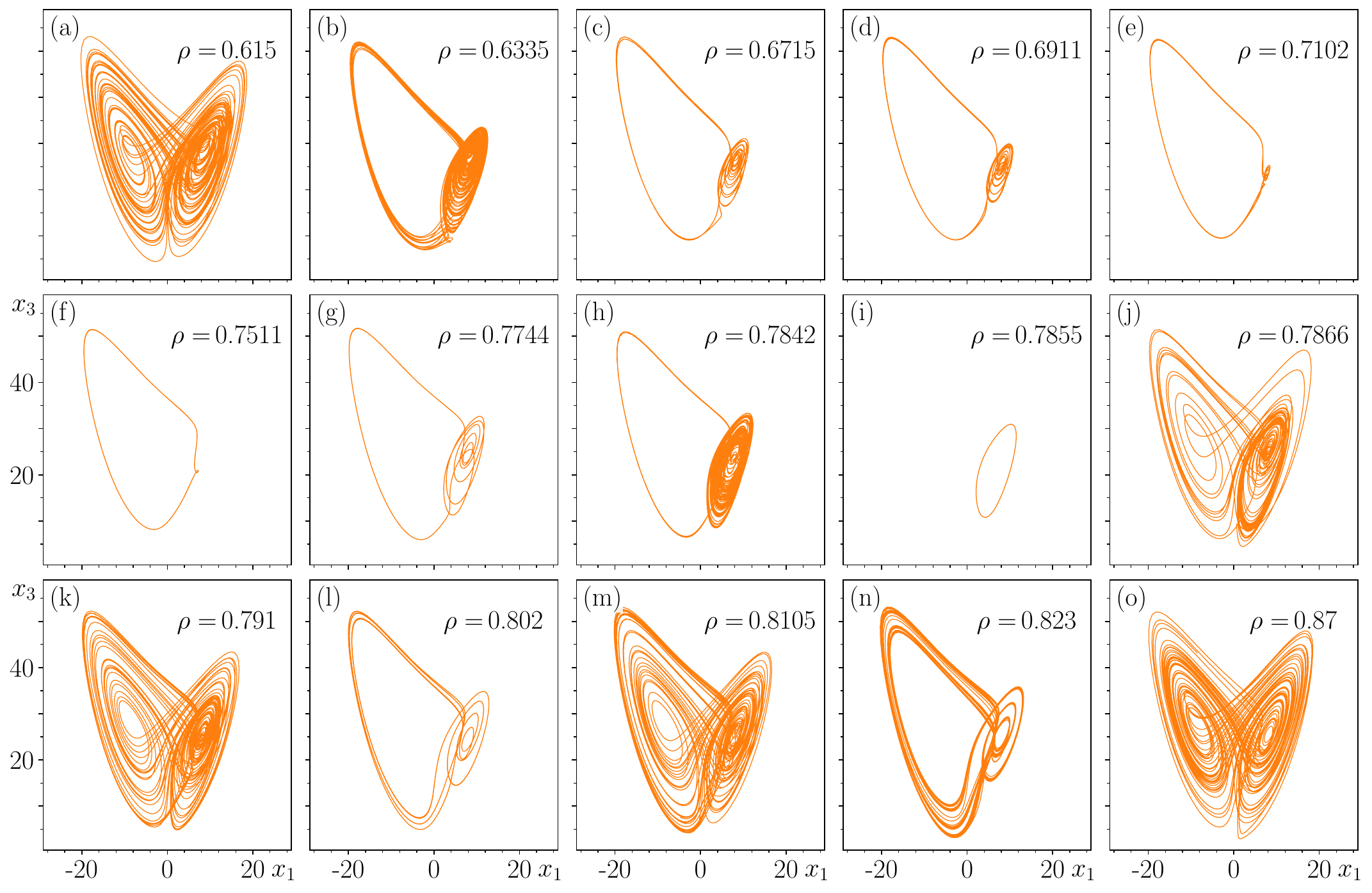}
    \caption{Snapshots of the different routes for reconstruction of $\pzl$ shown in Fig.~\ref{fig:Lor_M49_Bifexample_rho_x3_} at specified values of $\rho$ in panels (a)--(e).}
    \label{fig:Lor_snapshots_M49_}
\end{figure*}

The coloured bars and bifurcation diagram in Fig.~\ref{fig:Lor_M49_Bifexample_rho_x3_} show that while the RC's reconstruction of $\pzl$ is reasonably good for $\rho = 0.6$, the reconstruction worsens before it improves for increasing $\rho$ from $0.6$ to $0.87$. Within this range of $\rho$ values the poor reconstructions of $\pzl$, $\tilde{\pzl}$, undergo several changes between chaotic and periodic dynamics. Fig.~\ref{fig:Lor_PD_snapshots_M15_} highlights some of these changes. 

The first major change in the RC's reconstruction of $\pzl$ occurs nearby $\rho = 0.625$, the transition from $\hat{\pzl}$ to $\tilde{\pzl}$ as shown in Figs.~\ref{fig:Lor_PD_snapshots_M15_}~(a) and (b) involves a break in the symmetry of the reconstructed attractor. However what is more interesting about the dynamics of $\tilde{\pzl}$ for this $\textbf{M}$, and the changes it undergoes for increasing $\rho$ as shown in panels (b) to (i), is that it provides a much clearer indication of the mechanisms that give rise to this interesting break in the symmetry of the RC's dynamics. From looking at these changes in the dynamics of $\tilde{\pzl}$, the prolonged winding motion nearby the point $(x_{1}, x_{3}) = (8,20)$ followed by the much larger orbit strongly resembles a saddle-focus homoclinic orbit, indicating that a Shilnikov bifurcation occurs. Figs.~\ref{fig:Lor_M49_Bifexample_rho_x3_} and \ref{fig:Lor_PD_snapshots_M15_} show that this attractor exhibits periodic dynamics for a wide range of $\rho$ values and also undergoes a period-doubling route to chaos nearby $\rho = 0.675$. While this attractor alternates between chaotic and periodic dynamics, what is most interesting here is that the chaotic dynamics are largely confined to the portion of the trajectory that revolves around the saddle-focus and while that the homoclinic orbit portion of the trajectory can also be chaotic (like in (b)) it is mainly periodic (like in (c)--(h)) depending on the value of $\rho$. In other words, in (b) there are multiple routes in and out of the saddle-focus portion of the trajectory on the attractor but in (c)--(h) there is typically only one or two routes in and out. For a small range of $\rho$ values nearby $0.785$ we see in Fig.~\ref{fig:Lor_M49_Bifexample_rho_x3_}~(i) that there is no homoclinic orbit component to $\tilde{\pzl}$. The growth of the saddle-focus portion of the attractor from (f)--(h) indicates that the transition from (h) to (i) may be due to a collision with a saddle. Increasing $\rho$ further we find that the reconstruction alternates between a more Lorenz-like chaotic attractor and the saddle-focus homoclinic orbit which also alternates between chaotic and periodic dynamics as shown in Figs.~\ref{fig:Lor_M49_Bifexample_rho_x3_}~(j)--(o). For $\rho > 0.825$ we see that reasonably good reconstructions of $\pzl$, like in Fig.~\ref{fig:Lor_M49_Bifexample_rho_x3_}~(o), become more and more prevalent.

\subsubsection{Reconstruction routes found from \texorpdfstring{$\textbf{M}_{27}$}{TEXT}}

\begin{figure*}
    \centering
        \includegraphics[width=0.8\linewidth]{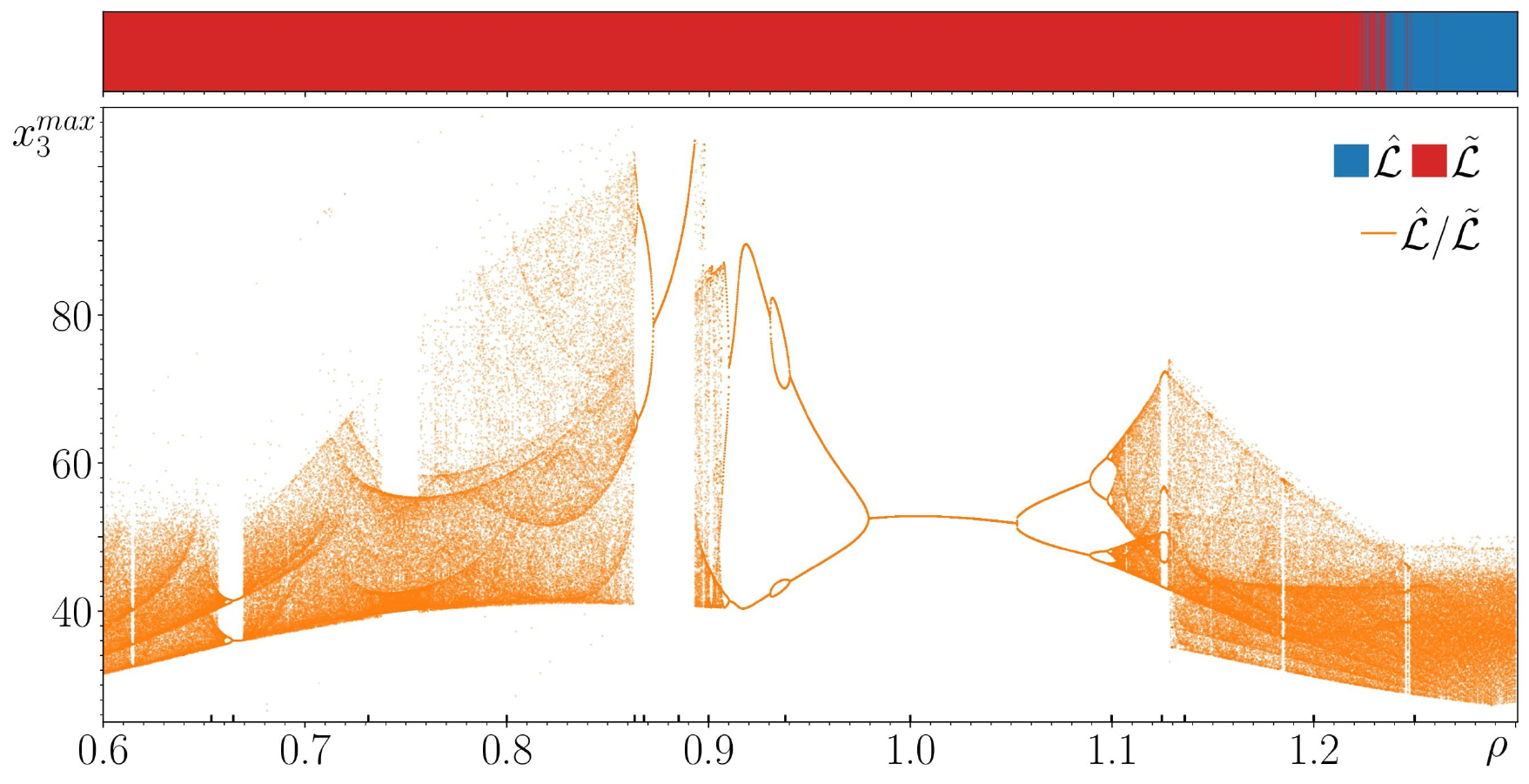}
    \caption{Period-doubling route for reconstruction of $\pzl$ for $\textbf{M}_{27}$. Plotting changes in local maxima of $x_{3}$ in $\mathbb{P}$ for $x_{1}>0$ for changes in $\rho$ for all attractors found for $\rho \in \left[0.6,1.3\right]$. Inward ticks indicate $\rho$ values used to generate Fig.~\ref{fig:Lor_PD_snapshots_M15_}. Coloured bars above bifurcation diagram indicate the algorithm's classification of the RC's output for a given $\rho$.}
    \label{fig:Lor_M15_PDBifexample_rho_x3_}
\end{figure*}

\begin{figure*}
    \centering
    \includegraphics[width=\textwidth]{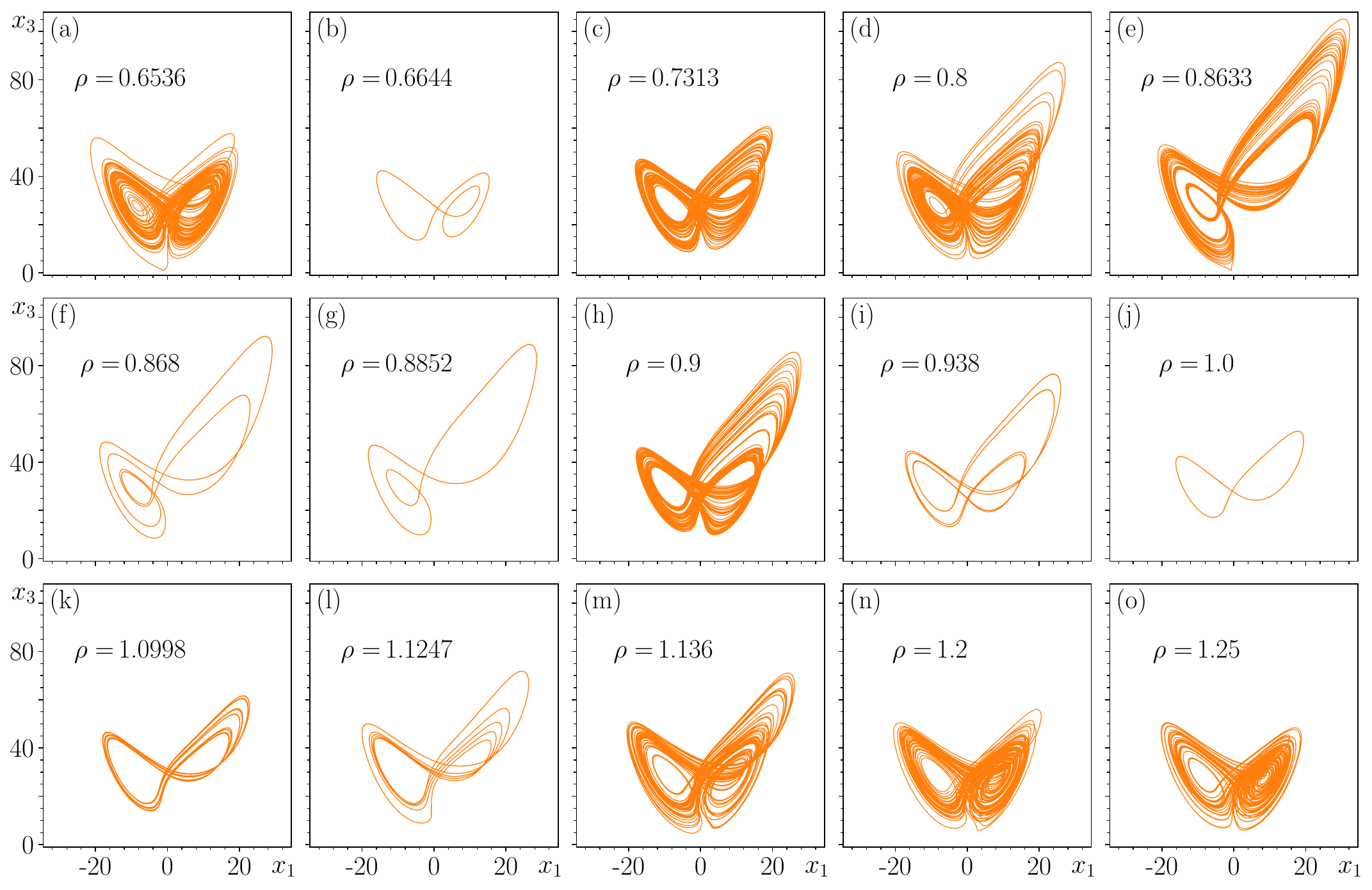}
    \caption{Snapshots of the different routes for reconstruction of $\pzl$ shown in Fig.~\ref{fig:Lor_M15_PDBifexample_rho_x3_} at specified values of $\rho$ in panels (a)--(e). Main types of transitions: homoclinic orbit of saddle focus - period-doubling, interior crisis describes rapid growth of chaotic attractor/limit cycle}
    \label{fig:Lor_PD_snapshots_M15_}
\end{figure*}

The coloured bars and bifurcation diagram in Fig.~\ref{fig:Lor_M49_Bifexample_rho_x3_} show that the RC's reconstruction of $\pzl$ improves by increasing $\rho$, however there is by no means a gradual improvement from the poor reconstruction, $\tilde{\pzl}$, to the good reconstruction, $\hat{\pzl}$, as $\tilde{\pzl}$ undergoes a wide-range of changes from $\rho = 0.6$ to $1.25$.

The first of the major changes in the dynamics of $\tilde{\pzl}$ occurs around $\rho = 0.66$ where the Lorenz-like chaotic attractor in Fig.~\ref{fig:Lor_snapshots_M49_}~(a) becomes a limit cycle in (b). What is interesting about this limit cycle is that there is a clear break in the symmetry of the RC's dynamics, two round trips of the right-hand wing occur before a single trip around the left-hand wing. These periodic dynamics do not persist for many $\rho$ values as $\tilde{\pzl}$ becomes chaotic again near $\rho = 0.67$. 
However, it is after the return to chaos that we see break in the symmetry of the RC's dynamics becomes much more pronounced, the range of values taken by $x_{3}^{max}$ for the right-hand wing of $\tilde{\pzl}$ increases as $\rho$ is increased to $0.8633$. These changes in $\tilde{\pzl}$ are further highlighted by Figs.~\ref{fig:Lor_snapshots_M49_}~(c)--(e), it is also worth noting in (e) that we indirectly observe the presence of a saddle point given the sharp bend in the trajectory nearby $(x_{1}, x_{3}) = (0,0)$ thereby indicating that the attractor is increasing in volume and approaching a saddle point. This may be responsible for the return to periodic dynamics in (f) and (g) that we see as $\rho$ is increased further. In contrast to the limit cycle shown previously in (b), we see in (f) and (g) there are more round trips on the left-hand wing. $\tilde{\pzl}$ alternates between chaotic and periodic dynamics for increasing $\rho$ and the range of values taken by $x_{3}^{max}$ for the right-hand wing of $\tilde{\pzl}$ also starts to decrease. The first instances of good reconstructions of $\pzl$, namely $\hat{\pzl}$, emerge near $\rho = 1.215$ and the RC maintains a reasonably good reconstruction of $\pzl$ for $\rho > 1.24$.

\section{Examining the continuous transition between \texorpdfstring{$\hat{\pzl}$}{TEXT} and \texorpdfstring{$\hat{\pzh}$}{TEXT}}\label{ssec:CloserExam_L_H}

\begin{figure*}[t]
    \centering
    \includegraphics[width=0.73\linewidth]{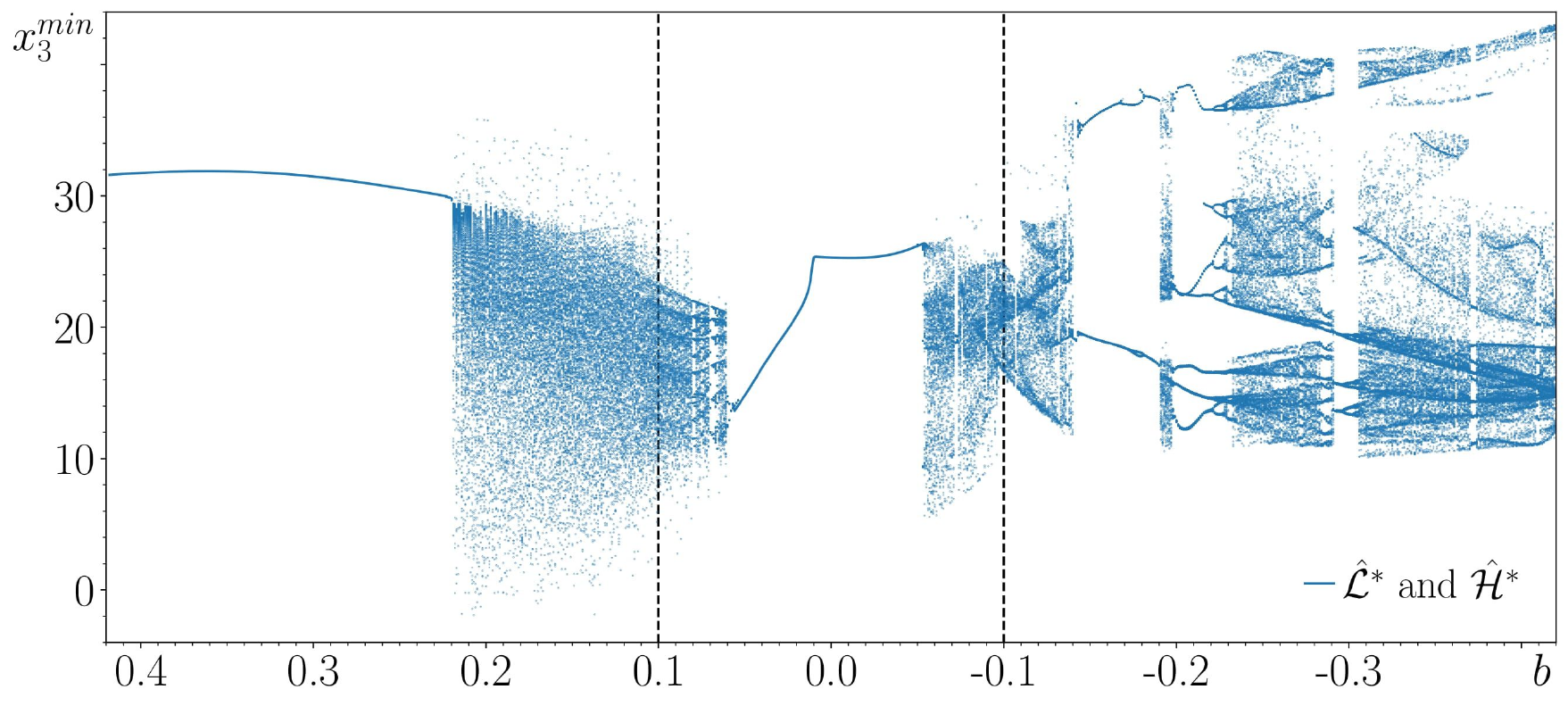}
    \caption{Tracking changes with respect to $b$ in the local minima of the $x_{3}$ variable of the GAs that form after training the parameter-aware RC to perform task (iii) with $b=\pm0.1$.}
    \label{fig:Bif_b01_Lor_Hal_rho12_sigma02_x3}
\end{figure*}

\begin{figure*}
    \centering
    \includegraphics[width=0.85\linewidth]{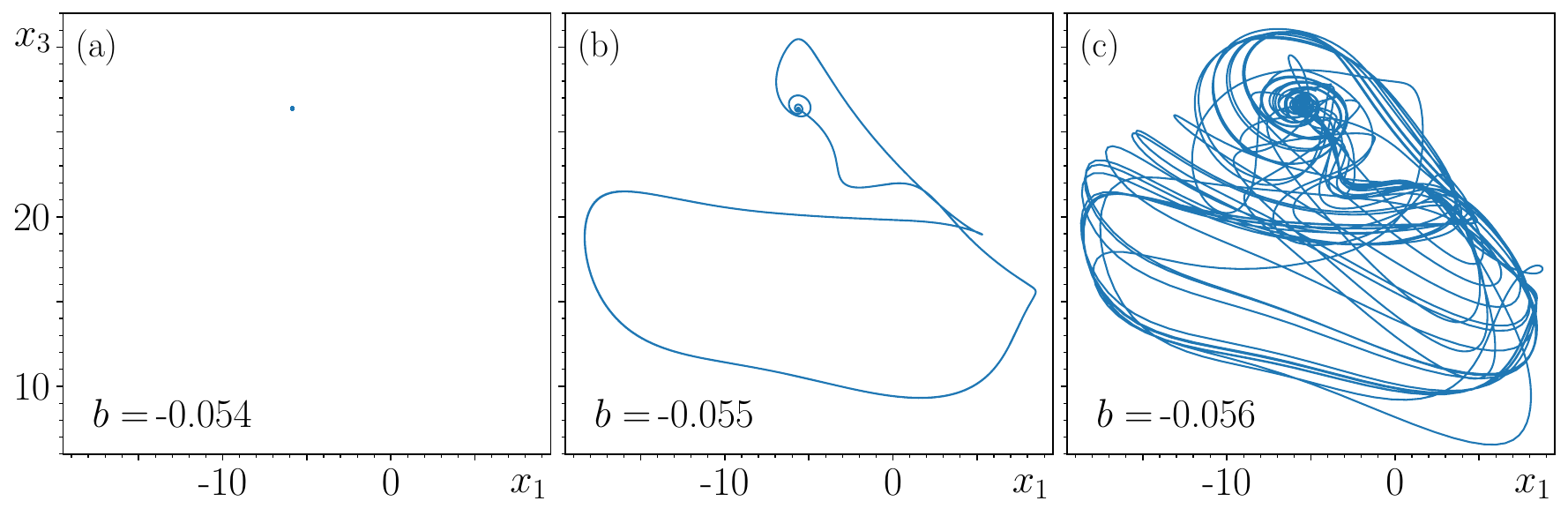}
    \caption{Illustrating transition from $\hat{\pzh}^{*}$ to fixed point for decreasing $b$ as shown in Figs.~\ref{fig:Task3_attrecon_b04_03_02_01_}~(d) and ~\ref{fig:Bif_b01_Lor_Hal_rho12_sigma02_x3}.}
    \label{fig:LorHal_b01_HalToFP_}
\end{figure*}

\begin{figure*}
    \centering
    \includegraphics[width=0.85\linewidth]{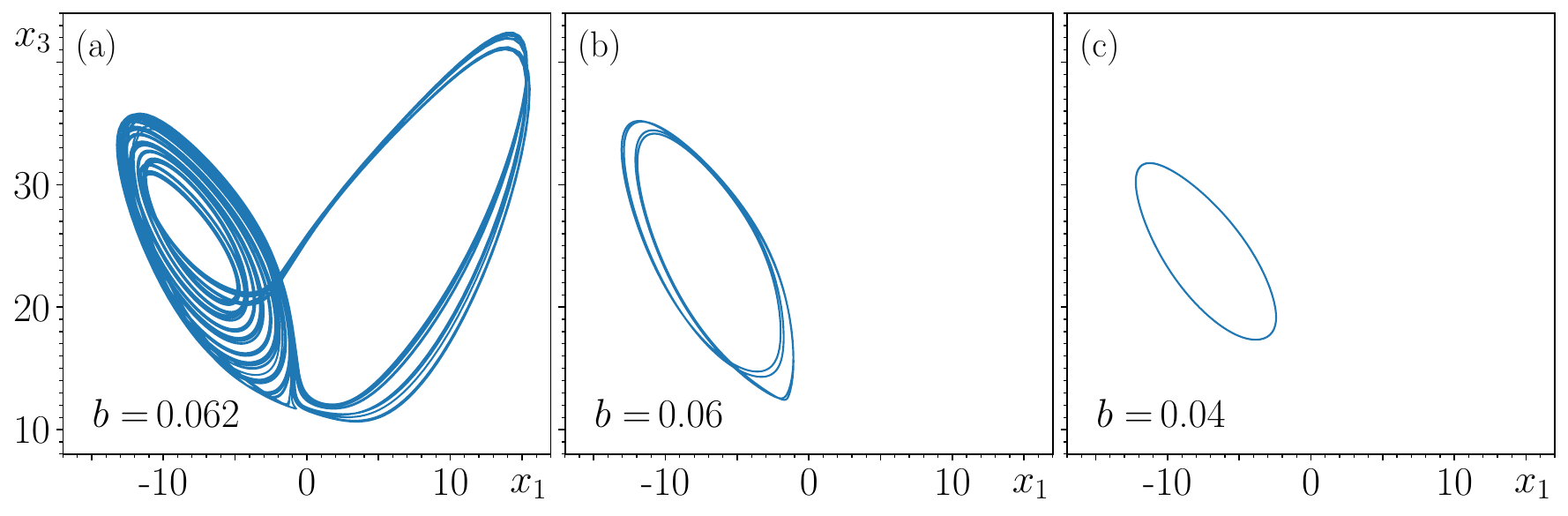}
    \caption{Illustrating transition from $\hat{\pzl}^{*}$ to limit cycle for decreasing $b$ as shown in Figs.~\ref{fig:Task3_attrecon_b04_03_02_01_}~(d) and ~\ref{fig:Bif_b01_Lor_Hal_rho12_sigma02_x3}.}
    \label{fig:LorHal_b01_LorToLC_}
\end{figure*}

To more clearly illustrate that in Fig.~\ref{fig:Task3_bbifplot04030201_compare_}~(d) the transition in $\hat{\pzh}^{*}$ from chaotic to the fixed point dynamics near $b=-0.05$ happens via a bifurcation, in Fig.~\ref{fig:Bif_b01_Lor_Hal_rho12_sigma02_x3} we re-present Fig.~\ref{fig:Task3_bbifplot04030201_compare_}~(d) from the perspective of the changes in the local minima of the $x_{3}$-variable. Also, to help further illustrate this bifurcation and the bifurcation from  the limit cycle near $b=0.055$, in Figs.~\ref{fig:LorHal_b01_HalToFP_} and \ref{fig:LorHal_b01_LorToLC_}, we present snapshots of the RC's dynamics at different $b$ values before and after these bifurcations take place. 


When looking at Figs.~\ref{fig:Bif_b01_Lor_Hal_rho12_sigma02_x3} and \ref{fig:LorHal_b01_HalToFP_} together we see that there is no direct transition from chaotic attractor to fixed point, but rather the chaotic attractor that exists for $b=-0.056$ transitions to a limit cycle for $b=-0.055$. Prior to this transition, we see in Fig.~\ref{fig:LorHal_b01_HalToFP_}~(c) that significant portions of the trajectory on this chaotic attractor are concentrated around where the limit cycle subsequently emerges. What is common among the limit cycle and chaotic attractor is that the trajectory around each attractor appears to contain a point nearby $(x_{1}, x_{3}) = (-6, 26)$ and the state of the system spends a relatively large duration of time in the vicinity of this point. This suggests these attractors are of the saddle-focus homoclinic orbit type. Furthermore, it is this point (the saddle-focus) which becomes the stable fixed point for $b = - 0.054$ shown in Fig.~\ref{fig:LorHal_b01_HalToFP_}~(a). 

In Fig.~\ref{fig:LorHal_b01_LorToLC_} we take a closer look at the transition in $\hat{\pzl}^{*}$ from Lorenz-like dynamics to the limit cycle near $b=0.06$ in Fig.~\ref{fig:Bif_b01_Lor_Hal_rho12_sigma02_x3}. 
In Fig.~\ref{fig:LorHal_b01_LorToLC_}~(a) we see that the state of the RC spends much larger durations of time in the left-hand wing of the attractor. The sharp corner in the trajectory nearby $(x_{1}, x_{3}) = (-1, 12)$ indicates there is a saddle nearby, as flow converges towards this point along the $x_{3}$ direction and diverges away from this point along the $x_{1}$ direction. 
When decreasing $b$ from $0.062$ in Fig.~\ref{fig:LorHal_b01_LorToLC_}~(a) to $0.06$ in Fig.~\ref{fig:LorHal_b01_LorToLC_}~(b), it appears that the Lorenz-like attractor and saddle collide, which results in the right-hand wing disappearing and subsequently results in the period-4 limit cycle seen in Fig.~\ref{fig:LorHal_b01_LorToLC_}~(b). 
We choose to show the change in the dynamics of this limit cycle from $b=0.06$ in Fig.~\ref{fig:LorHal_b01_LorToLC_}~(b) to $b=0.04$ in Fig.~\ref{fig:LorHal_b01_LorToLC_}~(b) in order to compliment the dynamics shown in Fig.~\ref{fig:Bif_b01_Lor_Hal_rho12_sigma02_x3}, that the limit cycle continues to shrink until it becomes point-like and subsequently becomes the fixed point shown explicitly in Fig.~\ref{fig:LorHal_b01_HalToFP_}~(a), which suggests that a Hopf bifurcation is responsible for this transition.

While we do not see regions of bistability between the GAs from the reconstructions of the different limit cycles in task (ii), the results presented in Figs.~\ref{fig:Task3_bbifplot04030201_compare_}~(d) and \ref{fig:Bif_b01_Lor_Hal_rho12_sigma02_x3} indicates that such a scenario is also possible. This leads us to make the following, more general, statement. If the RC fails to interpolate between reconstructed attractors via a sequence of continuous transitions then two options are possible: there is either a region of bistability between GAs that arise from the reconstructed attractors, or there is an UA that coexists with GAs or reconstructed attractors at different parameter settings.

\section*{DATA AVAILABILITY}
The data that support the findings of this study are available from the corresponding author upon reasonable request.





%


\end{document}